\documentclass[preprint,10pt]{elsarticle}
\usepackage{CJKutf8}
\usepackage{algpseudocode}  
\usepackage{subfigure}

\usepackage{algorithmicx}  
\usepackage[colorlinks, linkcolor=blue, anchorcolor=blue, citecolor=blue]{hyperref}
\usepackage[margin=2.5cm]{geometry}
\usepackage{algorithm}
\usepackage{setspace}
\usepackage{lmodern}
\usepackage{amssymb}
\usepackage{amsmath}
\usepackage{amscd,amsxtra,latexsym}
\usepackage{graphicx}
\usepackage{subfigure}
\usepackage{epsf,epsfig}
\usepackage{bm,bbm}
\usepackage{color}
\usepackage[footnotesize,bf]{caption}
\usepackage{mathtools}
\usepackage{hhline}
\usepackage[T1]{fontenc}
\usepackage{placeins}
\usepackage{multirow,dcolumn}
\usepackage{mathrsfs}
\usepackage{verbatim}
\usepackage{comment}
\usepackage{epstopdf}
\usepackage{physics}
\usepackage{pxfonts} 
\usepackage{hyperref}
\usepackage[capitalise]{cleveref}
\usepackage{lineno}

\graphicspath{{figure/}}
\journal{arXiv}

\begin{document}
\begin{CJK}{UTF8}{gbsn}
\begin{frontmatter}



    \title{
    MHDnet: Physics-preserving learning for solving magnetohydrodynamics problems
    }

    \author[TJU,TJU1]{Xiaofei Guan}
    \ead{guanxf@tongji.edu.cn}
    \author[TJU]{Boya Hu}
    \author[LSEC,UCAS]{Shipeng Mao \corref{cor}}
    \ead{maosp@lsec.cc.ac.cn}
    \author[TJU]{Xintong Wang}
    \author[NPU]{Zihao Yang}
    \ead{yangzihao@nwpu.edu.cn}

    \address[TJU]{School of Mathematical Sciences, Tongji University, Shanghai 200092, China}
    \address[TJU1]{Key Laboratory of Advanced Civil Engineering Materials of Ministry of Education,  Tongji University, Shanghai 200092, China}
    \address[LSEC]{LSEC, Academy of Mathematics and Systems Science, Chinese Academy of Sciences, Beijing 100190, China}
    \address[UCAS]{School of Mathematical Sciences, University of Chinese Academy of Sciences, Beijing 100049, China}
    \address[NPU]{School of Mathematics and Statistics, Northwestern Polytechnical University, Xi'an 710129, China}
    \cortext[cor]{Corresponding Author}

    \begin{abstract}
    Designing efficient and high-accuracy numerical methods for complex dynamic incompressible magnetohydrodynamics (MHD) equations remains a challenging problem in various analysis and design tasks. This is mainly due to the nonlinear coupling of the magnetic and velocity fields occurring with convection and Lorentz forces, and multiple physical constraints, which will lead to the limitations of numerical computation. In this paper, we develop the MHDnet as a physics-preserving learning approach to solve MHD problems, where three different mathematical formulations are considered and named $B~formulation$, $A_1~formulation$, and $A_2~formulation$. Then the formulations are embedded into the MHDnet that can preserve the underlying physical properties and divergence-free condition. Moreover, MHDnet is designed 
    by the multi-modes feature merging with multiscale neural network architecture, which can accelerate the convergence of the neural networks (NN) by alleviating the interaction of magnetic fluid coupling across different frequency modes. 
    Furthermore, the pressure fields of three formulations, as the hidden state, can be obtained without extra data and computational cost. Several numerical experiments are presented to demonstrate the performance of the proposed MHDnet compared with different NN architectures and numerical formulations.
    
     
    \end{abstract}


    %


\end{frontmatter}

\section{Introduction}

The MHD system can be formulated by the fully nonlinear coupled equations including Maxwell equations of electromagnetism and the Navier-Stokes equations of fluid dynamics, which describe the macroscopic interaction of conductive fluids under electromagnetic induction. Here, the MHD equations have been considered in terms of two fundamental effects. First, the motion of conducting materials in the presence of a magnetic field induces an electric current and changes the existing electromagnetic field. Second, the current and magnetic fields generate the Lorentz force, which accelerates fluid particles in a direction perpendicular to the magnetic and current fields. In this paper, we are concerned with the dynamic incompressible MHD equations which have received considerable attention for their important applications in liquid metal modeling and plasma physics, such as metallurgical engineering, electromagnetic pumping, the stirring of liquid metals,  liquid metal cooling of nuclear reactors and induction-based flow measurement \cite{burgers_1958, Davidson2001,  Lifschitz1989, Moreau1990,Gerbeau2006MathematicalMF}.

 Let $\Omega$ be a bounded domain in $\mathbb{R}^{d}\ (d=2,3)$ with connected boundary $\partial \Omega$, and the following time-dependent incompressible MHD equations are considered:

\begin{align}\label{MHDeqs1}
    \boldsymbol{u}_{t}-R_{e}^{-1} \Delta \boldsymbol{u}+(\boldsymbol{u} \cdot \nabla) \boldsymbol{u}+\nabla p-\boldsymbol{J} \times \boldsymbol{B}=\boldsymbol{f} & \qquad \text { in } Q_{T},  \\ \label{MHDeqs2}
    \boldsymbol{B}_{t}+\operatorname{curl} \boldsymbol{E}=0                                                                                                       & \qquad \text { in } Q_{T},  \\\label{MHDeqs3}
    \boldsymbol{J}=\kappa(\boldsymbol{E}+\boldsymbol{u} \times \boldsymbol{B})                                                                                    & \qquad \text { in } Q_{T},  \\\label{MHDeqs4}
    \mu^{-1} \operatorname{curl} \boldsymbol{B}=\boldsymbol{J}                                                                                                    & \qquad \text { in } Q_{T},  \\\label{MHDeqs5}
    \operatorname{div} \boldsymbol{u}=0, \quad \operatorname{div} \boldsymbol{B}=0                                                                                & \qquad \text { in } Q_{T},  \\\label{MHDeqs6}
    \boldsymbol{u}(0)=\boldsymbol{u}^{0}, \quad \boldsymbol{B}(0)=\boldsymbol{B}^{0}                                                                              & \qquad \text { in } \Omega,
\end{align}
where $Q_{T}=\Omega \times (0,T)$, $T>0$ is a given finite final time, $\boldsymbol{u}$ denotes the
velocity field, $p$ is the pressure, $\boldsymbol{B}$ is the magnetic induction, $\boldsymbol{J}$ is the electric current
density, $\boldsymbol{E}$ is the electric field, $R_{e}$ is the hydrodynamic Reynolds number, $\kappa$ is the
electric conductivity, $\mu$ is the magnetic permeability, and $\boldsymbol{f}$ is a given forcing term. To make the system of equations well-posed, the no-slip and perfectly conducting wall conditions are defined as
\begin{equation}
    \boldsymbol{u}=0, \qquad \boldsymbol{n} \times \boldsymbol{B}=0\qquad \text { in } \partial\Omega\times(0,T).
    \label{MHDeqs7}
\end{equation}

 In the last several decades, various numerical methods for incompressible MHD problems have been
extensively developed in the literature.  Let us review the references and try to summarize them, which is unlikely to be complete and accurate, of course.
The most commonly used mathematical model for incompressible MHD
system \crefrange{MHDeqs1}{MHDeqs6} is  the so-called $\bm{B}$ based MHD equations which are described by \cref{eq:mhdusual}  in the next section. There is a lot of research on numerical methods for this model, especially represented by finite element methods.  Interested readers are referred to \cite{Gunzburger1991, Armero1996, Gerbeau2000,  Wiedmer2000,  He2015,  Su2019}   for the discretization by using standard Lagrange $\bm{H}^1$ finite element spaces to approximate both the hydrodynamic unknowns and the magnetic induction. 
 However, it is known that the nodal finite element method discretizations of the magnetic operator cannot be correctly approximated when the magnetic induction components may have regularity below $\bm{H}^1(\Omega)$, cf. \cite{Costabel2002,Hiptmair2002}, which may be frequently encountered in non-convex polyhedral or with a non $C^{1,1}$ boundary.  A possible way to overcome these difficulties is by virtue of  N$\acute{\mbox{e}}$d$\acute{\mbox{e}}$lec finite elements for the magnetic induction $\bm{B}$, which leads to another natural formulation and is valid for non-smooth magnetic solution, see \cite{Schoetzau2004, Prohl2008, Gao2019, Ding2020ConvergenceAO,ding2022convergence}  and the references therein.
In addition, the preserving physical properties of numerical schemes for the original MHD system should be carefully considered. Numerous studies have been extensively developed to overcome these difficulties, such as the special treatments of constraint $ \operatorname{div} \boldsymbol{B}=0 $ \cite{article}, i.e., the preservation of the Gauss’s law of discrete magnetic induction. This condition means that there is no source of the magnetic field and it guarantees that no magnetic monopole exists.  The small perturbations to this condition will introduce a strong non-physical force and may cause huge errors in the computation of MHD, see e.g., \cite{BB:JCP:1980,Toth2000}.  This divergence-free condition also plays an important role in the development of a simulation code for fusion reactor blanket \cite{Abdou2001,Abdou2005}.  In \cite{Hu2017},  the authors proposed magnetic-electric formulations for incompressible MHD, in which the electrical field $\bm{E}$ was approximated by $\mathrm{N\acute{e}d\acute{e}lec}$ edge elements and magnetic induction $\bm{B}$ was approximated by Raviart-Thomas face elements. In this way, it naturally preserves the important Gauss's law. Recently, the third author of this paper and his coauthors proposed and studied a fully finite element scheme for the incompressible MHD system in \cite{Hiptmair2018, Ding2021} by introducing a new vector potential formulation. Concerning the preserving of mass conservation, i.e., $\operatorname{div} \boldsymbol{u}=0 $ for the incompressible flows, the readers are referred to \cite{Greif2010,Ding2022} and the references therein.
The above methods are all based on traditional numerical methods 
and rely on specially designed schemes and high-quality generated meshes to obtain approximate solutions, which can be computationally expensive and may not be suitable for complex PDEs or high-dimensional problems.

It is worth noting that, compared with traditional discrete methods, deep learning-based methods have great potential in solving partial differential equations (PDEs). Neural networks can directly approximate the solution of the PDEs without discretization, and they can learn the underlying physics of the system from the given data and generalize to unseen inputs. By combining the physics-based knowledge of the PDEs with data-driven neural network models, it is possible to solve PDEs more efficiently and accurately. The methods can mainly be divided into Deep Ritz method \cite{EY:CMS:2018}, Deep Galerkin method \cite{SS:JCP:2018}, and Physics-informed neural network (PINN) \cite{RPK:JCP:2019}, etc.
\cite{sukumarExactImpositionBoundary2022} introduced a new geometric perception method with R-functions and transfinite interpolation, which can be applied to boundary conditions on complex (affine, curved, and multi-connected) geometries. \cite{zangWeakAdversarialNetworks2020} reformulated the weak solution and test function of PDEs into the original and adversarial parts of neural networks, respectively. During training, the parameters of these two subnetworks are alternately updated to achieve optimal results, which is efficient in high-dimensional problems.
\cite{jinNSFnetsNavierStokesFlow2021a} developed the NSFnets by considering two different mathematical formulations of the Navier-Stokes equations, where pressure fields are obtained as a hidden state without additional data and computational cost.
In summary, neural network methods for solving PDEs have many advantages such as efficiency \cite{MN:2021,GWW:2023}, nonlinear modeling \cite{GZR:u:2019,YZCW:2020}, scalability \cite {RYK:S:2020,YWGH:2022}, and adaptability \cite{ZLGK:JCP:2019}.
However, when dealing with complex coupling PDEs, especially for the underlying problems with high-frequency modes, these deep learning methods still face challenges. 
 The fitting dynamics of deep neural networks show characteristics of frequency principles and spectral biases \cite{jacot_neural_2020}. Specifically, the neural network will first fit the low-frequency components of the target function with gradient-based training algorithms, and then fit the high-frequency components. This may ultimately lead to failure in training optimization.
To address these issues, \cite{FrequencyPrinciple} proposed a multiscale neural network method, and Fourier feature embeddings \cite{liDeepDomainDecomposition2022a,tancikFourierFeaturesLet2020} can also be used to design effective neural tangent kernel (NTK) feature spaces. 

The goal of this paper is to develop a novel physics-preserving learning framework, named MHDnet, for solving MHD problems with three different mathematical formulations. One is $B~formulation$ (see \cref{eq:mhdusual} in \cref{sec:cons}), another two novel formulations are $A_1~formulation$ (see \cref{eq:mao1} in \cref{sec:cons}), and $A_2~formulation$ (see \cref{eq:mao22} in \cref{sec:cons})
where the $\operatorname{div}\boldsymbol{B}=0$, $\operatorname{div}\boldsymbol{u}=0 $ are naturally satisfied for $A_2$ formulation.
Due to the complexity of the MHD system, the different loss terms from different physical laws are naturally tied to measurement units and have significant differences in magnitude. Therefore, the interaction between different loss terms will lead to model training failure, which is related to numerical stiffness with unbalanced back-propagated gradients. Although the learning rate annealing algorithm \cite{Gradientflow2021} has successfully applied to the PINN with complex loss functions, it is still difficult and inefficient to solve the coupling MHD equations. Thus, by combining multiscale neural networks \cite{CiCP-28-1970} with novel multi-modes features \cite{ZYG:AAMM:2023,YHFG:SJSC:2022}, MHDnet can efficiently balance the coupling loss terms of MHD equations with the high-accuracy approximation. Thus, this method can also be robust for MHD problems with missing or noisy data, and the pressure fields as the hidden state can be obtained by $A_1$ and $A_2$ formulations of MHDnet with high accuracy. To our knowledge, MHDnet is the first attempt to solve three-dimensional MHD equations under the deep learning framework, and it can be easily extended to more complex coupled physical problems.

The rest of the paper is organized as follows. In Section \cref{sec:pinn}, a short analysis of the PINN for the MHD system is given. Then three formulations of the MHD equations are derived for MHDnet, which combines the multiscale NN architecture and the multi-modes feature. Moreover, a detailed implementation algorithm is given for the proposed MHDnet in \cref{sec:cons}. In Section \cref{sec:num}, several numerical experiments are presented to demonstrate the performance of the proposed method in 2D and 3D cases. Finally, some conclusions and remarks are given in Section \cref{sec:con}.

\section{Physics-informed neural networks (PINNs) for MHD equations}\label{sec:pinn}
In this subsection, we first briefly review the main idea of the PINN approach for solving MHD equations \crefrange{MHDeqs1}{MHDeqs7}. We consider the general operator forms of MHD equations as follows:
\begin{align}\label{MHDeqs8}
    \mathcal{L}[\boldsymbol{u(x},t),\boldsymbol{B(x},t)] & = \boldsymbol{\hat{f}},\qquad \text { in } Q_{T},                     \\\label{MHDeqs9}
    \mathcal{I}[\boldsymbol{u(x},0),\boldsymbol{B(x},0)] & = \boldsymbol{\hat{g}},\qquad \text { in } \Omega,                    \\\label{MHDeqs10}
    \mathcal{B}[\boldsymbol{u(x},t),\boldsymbol{B(x},t)] & = \boldsymbol{\hat{h}},\qquad \text { in } \partial\Omega\times(0,T),
\end{align}
where $\mathcal{L}$ denotes the general differential operator, and $\boldsymbol{\hat{f}}$ is the source or force vector. $\mathcal{I}$ represents the initial conditions operator defined at the initial time $t=0$, and $\boldsymbol{\hat{g}}$ is the initial value. $\mathcal{B}$ is the boundary operator, which can be an identity operator for the Dirichlet boundary condition, a differential operator for the Neumann boundary condition, or a mixed identity-differential operator for the Robin boundary condition, and the desired values of $\boldsymbol{\hat{h}}$ is defined at the domain boundary $\partial \Omega$.

\begin{figure}[ht]
    \centering
    \begin{minipage}[t]{0.8\textwidth}%
        \centering
        \includegraphics[width=\textwidth]{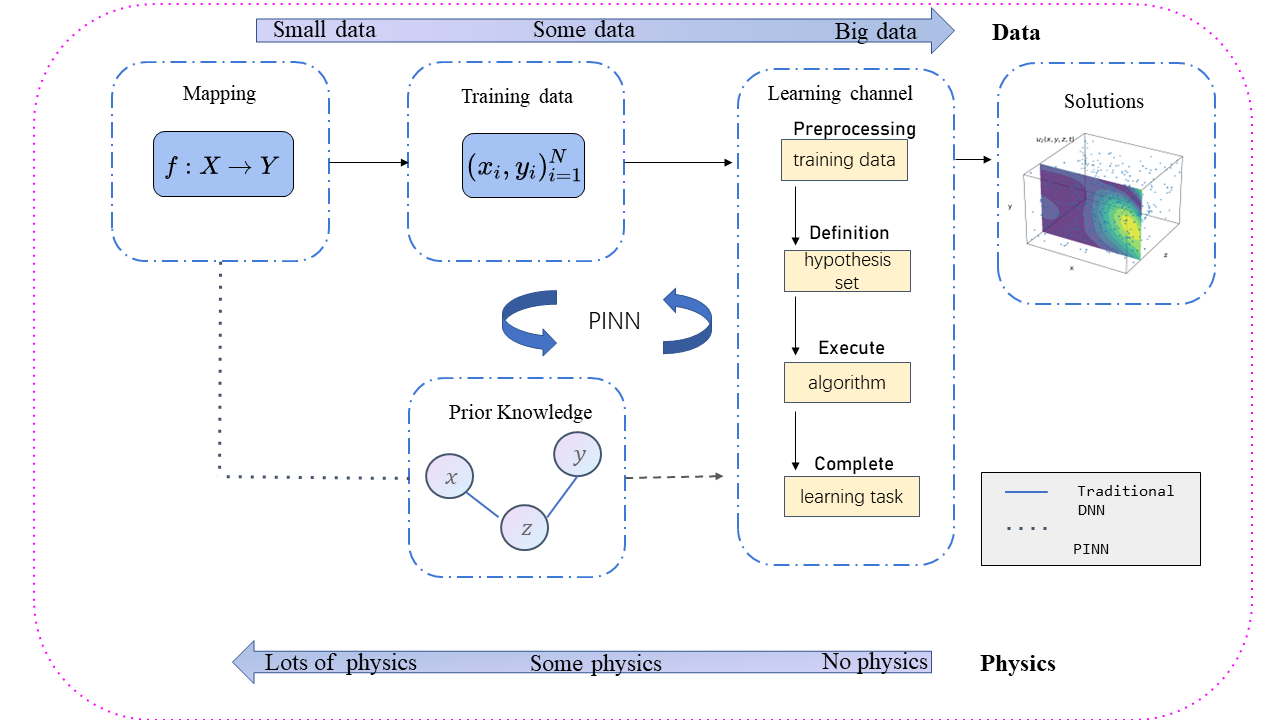}
        \caption{PINN serves as a bridge for fusing data-driven methods with knowledge graphs and knowledge reasoning. }%
        \label{fig1:pipline}
    \end{minipage}
\end{figure}

MHDnet is designed by use of the PINN framework. The difference and connections between PINN and traditional data-driven DNN are shown in \cref{fig1:pipline}. Generally speaking, machine learning (ML) aims to learn a specific mapping relationship from inputs X to outputs Y, where the known data information or prior knowledge is fed into the learning pipeline to complete the regression task. The learning pipeline mainly consists of four parts: training samples, prior knowledge, ML algorithms, and regression tasks. However, the traditional ML methods generally only use sample data as known information, but PINN-based learning pipelines can fuse different information sources, including sample data and prior knowledge. Here, prior knowledge exists independently of the learning task and includes some physical laws \cite{JKK:CMAME:2020}, knowledge graphs \cite{GZW:CMAME:2022}, empirical information \cite{wangRespectingCausalityAll2022}, etc., and it is usually provided and expressed as explicit equations and boundary conditions. Therefore, to a certain extent, PINN can avoid the waste of sample data with the additional constraint. In contrast to the data-driven approach, PINN is not a purely 'knowledge-driven' approach but can be seen as a bridge between the data-driven and traditional knowledge-driven approaches.

Then, the solutions $\boldsymbol{u}(\boldsymbol{x},t),\boldsymbol{B}(\boldsymbol{x},t)$ of the MHD equations can be approximated by the PINN surrogate $\boldsymbol{\tilde{u}}(\boldsymbol{x},t;\theta),\boldsymbol{\tilde{B}}(\boldsymbol{x},t;\theta)$ with DNN parameters $\theta$,
and the $\boldsymbol{\tilde{u}}(\boldsymbol{x},t;\theta),\boldsymbol{\tilde{B}}(\boldsymbol{x},t;\theta)$ can be obtained
by MHD equations \crefrange{MHDeqs8}{MHDeqs10} integrated with automatic differentiation and some machine learning packages, such as PyTorch or TensorFlow, etc. Assume that we have the training data set $D = (D_{\boldsymbol{\hat{f}}}, D_{\boldsymbol{\hat{g}}}, D_{\boldsymbol{\hat{h}}})$ defined as follows
\begin{equation}
    D_{\boldsymbol{\hat{f}}}=\big\{\boldsymbol{x}_i,t_i,\boldsymbol{\hat{f}}_i \big\}_{i=1}^{N_{\boldsymbol{\hat{f}}}},\quad
    D_{\boldsymbol{\hat{g}}}=\big\{\boldsymbol{x}_i,t_i,\boldsymbol{\hat{g}}_i \big\}_{i=1}^{N_{\boldsymbol{\hat{g}}}},\quad
    D_{\boldsymbol{\hat{h}}}=\big\{\boldsymbol{x}_i,t_i,\boldsymbol{\hat{h}}_i \big\}_{i=1}^{N_{\boldsymbol{\hat{h}}}},
\end{equation}
where the coordinates $\boldsymbol{x},t$ are the input of $\boldsymbol{\tilde{u}}, \boldsymbol{\tilde{B}}$, which have the same dimension with $\boldsymbol{u},\boldsymbol{B}$. Furthermore, the DNN parameters $\theta$ will be trained to reach the optimal set $\hat{\theta}$, by using the gradient descent algorithm to minimize the following loss function $Loss(\theta)$:

\begin{equation}\label{eq:loss1}
    \displaystyle Loss(\theta)=\lambda_{\boldsymbol{\hat{f}}} L_{\boldsymbol{\hat{f}}}+\lambda_{\boldsymbol{\hat{g}}} L_{\boldsymbol{\hat{g}}}+\lambda_{\boldsymbol{\hat{h}}} L_{\boldsymbol{\hat{h}}}+\lambda_{data} L_{data}.
\end{equation}
where
\begin{equation}
    \begin{array}{lll}
        \displaystyle L_{\boldsymbol{\hat{f}}} =  \frac{1}{N_{\boldsymbol{\hat{f}}}} \left\| \boldsymbol{\tilde{f}}(\boldsymbol{x}_i,t_i;\theta)-\boldsymbol{\hat{f}}_i \right\|_{Q_T}^{2},\quad
        L_{\boldsymbol{\hat{g}}} =  \frac{1}{N_{\boldsymbol{\hat{g}}}} \left\| \boldsymbol{\tilde{g}}(\boldsymbol{x}_i,0;\theta)-\boldsymbol{\hat{g}}_i \right\|_{\Omega}^{2},\quad \vspace{1ex} \\
        \displaystyle L_{\boldsymbol{\hat{h}}} =  \frac{1}{N_{\boldsymbol{\hat{h}}}} \left\| \boldsymbol{\tilde{h}}(\boldsymbol{x}_i,t_i;\theta)-\boldsymbol{\hat{h}}_i \right\|_{\partial \Omega \times(0, T]}^{2},\quad
        L_{data}=\frac{1}{N_{data}} \left\{\left[\boldsymbol{\tilde{u}}(\boldsymbol{x}_i,t_i;\theta)-\boldsymbol{u}_i\right]^2+\left[\boldsymbol{\tilde{B}}(\boldsymbol{x}_i,t_i;\theta)-\boldsymbol{B}_i\right]^2\right\},
    \end{array}
\end{equation}
and
\begin{equation}
    \begin{array}{lll}
        \boldsymbol{\tilde{f}}=\mathcal{L}_{\theta}[\boldsymbol{\tilde{u}}, \boldsymbol{\tilde{B}}],\quad
        \boldsymbol{\tilde{g}}=\mathcal{I}_{\theta}[\boldsymbol{\tilde{u}}, \boldsymbol{\tilde{B}}],\quad
        \boldsymbol{\tilde{h}}=\mathcal{B}_{\theta}[\boldsymbol{\tilde{u}}, \boldsymbol{\tilde{B}}].\quad
    \end{array}
\end{equation}

Notice that $L_{data}$ denotes the loss term for the observable data, this is because some physical phenomena can not be directly formulated by the analytical equations, and some observable data should be added. $\lambda_{\boldsymbol{\hat{f}}}$, $\lambda_{\boldsymbol{\hat{g}}}$, $\lambda_{\boldsymbol{\hat{h}}}$ and $\lambda_{data}$ represent the weight factor of different terms in the loss function, respectively. Since each term in the loss function is derived from the laws of physics, they are naturally bound to measurement units varied greatly in magnitude. Suitable scaling and trad-off of weight factors and loss functions can significantly speed up the convergence of PINN training \cite{multiobject}.

\section{Construction of MHDnet} \label{sec:cons}
\subsection{Three different mathematical formulations of MHD equations}
In this subsection, three different mathematical formulations of MHD equations will be introduced for the construction of MHDnet. Firstly, the most commonly used mathematical models of the incompressible MHD system \crefrange{MHDeqs1}{MHDeqs6} are considered by inserting \cref{MHDeqs3} and \cref{MHDeqs4} in \cref{MHDeqs2}, in which the electric field $\boldsymbol{E}$ and the current density $\boldsymbol{J}$ are eliminated. Then, the MHD system can be described by the following equations with the fields $(\boldsymbol{u}, p, \boldsymbol{B})$:
\begin{equation}
    \begin{array}{rc}
        \boldsymbol{u}_{t}-R_{e}^{-1} \Delta \boldsymbol{u}+(\boldsymbol{u} \cdot \nabla) \boldsymbol{u}+\nabla p-\mu^{-1} \operatorname{curl} \boldsymbol{B} \times \boldsymbol{B}=\boldsymbol{f} & \text { in } Q_{T}, \\
        \boldsymbol{B}_{t}+(\mu \kappa)^{-1}\operatorname{curl} \operatorname{curl} \boldsymbol{B}-\operatorname{curl}(\boldsymbol{u} \times \boldsymbol{B})=0                                     & \text { in } Q_{T}, \\
        \operatorname{div} \boldsymbol{u}=0, \quad \operatorname{div} \boldsymbol{B}=0                                                                                                             & \text { in } Q_{T}, \\
    \end{array}
    \label{eq:mhdusual}
\end{equation}
and \cref{eq:mhdusual} is named $B~formulation$.

Then, considering the preservation of physical conservation, two novel mathematical formulations $\boldsymbol{A_1}$ and $\boldsymbol{A_2}$ for the incompressible MHD system will be given. In \cite{Hiptmair2018, Ding2021}, a new vector potential formulation is first proposed, where the divergence-free property of $\boldsymbol{B}$ can be replaced by means of the magnetic vector potential $\boldsymbol{A_1}$ with a gauge such that
\begin{equation}\label{eq:Adefin}
    \boldsymbol{B}=\operatorname{curl} \boldsymbol{A}_1,\qquad \boldsymbol{E}=-\boldsymbol{A}_{1_t}.
\end{equation}
Inserting \cref{eq:Adefin} into \crefrange{MHDeqs1}{MHDeqs6} with a simple calculation, the incompressible MHD problem can be given as follows
\begin{equation}
    \begin{array}{rc}\vspace{2mm}
        \boldsymbol{u}_{t}-R_{e}^{-1} \Delta \boldsymbol{u}+(\boldsymbol{u} \cdot \nabla) \boldsymbol{u}+\nabla p-\boldsymbol{J}\times  \operatorname{curl} \boldsymbol{A}_1 =\boldsymbol{f} & \text { in } Q_{T}, \\\vspace{2mm}
        \boldsymbol{A}_{1_t}+\operatorname{curl} \boldsymbol{A}_1\times\boldsymbol{u} +R_{m}^{-1}\operatorname{curl}\operatorname{curl}\boldsymbol{A}_1 =0                                   & \text { in } Q_{T}, \\
        \operatorname{div} \boldsymbol{u}=0, \quad \boldsymbol{u}(0)=\boldsymbol{u}^{0}, \quad \boldsymbol{A}_1(0)=\boldsymbol{A}^{0}                                                        & \text { in } Q_{T}, \\
    \end{array}
    \label{eq:mao}
\end{equation}
where $R_{m}$ is the magnetic Reynolds number, and the electric current density satisfies
\begin{equation}\label{eq:Jdefin}
    \boldsymbol{J}=\kappa(\boldsymbol{E}+\boldsymbol{u}\times\boldsymbol{B})=-\kappa(\boldsymbol{A}_{1_t}+\operatorname{curl} \boldsymbol{A}_1\times\boldsymbol{u}).
\end{equation}
Substituting \cref{eq:Jdefin} into \cref{eq:mao}, there holds
\begin{equation}
    \begin{array}{rc}\vspace{2mm}
        \boldsymbol{u}_{t}-R_{e}^{-1} \Delta \boldsymbol{u}+(\boldsymbol{u} \cdot \nabla) \boldsymbol{u}+\nabla p+\kappa(\boldsymbol{A}_{1,t}+\operatorname{curl} \boldsymbol{A}_1 \times\boldsymbol{u})\times  \operatorname{curl} \boldsymbol{A}_1 =\boldsymbol{f} & \text { in } Q_{T}, \\\vspace{2mm}
        \boldsymbol{A}_{1_t}+\operatorname{curl} \boldsymbol{A}_1\times\boldsymbol{u} +R_{m}^{-1}\operatorname{curl}\operatorname{curl}\boldsymbol{A}_1 =0                                                                                      & \text { in } Q_{T}, \\
        \operatorname{div} \boldsymbol{u}=0, \quad \boldsymbol{u}(0)=\boldsymbol{u}^{0}, \quad \boldsymbol{A}_1(0)=\boldsymbol{A}_1^{0}                                                                                                         & \text { in } Q_{T}, \\
    \end{array}
    \label{eq:mao1}
\end{equation}
where \cref{eq:mao1} is named $A_1~formulation$, which keeps the energy stability (see Theorem 3 in \cite{Ding2021}).


Furthermore, considering $\operatorname{div}\boldsymbol{u}=0 $ of the MHD system, and for any vector field $\boldsymbol{u}$, we define
\begin{equation}
    \boldsymbol{u}=\operatorname{curl} \boldsymbol{A}_2,
    \label{eq:A2-1}
\end{equation}
then we have
\begin{equation}
    \operatorname{div} \boldsymbol{u}=\boldsymbol{0}.
    \label{eq:A2-2}
\end{equation}
Substituting \cref{eq:A2-1} and \cref{eq:A2-2} into \cref{eq:mao1}, the novel $A_2~formulation$ can be defined as follows 
\begin{equation}
    \begin{array}{rc}\vspace{2mm}
    \operatorname{curl}\boldsymbol{A}_{2_t}-R_{e}^{-1} \Delta \operatorname{curl} \boldsymbol{A}_2+(\operatorname{curl} \boldsymbol{A}_2 \cdot \nabla) \operatorname{curl} \boldsymbol{A}_2+\nabla p+\kappa(\boldsymbol{A}_{1_t}+\operatorname{curl} \boldsymbol{A}_1  \times \operatorname{curl} \boldsymbol{A}_2 )\times  \operatorname{curl} \boldsymbol{A}_1 =\boldsymbol{f}  &\text { in } Q_{T},  \\
    \boldsymbol{A}_{1_t}+\operatorname{curl} \boldsymbol{A}_1\times\operatorname{curl} \boldsymbol{A}_2 +R_{m}^{-1}\operatorname{curl}\operatorname{curl}\boldsymbol{A}_1 =0                                                                                                                                                                              &\text { in } Q_{T}, \\
    \quad \boldsymbol{A}_2(0)=\boldsymbol{A}_2^{0}, \quad \boldsymbol{A}_1(0) 
 =\boldsymbol{A}_1^{0}            &\text { in } Q_{T}. 
    \end{array}
    \label{eq:mao22}
\end{equation}
It can both keep the properties of incompressibility and energy stability, which can alleviate the unbalanced gradient. Then, the MHD system can be solved with the 
 following homogeneous Dirichlet boundary conditions:
\begin{equation}
    \boldsymbol{A}_2=0,\qquad \boldsymbol{A}_1\times\boldsymbol{n}=0,\quad \text { on } \partial\Omega\times(0,T).
    \label{eq:maobd1}
\end{equation}

Finally, the loss terms in \cref{eq:loss1} can be rewritten, and $L_{\boldsymbol{\hat{f}}}$ represents the residual of the first two equations of  \cref{eq:mhdusual}, \cref{eq:mao1}, and \cref{eq:mao22}, respectively. Similarly, $L_{\boldsymbol{\hat{g}}}$ is the corresponding initial conditions, and  $L_{\boldsymbol{\hat{h}}}$ denotes \cref{eq:maobd1}. 

\begin{figure}[htbp]
    \centering
    \begin{minipage}[t]{0.85\textwidth}
        \centering
        \includegraphics[width=\textwidth]{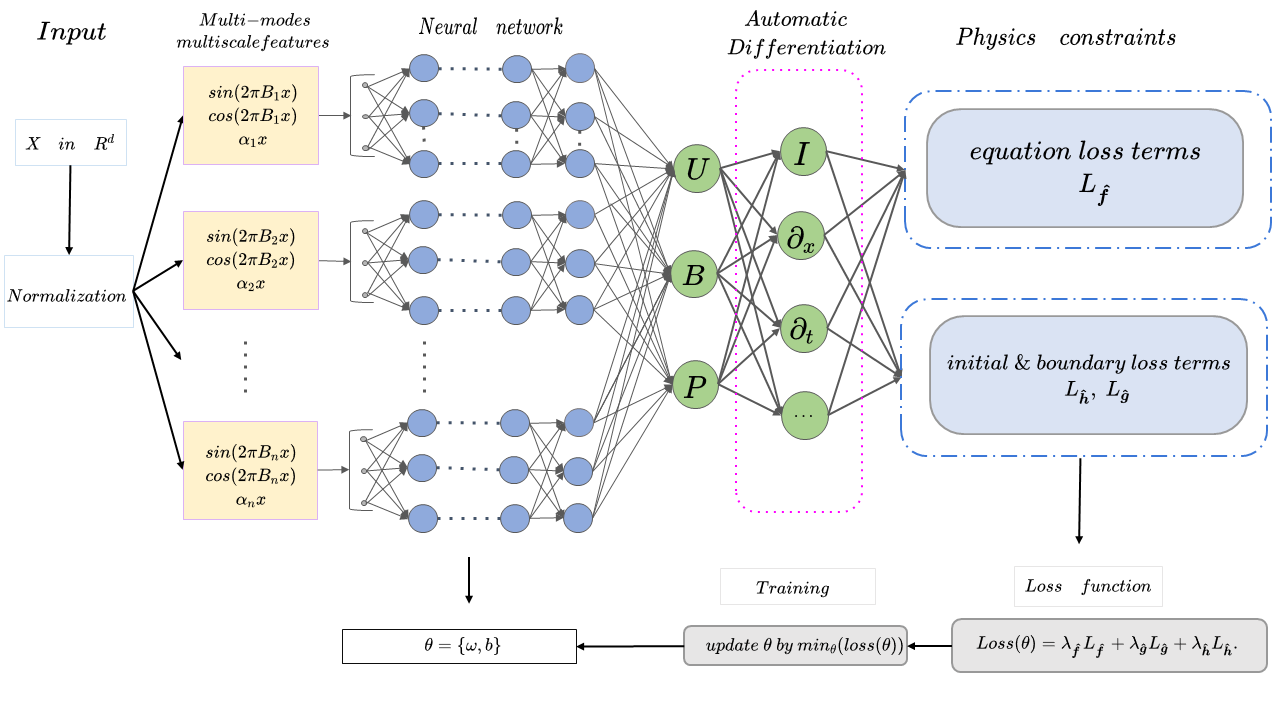}
        \caption{ The framework of MHDnet }
        \label{multiscale}
    \end{minipage}
\end{figure}

\subsection{Multi-modes feature embedding method}
In addition to choosing appropriate mathematical numerical formulations as the loss function, the feature embedding method can also alleviate the convergence difficulty of neural networks when dealing with complex MHD systems. This can be clarified by the neural tangent kernel (NTK) theory \cite{jacot_neural_2020} with spectral biases \cite{Neal1996PriorsFI,lee_deep_2018} or frequency principles \cite{FrequencyPrinciple}. Under the NTK theory, the kernel can be defined as follows 
\begin{equation}
    \boldsymbol{K}_{i j}=\boldsymbol{K}(\boldsymbol{x}_{i}, \boldsymbol{x}_{j})=E_{\theta} \big\langle \frac{\partial NN(\boldsymbol{x}_{i}, \boldsymbol{\theta})}{\partial \boldsymbol{\theta}}, \frac{\partial NN(\boldsymbol{x}_{j}, \boldsymbol{\theta})}{\partial \boldsymbol{\theta}}\big\rangle,
\end{equation}
where $NN (\boldsymbol{x},\boldsymbol{\theta})$ is a fully connected neural network whose weight $\theta$ follows a Gaussian distribution $N(0,1)$. According to the asymptotic conditions in \cite{lee_wide_2020} and the spectral decomposition $\boldsymbol{K}=\boldsymbol{Q} \boldsymbol{\Lambda} \boldsymbol{Q}^{T}$, the training error can be decomposed in the eigenspace of the NTK
\begin{equation}
    \begin{aligned}
        \hat{\boldsymbol{y}}^{(m)}-\boldsymbol{y}  =
     \sum_{i=1}^{N}\left(e^{-\lambda_{i} t} \boldsymbol{q}_{i}^{T} \boldsymbol{Y}_{\text {train }}\right) \boldsymbol{q}_{i}.
    \end{aligned}
    \label{eq:convergence}
\end{equation}
Here, $Q$ is an orthogonal matrix whose ith column is the eigenvector $q_i$ of K and $\Lambda$ is a diagonal matrix whose diagonal entries $\lambda_i$ are the corresponding eigenvalues. $\{\boldsymbol{X}_{\text {train }}, \boldsymbol{Y}_{\text {train }}\}$ is a given dataset, where $\boldsymbol{X}_{\text {train }} = \{x_i\}_{i=1}^N$ are inputs and $\boldsymbol{Y}_{\text {train }} =\{y_i\}_{i=1}^N$ are the corresponding outputs. $\hat{\boldsymbol{y}}^{(m)}$ denotes the prediction of DNN at iteration m. From \cref{eq:convergence}, it can be concluded that the DNN prefers to approximate the target function in the eigendirections of NTK corresponding to larger eigenvalues. However, the eigenvalues decrease as the frequency of the corresponding eigenfunctions increases, and then the high-frequency components of the target function exhibit a considerably slower convergence rate. 

The natural idea is to design the new feature space of NTK, where the different frequencies of the target function can be easily learned by most of the eigenvalues at the same level. The possible approach is feature embedding methods denoted by $\gamma$, which maps input points $\boldsymbol{x} \in [0, 1)^d$ to a new feature space. Several commonly used paradigms have been summarized:

\begin{itemize}
    \item  Original: $\gamma(\boldsymbol{ x})=[\cos (2 \pi \boldsymbol{x} ), \sin (2 \pi \boldsymbol{x})]^{T}$. $\boldsymbol{x} \in R^d$ is the input coordinates. It is inspired by random feature mapping $cos(\omega \boldsymbol{x} + b)$ where $\omega \in R^d$ and $b \in R$ are random variables \cite{mildenhallNeRFRepresentingScenes2020}.
    
    \item  Gaussian: $\gamma(\boldsymbol{x})=[\cos (2 \pi \boldsymbol{B x}), \sin (2 \pi \boldsymbol{B x})]^{\mathrm{T}}$. Here, each entry in the matrix $\boldsymbol{B} \in \mathbb{R}^{m \times d}$(m is a given integer) is sampled from a normal distribution with mean 0 and variance $\sigma^2$, where $\sigma$ is a hyperparameter \cite{liDeepDomainDecomposition2022a}. 

    \item Multiscale: $\gamma(\boldsymbol{x})= \alpha * \boldsymbol{x}$. The high-frequency parts of the target function can be transformed to low-frequency by scaling the input data with a factor $\alpha$, which is usually taken as the power of 2 \cite{tancikFourierFeaturesLet2020}.

    \item  Generalized positional encoding: $\gamma(\boldsymbol{x})=\left[\ldots, \cos \left(2 \pi \sigma^{j / m} \boldsymbol{x}\right), \sin \left(2 \pi \sigma^{j / m} \boldsymbol{x}\right), \ldots\right]^{T}$ for $j=0, \ldots, m-1$. This mapping is deterministic, and the selection of hyperparameters $\sigma$ depends on prior knowledge. For example, the famous NeRF \cite{mildenhallNeRFRepresentingScenes2020} can be regarded as a special case of position coding:
\begin{equation}
              \gamma(\boldsymbol{x})=\left[\sin (\boldsymbol{x}), \cos (\boldsymbol{x}), \ldots, \sin \left(2^{L-1} \boldsymbol{x}\right), \cos \left(2^{L-1} \boldsymbol{x}\right)\right]^{T},
\end{equation}
          and the corresponding matrix form can be expressed as follows
\begin{equation}
              \boldsymbol{P}=\left[\begin{array}{llllllllll}
                      1 & 0 & 0 & 2 & 0 & 0 &        & 2^{L-1} & 0       & 0       \\
                      0 & 1 & 0 & 0 & 2 & 0 & \cdots & 0       & 2^{L-1} & 0       \\
                      0 & 0 & 1 & 0 & 0 & 2 &        & 0       & 0       & 2^{L-1}
                  \end{array}\right]^{T},
              \gamma(\boldsymbol{v})=\left[\begin{array}{l}
                      \sin (\boldsymbol{P} \boldsymbol{x}) \\
                      \cos (\boldsymbol{P} \boldsymbol{x})
                  \end{array}\right].
\end{equation}
\end{itemize}

Inspired by the above methods, we propose a novel \textit{multi-modes feature embedding method},  defined by $\gamma(\boldsymbol{x})=[\cos (2 \pi \boldsymbol{B} \boldsymbol{x}), \alpha*\boldsymbol{x}, \sin (2 \pi \boldsymbol{B} \boldsymbol{x})]^{T}$ whose hyperparameters are similar to Gaussian type. Based on this analysis and discussion, MHDnet has been constructed by the combination of multi-modes feature, multiscale network, and three new mathematical formulations, as shown in \cref{multiscale}. Then, the algorithm procedure of MHDnet is given in \cref{code:algorithm1}.

\begin{algorithm}
    \caption{Algorithm procedure of MHDnet}
    \begin{algorithmic}[1]
        \item \textbf{Input:} $N_r/N_b/N_a$: number of region/boundary/initial collocation points:$\{x_r^i,y_r^i,t_r^i\}_{i=1}^{N_r}$,
        $\{x_b^i, y_b^i ,t_b^i,u_b^i,b_b^i\}_{i=1}^{N_b}$, $\{x_a^i , y_a^i ,t_a^i,u_a^i,b_a^i\}_{i=1}^{N_a}$, $Q_{T} \triangleq \Omega \times[0, T]$. $M$: number of sub-networks. $\{\sigma_i\}_{i=1}^{M}$, Standard deviation of matrix $\{\boldsymbol{B_i}\}_{i=1}^{M}$ sampled from the normal distribution in Fourier features of the mixed structure.
        \item        \textbf{Initialize:}
        Multi-scale Network architectures $u_{\boldsymbol{\theta}}:Q_T \rightarrow \Omega^n$ and parameters $\boldsymbol{\theta}$.\\
        \item \textbf{while} not converged \text{do:} \\

        \qquad \textbf{for} $i=1, \ldots, M$ \textbf{do:}\\
        \qquad \qquad   Preprocessing sampled data using multi-modes features 
        \\
        \qquad \qquad Calculate the loss $Loss(\theta)$ of forward propagation process.\\
        \qquad \textbf{end for} \\
        \qquad Update $\boldsymbol{\theta} \leftarrow \boldsymbol{\theta}-\eta \nabla_{\boldsymbol{\theta}} L$ where $\nabla_{\boldsymbol{\theta}} L$  is approximated using  sample points. \\
            \item  \textbf{Output:} solution  $u_{\boldsymbol{\theta}}(x,y,t)$ in $Q_{T}$;
    \end{algorithmic}
            \label{code:algorithm1}
\end{algorithm}

\section{Numerical experiments}\label{sec:num}
\par In this section, the proposed MHDnets are applied to simulate different incompressible MHD problems, including 2D,3D cases. The comparisons between MHDnet and traditional PINN with three different mathematical formulations are presented, and the influence of different Reynolds numbers is investigated on the accuracy of MHDnet solutions. Additionally, the 2D steady-state case with missing or noisy boundary conditions is also pursued to demonstrate the robustness of MHDnet. The performance of MHDnet is evaluated by the following definition of the relative $L_2$ errors,
\begin{equation}
 \epsilon_{u_i}=\frac{\Vert \hat{u_i}-u_i \Vert _2}{\Vert u_i \Vert _2},
    \epsilon_{b_i}=\frac{\Vert \hat{b_i}-b_i \Vert _2}{\Vert b_i \Vert _2},
    \epsilon_p=\frac{\Vert \hat{p}-p \Vert _2}{\Vert p \Vert _2},   
\end{equation}
where $u_i$, $b_i$, and $p$ refer to the components of velocity $(u_1, u_2, u_3)$, magnetic field $(b_1, b_2, b_3)$, and pressure $p$, respectively. Their hat notations denote the corresponding predicted solutions of MHDnet.

\par In order to optimize the loss function of machine learning models, classical optimization methods such as Stochastic Gradient Descent (SGD) and Adam \cite{adam}, are frequently used, which are based on first-order methods. However, second-order methods like Newton's method have a faster convergence speed while the computation of the Hessian matrix inversion required is computationally expensive. To address this issue, the Limited-memory Broyden–Fletcher–Goldfarb–Shanno (L-BFGS) method \cite{LBFGS} is used. This method is computationally efficient and has a fast convergence speed due to its limited memory usage compared to other second-order methods. In this paper, all the examples are pre-trained with the Adam optimizer, and subsequently, further gradient descent is performed using the L-BFGS method. This hybrid approach combines the advantages of both first and second-order methods to maximize the optimization benefits with minimizing computational expenses.

\subsection{2D steady case}
\label{Sect:2Dsteady}

 \par  We use the following two-dimensional steady MHD equation to demonstrate the performance of MHDnet,
\begin{equation}
\begin{aligned}
& u_1(x, y)=1-e^{\eta x} \cos (2 \pi y), \\
& u_2(x, y)=\frac{\eta}{2 \pi} e^{\eta x} \sin (2 \pi y), \\
& p(x, y)=\frac{1}{2}\left(1-e^{2 \eta x}\right),\\
& b_1(x, y)=sin(y),b_2(x,y)=sin(x)
\end{aligned}
\end{equation}
 
where
\begin{equation*}
    \operatorname{Re}=40,\quad
    \operatorname{Rm}=1,\quad
    \mu=\kappa=1,\quad
 \nu=\frac{1}{\operatorname{Re}}, \quad
    \eta=\frac{1}{2 \nu}-\sqrt{\frac{1}{4 \nu^2}+4 \pi^2}.
\end{equation*}
The calculation domain is $[-0.5,1]\times [0.5,0.5]$, and the weight of boundary loss term is $\lambda_{\boldsymbol{\hat{g}}}=100$. Hyperparameters of multi-modes feature embedding method are chosen as M=4, and ${\sigma_i}=0.1*i,(i=1,\cdots,4)$. Latin hypercube sampling is utilized to ensure adequate coverage, where 100 points with fixed spatial coordinates at each boundary and 2500 points inside the domain are used. Two-stage optimization (Adam optimizer with a learning rate of $10^{-3}$, followed by L-BFGS-B) is used to train the neural network. A sub-network of MHDnet is set up with 4 hidden layers of 50 neurons each, where the number of subnets is $M=4$. Traditional PINN with four hidden layers of 200 neurons is trained for comparison under the same conditions.

\begin{table}[htb]
    \centering
    \begin{tabular}{lllllllll}
        \hline                          & $\epsilon_{u_1}$ & $\epsilon_{u_2}$ & $\epsilon_{b_1}$ & $\epsilon_{b_2}$ & $\epsilon_{p}$ \\
        \hline {$MHDnet-A_2$}                   & $4.62e-5 $       & $3.87e-4$        & $3.33e-5$        & $1.76e-5 $       & $2.66e-4$      \\
        \hline{$PINN-A_2 $}                     & $4.16e-3$        & $2.44e-2$        & $9.83e-3$        & $5.15e-3 $       & $1.07e-0 $     \\
        \hline {$MHDnet-A_1$}            & $2.30e-3 $       & $5.11e-2$        & $2.72e-2$        & $6.90e-3 $       & $3.51e-2$      \\
        \hline{$PINN-A_1$}                & $2.58e-3$        & $6.04e-2$        & $2.57e-2$        & $7.28e-3 $       & $4.03e-2 $     \\
        \hline {$MHDnet-B $}                    & $6.67e-5 $       & $6.04e-4$        & $1.12e-4$        & $6.13e-5 $       & $2.33e-4$      \\
        \hline{$PINN-B $}                     & $3.85e-1$        & $7.28e-1$        & $1.11e-1$        & $8.09e-2 $       & $1.10e-0 $     \\
        \hline
    \end{tabular}
    \caption{2D steady case: relative $L_2$ errors of the velocity field, magnetic field, and pressure field with the number of residual points inside the domain and on the boundary are 2500 and 400, respectively. }
    \label{tab:2d-2-1}
\end{table}

\begin{figure}[H]
    \centering
    \subfigure[Exact]{
        \begin{minipage}[t]{0.2\linewidth}
            \centering
            \includegraphics[width=1.5in]{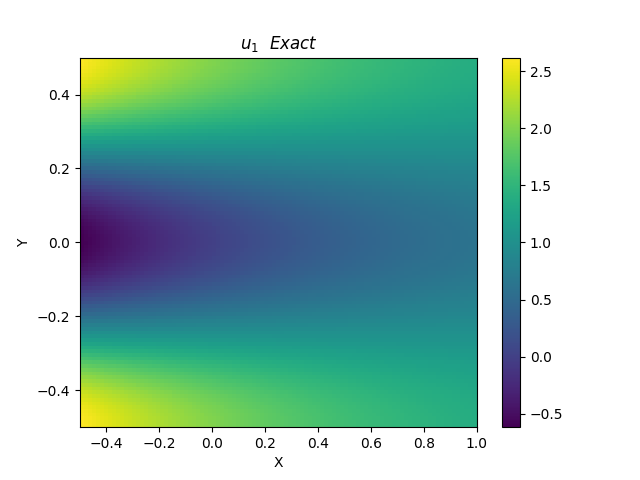}
            \hspace{0.2cm}
            \includegraphics[width=1.5in]{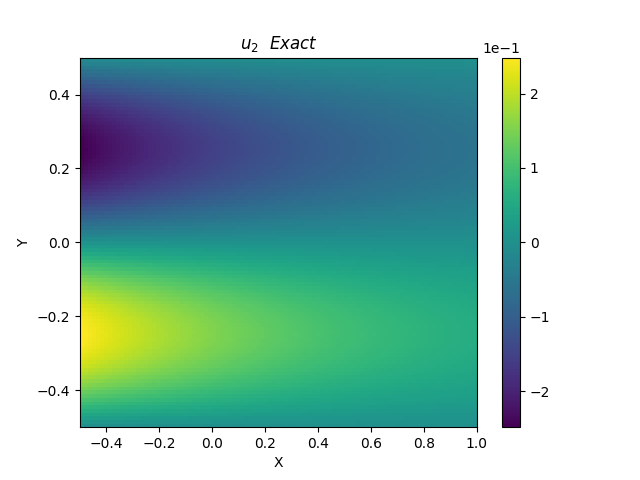}
            \hspace{0.2cm}
            \includegraphics[width=1.5in]{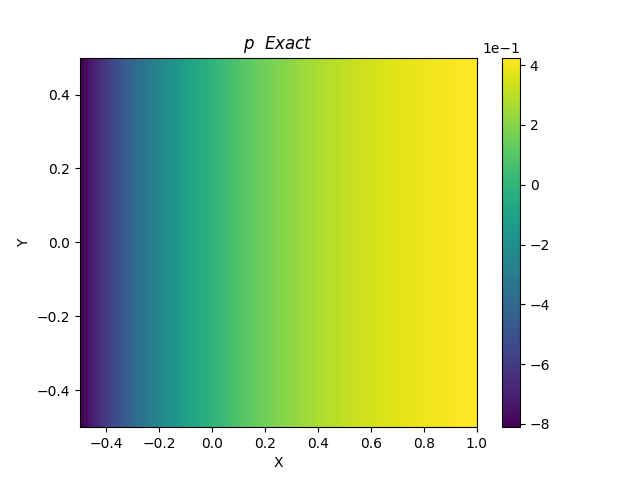}
            \hspace{0.2cm}
            \includegraphics[width=1.5in]{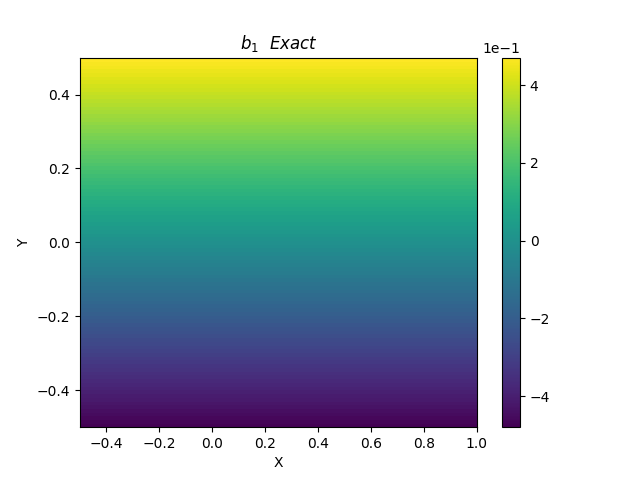}
            \hspace{0.2cm}
            \includegraphics[width=1.5in]{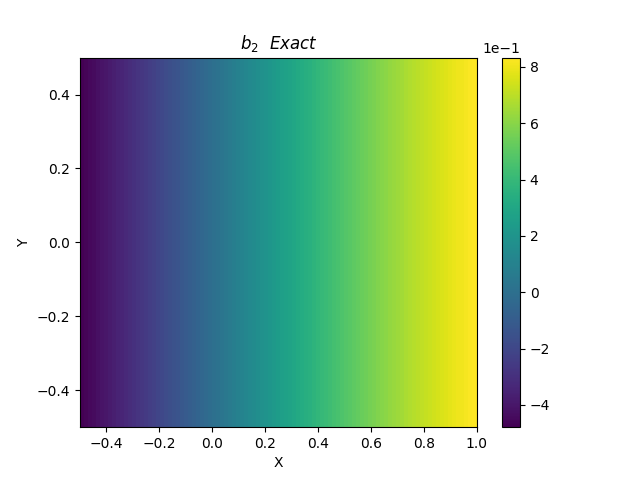}
            \hspace{0.2cm}
        \end{minipage}}%
    \subfigure[$MHDnet-A_2$]{
        \begin{minipage}[t]{0.2\linewidth}
            \centering
            \includegraphics[width=1.5in]{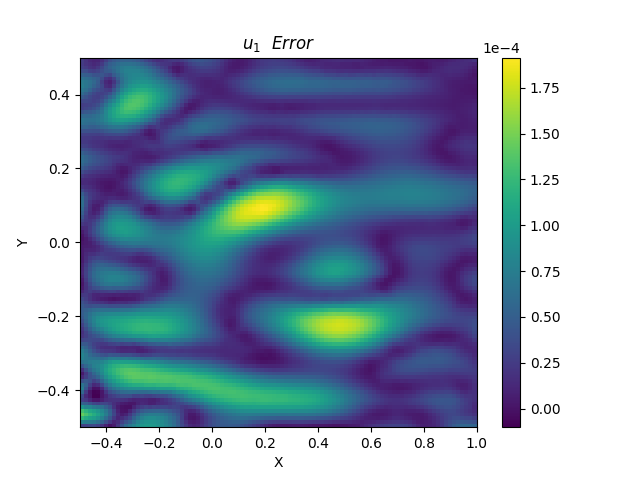}
            \hspace{0.2cm}
            \includegraphics[width=1.5in]{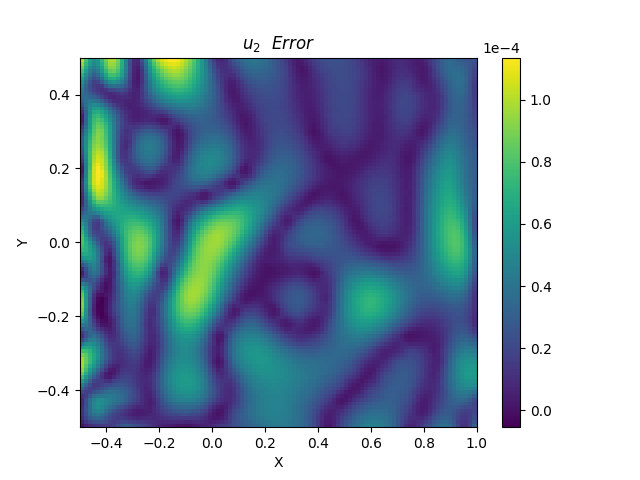}
            \hspace{0.2cm}
            \includegraphics[width=1.5in]{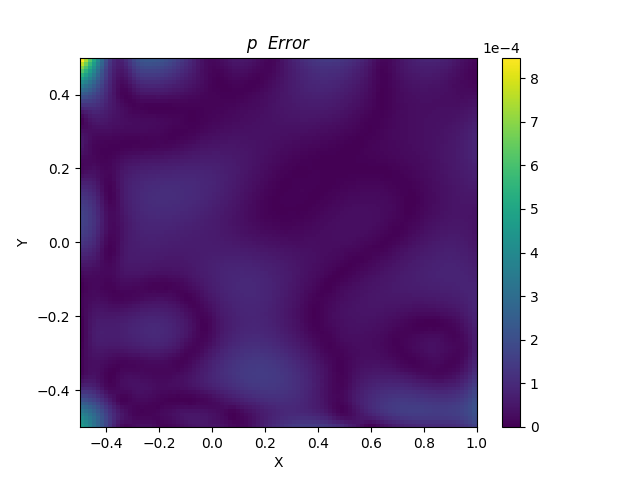}
            \hspace{0.2cm}
            \includegraphics[width=1.5in]{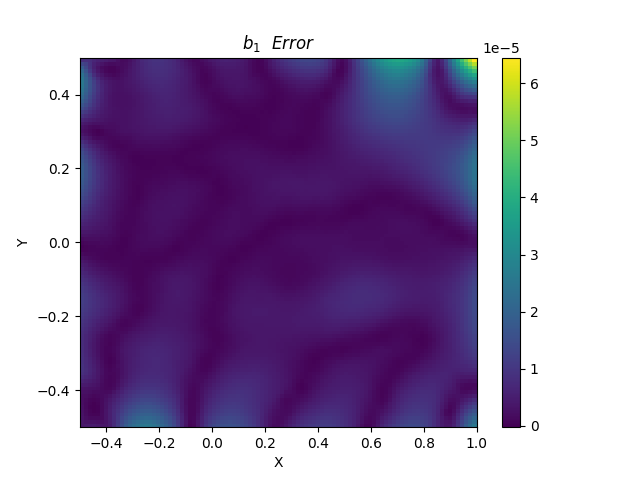}
            \hspace{0.2cm}
            \includegraphics[width=1.5in]{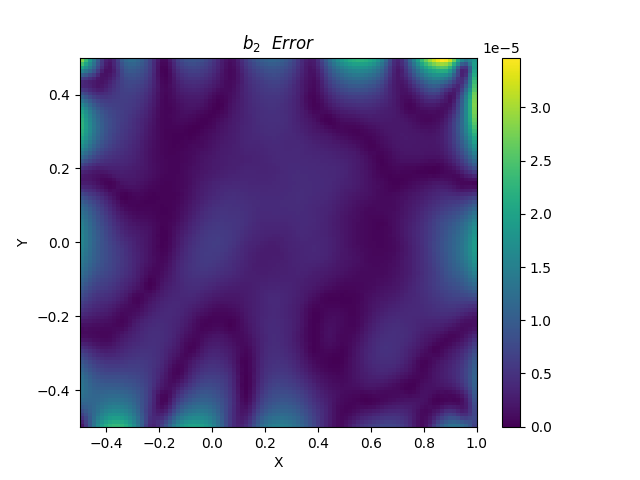}
            \hspace{0.2cm}
        \end{minipage}%
    }%
    \subfigure[$MHDnet-B$]{
        \begin{minipage}[t]{0.2\linewidth}
            \centering
            \includegraphics[width=1.5in]{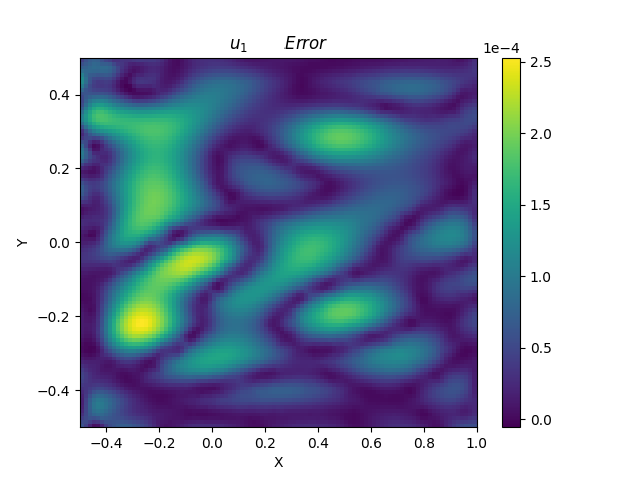}
            \hspace{0.2cm}
            \includegraphics[width=1.5in]{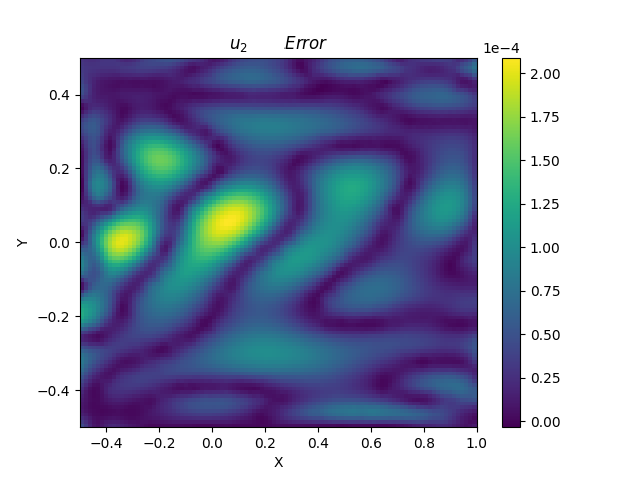}
            \hspace{0.2cm}
            \includegraphics[width=1.5in]{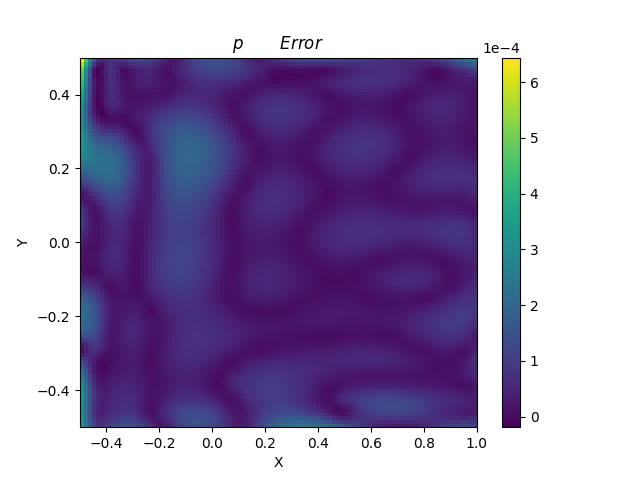}
            \hspace{0.2cm}
            \includegraphics[width=1.5in]{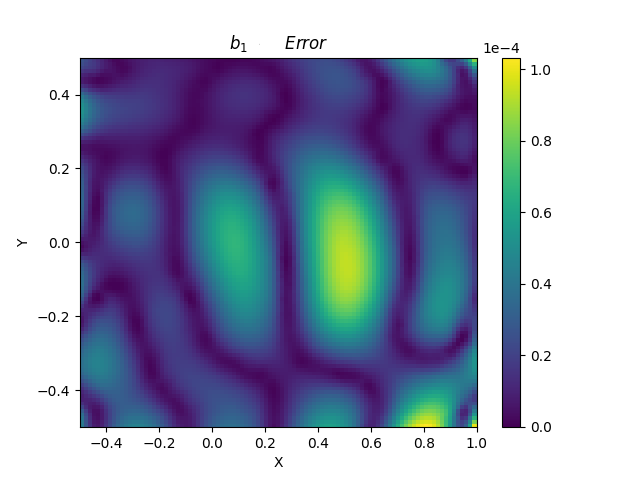}
            \hspace{0.2cm}
            \includegraphics[width=1.5in]{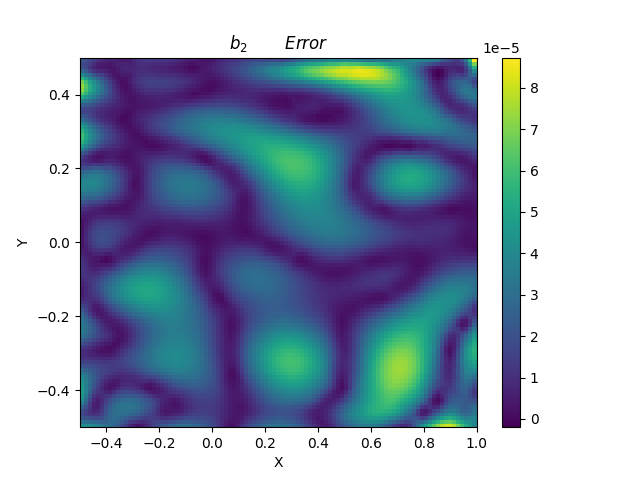}
            \hspace{0.2cm}
        \end{minipage}%
    }%
    \subfigure[$PINN-A_2$]{
        \begin{minipage}[t]{0.2\linewidth}
            \centering
            \includegraphics[width=1.5in]{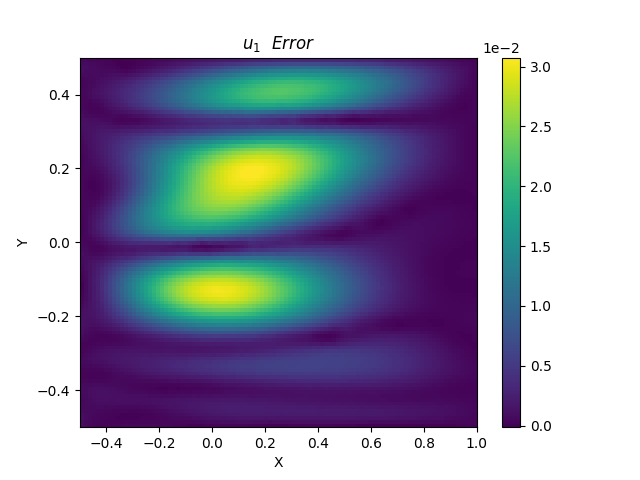}
            \hspace{0.2cm}
            \includegraphics[width=1.5in]{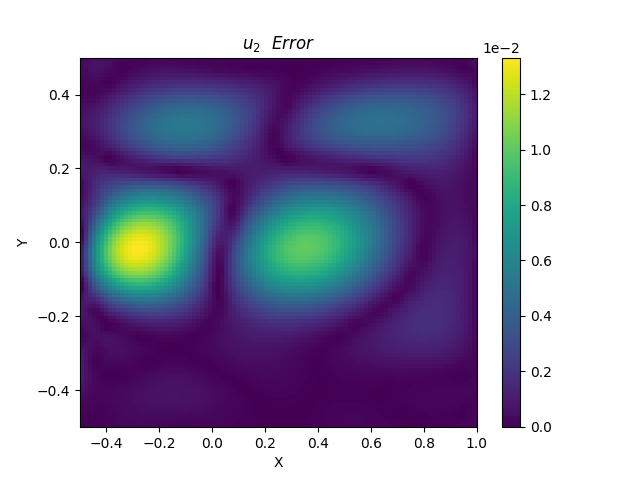}
            \hspace{0.2cm}
            \includegraphics[width=1.5in]{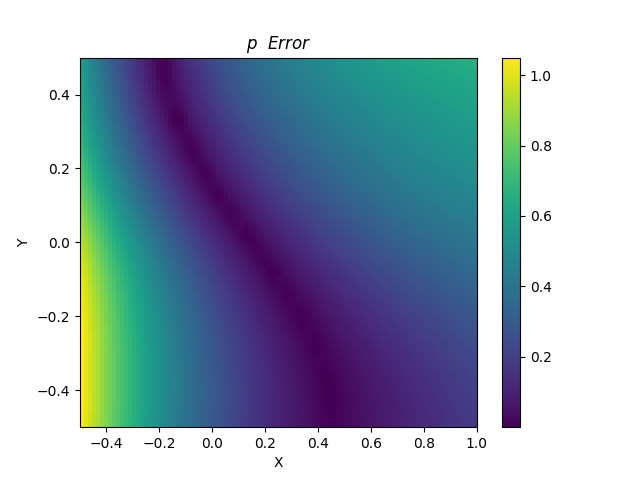}
            \hspace{0.2cm}
            \includegraphics[width=1.5in]{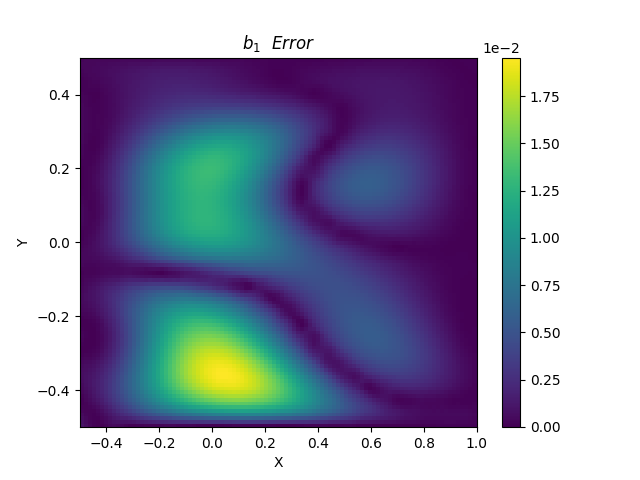}
            \hspace{0.2cm}
            \includegraphics[width=1.5in]{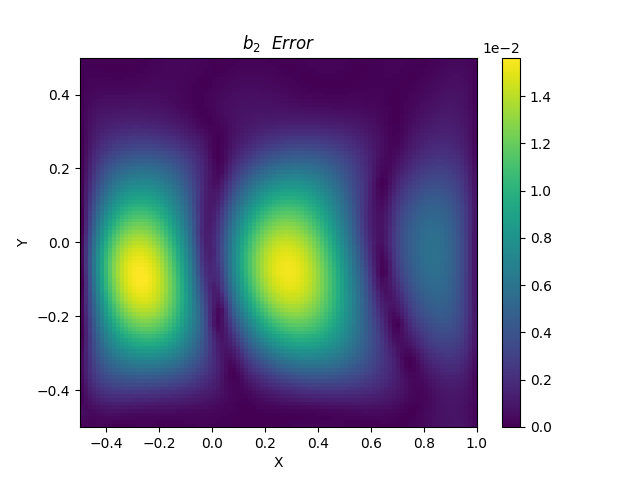}
            \hspace{0.2cm}
        \end{minipage}%
    }%
    \subfigure[$PINN-B$]{
        \begin{minipage}[t]{0.2\linewidth}
            \centering
            \includegraphics[width=1.5in]{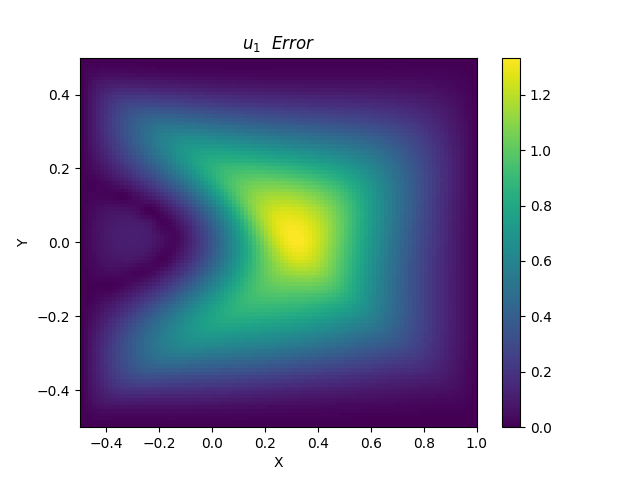}
            \hspace{0.2cm}
            \includegraphics[width=1.5in]{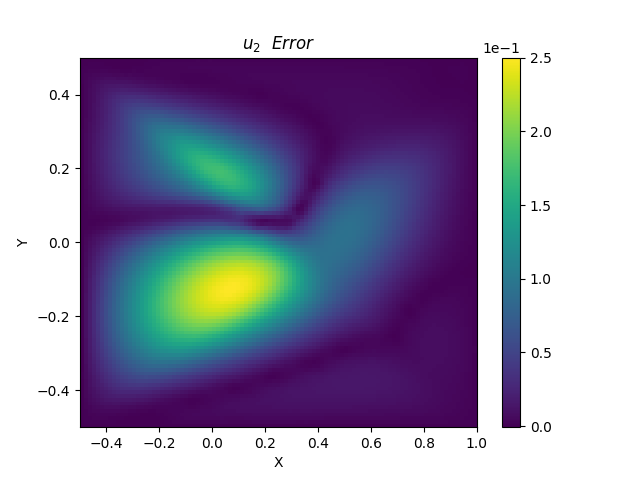}
            \hspace{0.2cm}
            \includegraphics[width=1.5in]{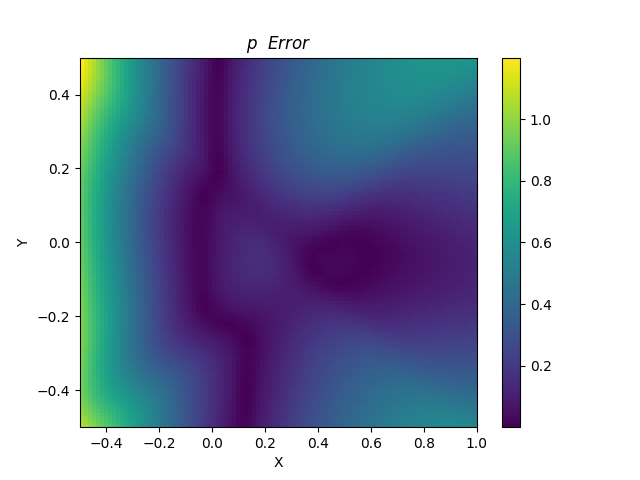}
            \hspace{0.2cm}
            \includegraphics[width=1.5in]{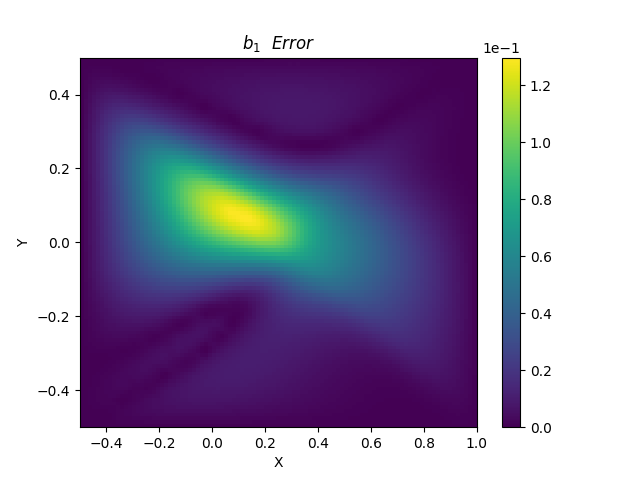}
            \hspace{0.2cm}
            \includegraphics[width=1.5in]{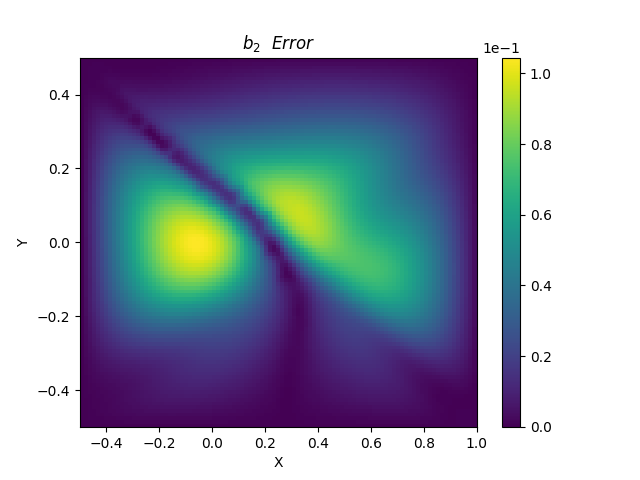}
            \hspace{0.2cm}
        \end{minipage}%
    }%
    \caption{2D steady case: exact solution and absolute error distribution of velocity field, pressure field, and magnetic field predicted by different models. }
    \vspace{-0.2cm}
    \label{fig:2d-2-up}
\end{figure}

\cref{tab:2d-2-1} presents the comparisons between MHDnet and PINN with three different mathematical formulations ($B$, $A_1$, $A_2$) as loss functions. The results show that MHDnet outperforms traditional PINN by several orders of magnitude in terms of the relative $L_2$ error. Although the relative error of PINN significantly decreases from '$B~formulation$' to '$A_2~formulation$', it is still not better than that of MHDnet. This is because the multi-modes multiscale framework of MHDnet can better capture the coupled characteristic of complex multiphysics fields, and the incompressible condition can be naturally incorporated into the fitting process of MHDnet. Moreover, for the hidden pressure field, the MHDnet framework with '$A_2~formulation$' exhibit the lowest relative $L_2$ error. In summary, the $MHDnet-A_2$ yields the best results, and the absolute error distribution of $MHDnet-A_2$ is better than others, which are displayed in \cref{fig:2d-2-up}.

\begin{figure}[htb]
    \centering
    \subfigure[$MHDnet-A_2$ ]{\includegraphics[width=0.3\hsize]{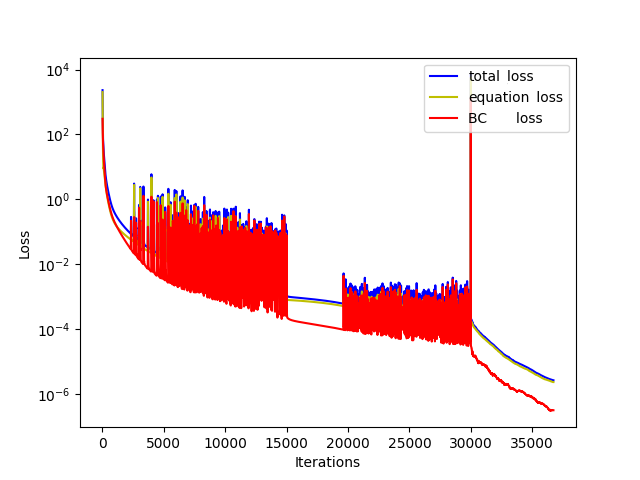}}
    \subfigure[$MHDnet-A_1$ ]{\includegraphics[width=0.3\hsize]{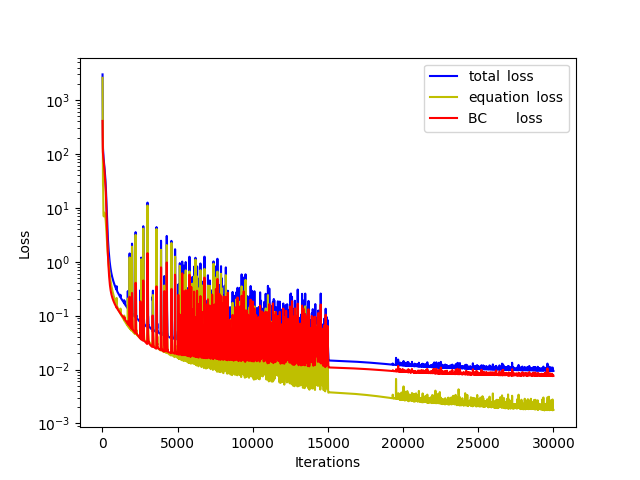}}
    \subfigure[$MHDnet-B$ ]{\includegraphics[width=0.3\hsize]{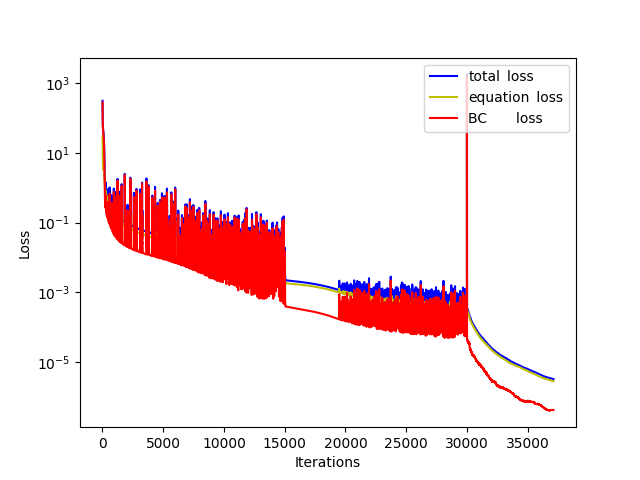}}  \\
    \subfigure[$PINN-A_2$ ]{\includegraphics[width=0.3\hsize]{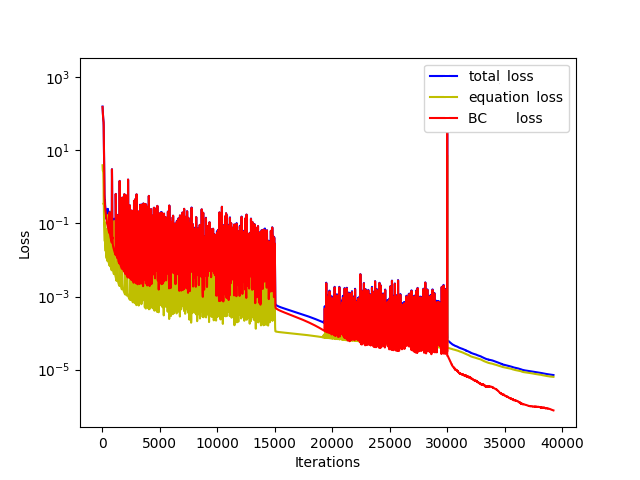}}
    \subfigure[$PINN-A_1$ ]{\includegraphics[width=0.3\hsize]{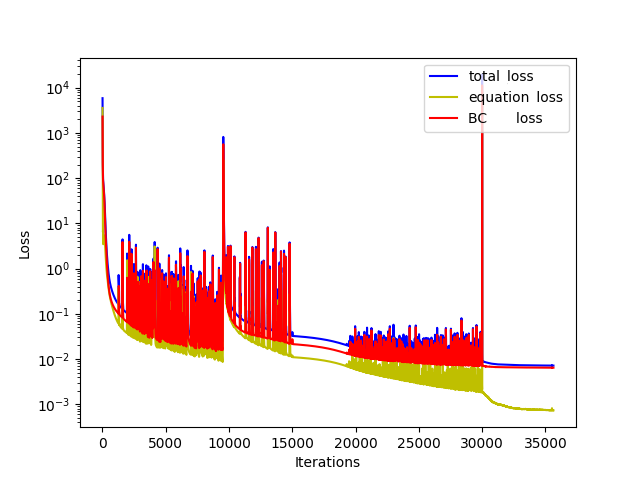}}
    \subfigure[$PINN-B$]{\includegraphics[width=0.3\hsize]{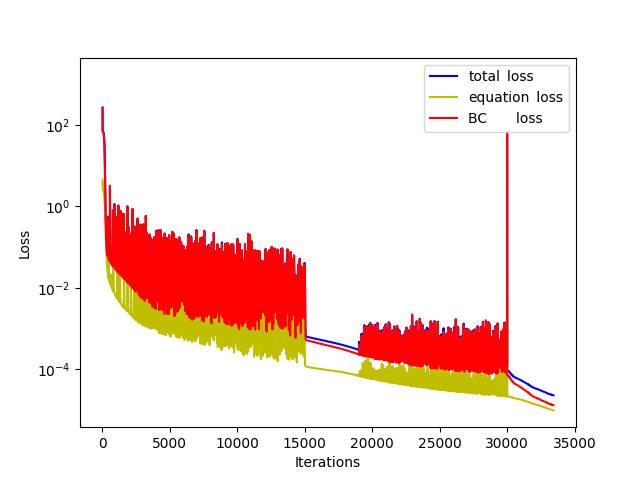}}
    \caption{2D steady case: convergence of total loss function, equation loss function, initial condition and boundary condition loss function (BC $\&$ IC loss).}
    \label{fig:2d-2-loss}
\end{figure}

\cref{fig:2d-2-loss} shows the convergence of loss functions for the three formulations, which can be divided into two parts. The oscillating part is obtained by the Adam optimizer, and the steady decline part is obtained by the L-BFGS optimizer. The $MHDnet-A_2$ achieved a minimum loss level of $1e-6$, which is better than the other architectures, and has fewer iterations than $PINN-A_2$. The \cref{fig:2d-2-loss} (b), (e), (f) also show that the boundary loss is the majority of the total loss for the $PINN-B$ at the end of the training process, which indicates that the network does not fit boundary conditions effectively. This is also the reason why the $MHD-A_1$ and $PINN-A_1$ have similar error levels and $MHDnet-A_2$ can converge faster, as shown in \cref{tab:2d-2-1}.

\begin{figure}[htp]
    \centering
    \subfigure[PINN]{\includegraphics[width=0.9\hsize]{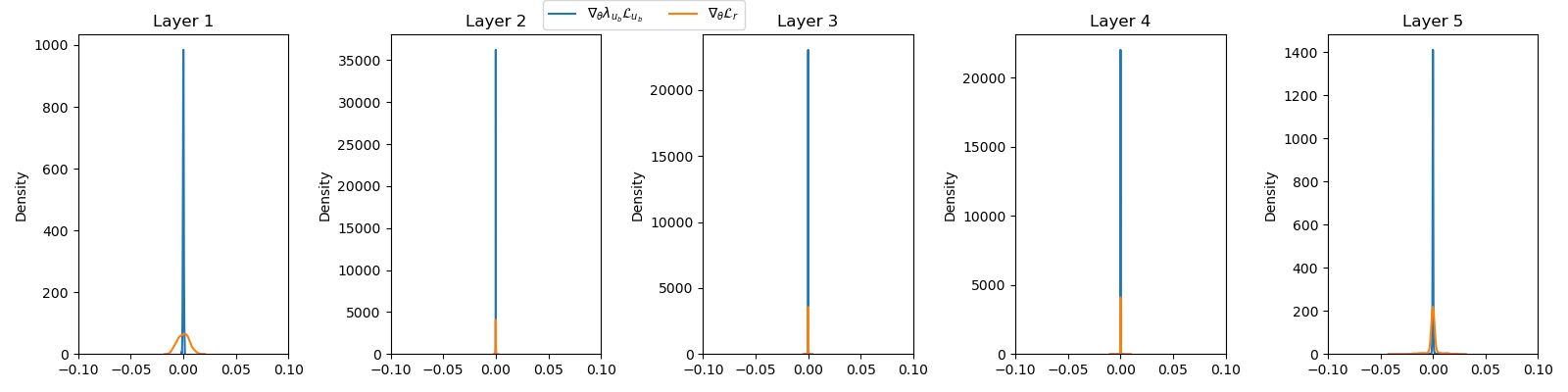}}\\
    \subfigure[subnet1]{\includegraphics[width=0.9\hsize]{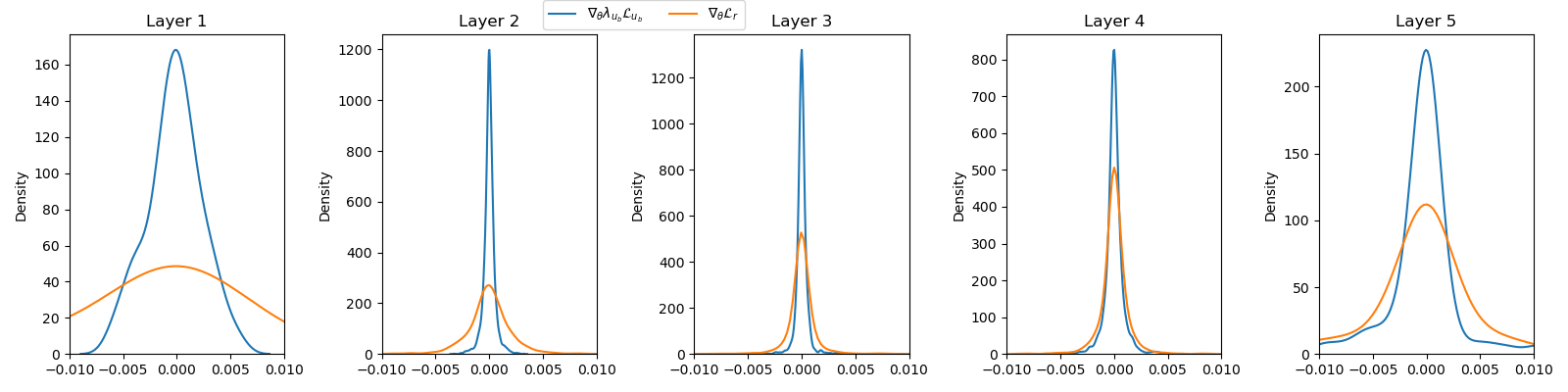}}
    \subfigure[subnet2]{\includegraphics[width=0.9\hsize]{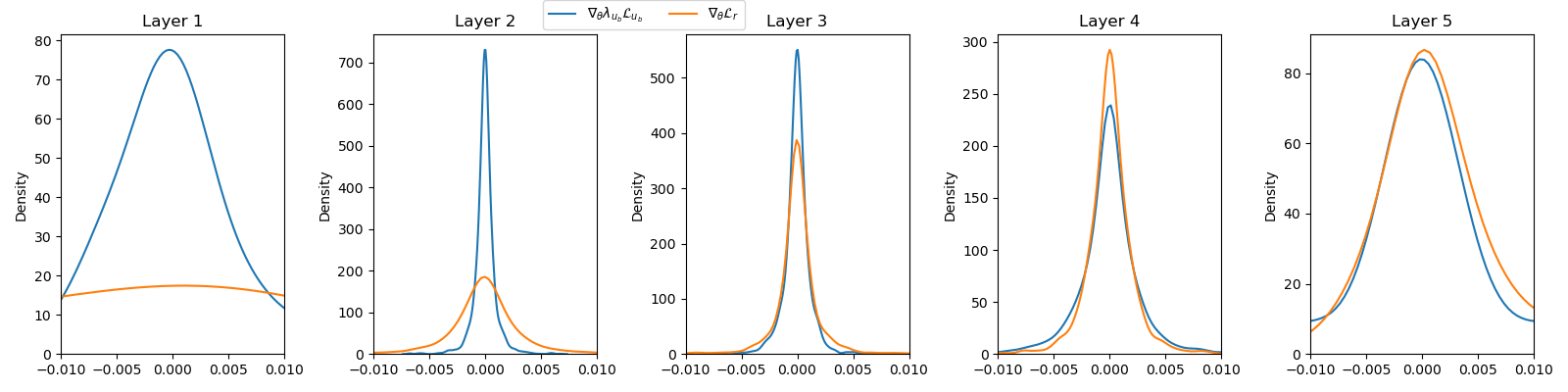}}  \\
    \subfigure[subnet3 ]{\includegraphics[width=0.9\hsize]{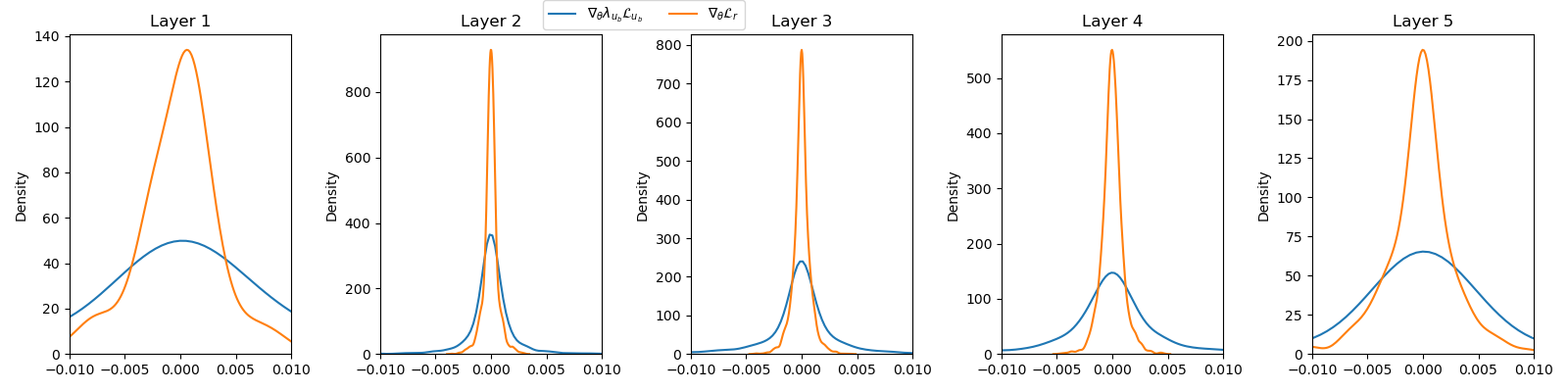}}
    \subfigure[subnet4 ]{\includegraphics[width=0.9\hsize]{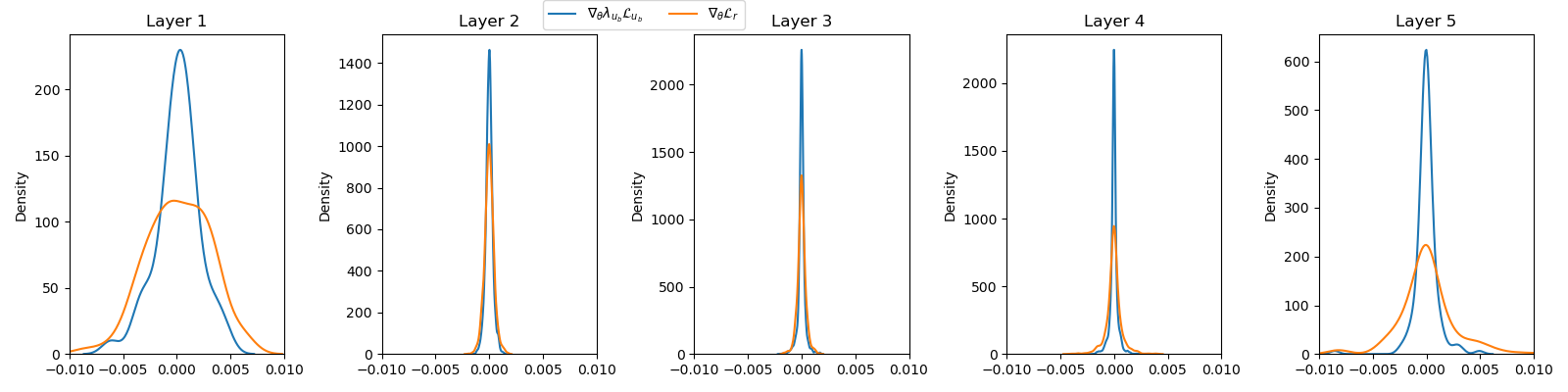}}
    \caption{2D steady case: histogram of back-propagation gradient for different loss terms in each layer of PINN and different subnetworks of MHDnet.}
    \label{fig:2d-2-gradient}
\end{figure}

\cref{fig:2d-2-gradient} displays the backpropagation gradient histogram of $MHDnet-A_2$ and $PINN-A_2$ at the 30000th training epoch. The yellow lines represent the gradient contribution of the equation loss $\nabla_\theta L_{\boldsymbol{\hat{f}}}$, while the blue lines are the gradient contribution of the boundary loss $\nabla\theta \lambda_{\boldsymbol{\hat{g}}} L_{\boldsymbol{\hat{g}}}$. It is obviously seen that, for $PINN-A_2$, the gradient of the hidden layer concentrates  around zero, leading to the gradient vanishing. By contrast, the gradient distribution in $MHDnet-A_2$ is more scattered and the parameters can be updated continuously. This also explains why $MHDnet-A_2$ can converge to a very small value.
 
\begin{figure}[htb]
    \centering
    \subfigure[Exact ]{
        \begin{minipage}[t]{0.33\linewidth}
            \centering
            \includegraphics[width=2 in]{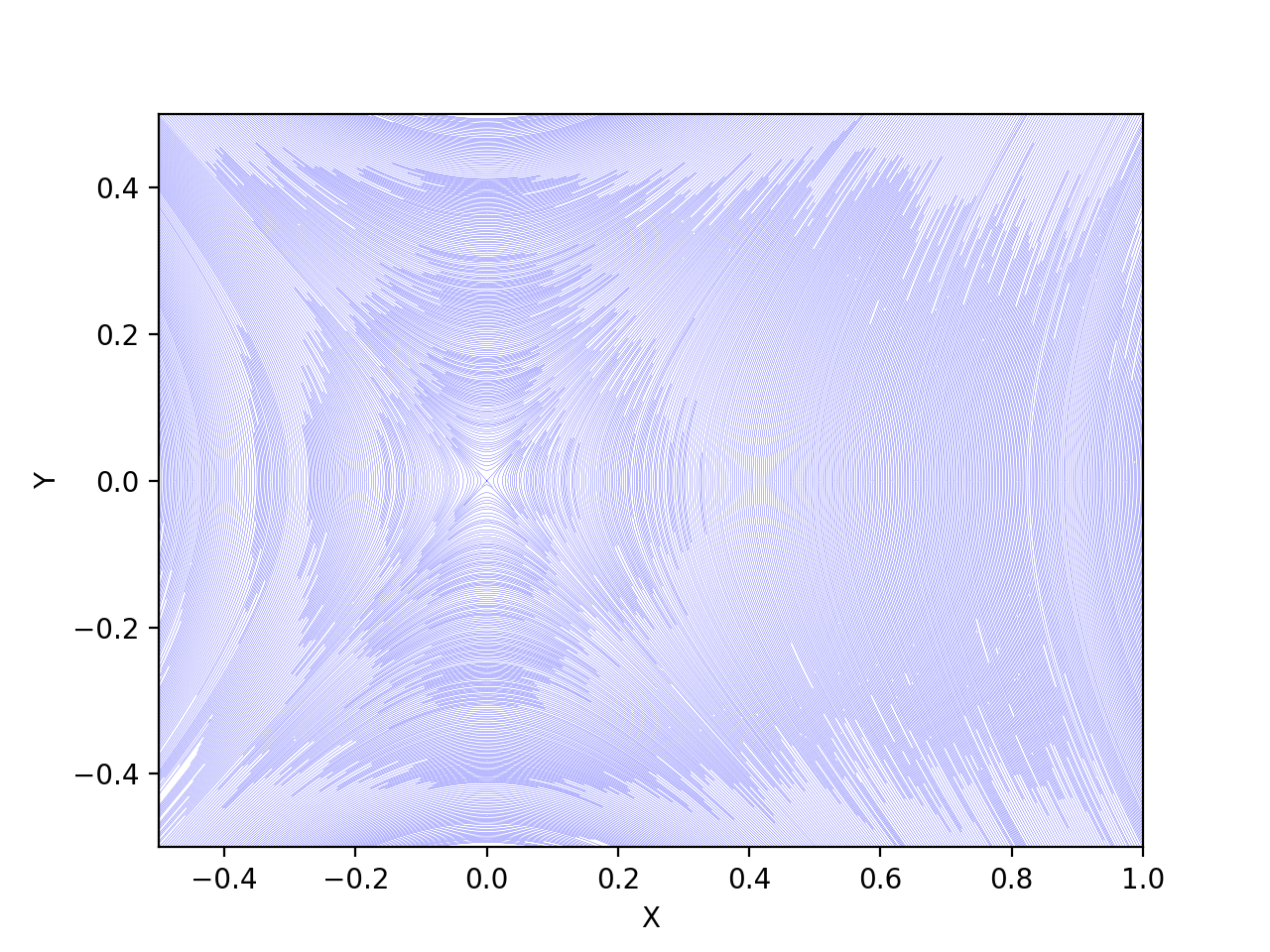}
            \hspace{0.2cm}
            \includegraphics[width=2 in]{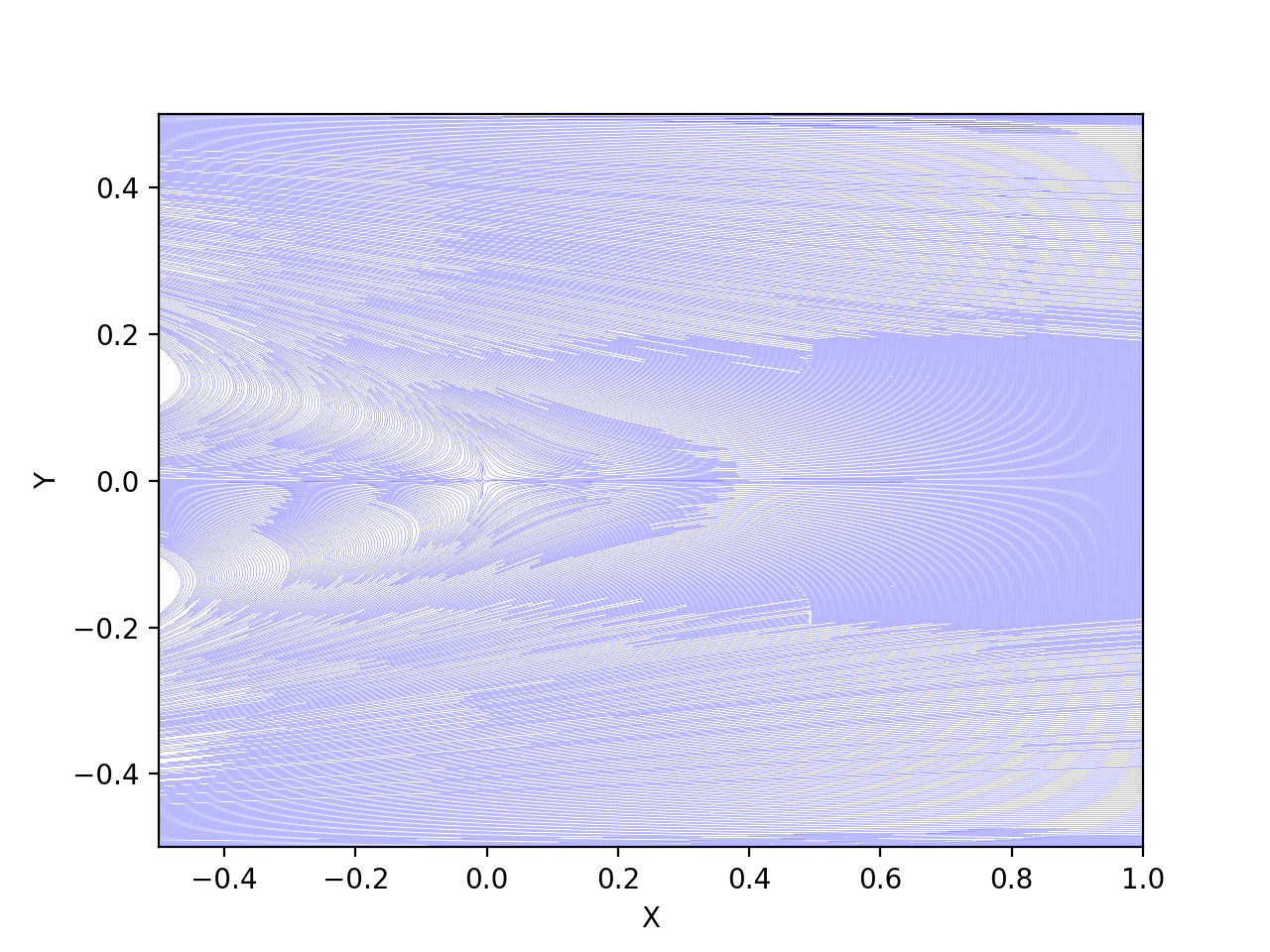}
            \hspace{0.2cm}
        \end{minipage}}%
    \subfigure[$MHDnet-A_2$]{
        \begin{minipage}[t]{0.33\linewidth}
            \centering
            \includegraphics[width=2 in]{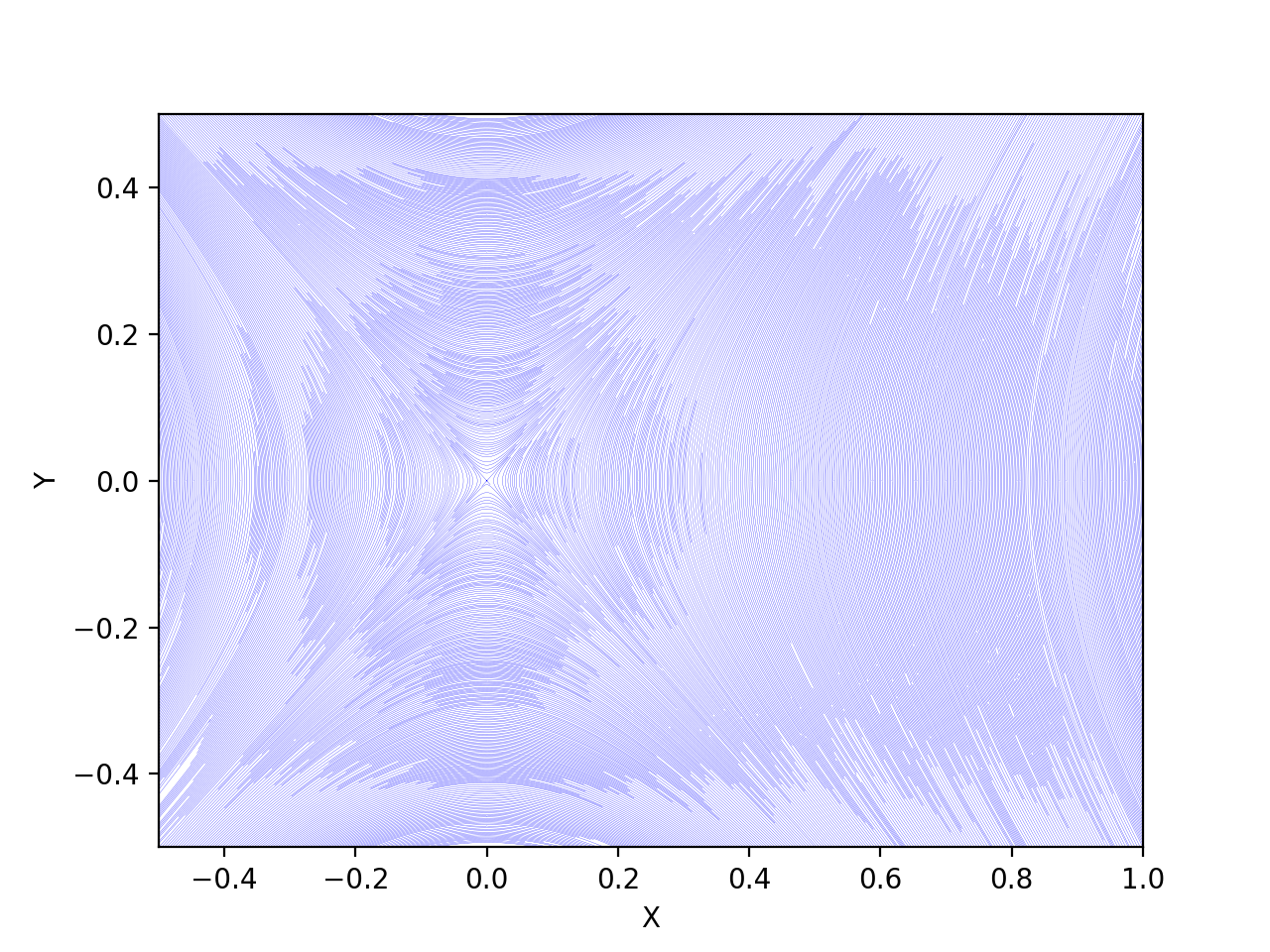}
            \hspace{0.2cm}
            \includegraphics[width=2 in]{new/2d-2-u-stream-MHDnet.png}
            \hspace{0.2cm}
        \end{minipage}%
    }%
    \subfigure[$PINN-A_2$]{
        \begin{minipage}[t]{0.33\linewidth}
            \centering
            \includegraphics[width=2 in]{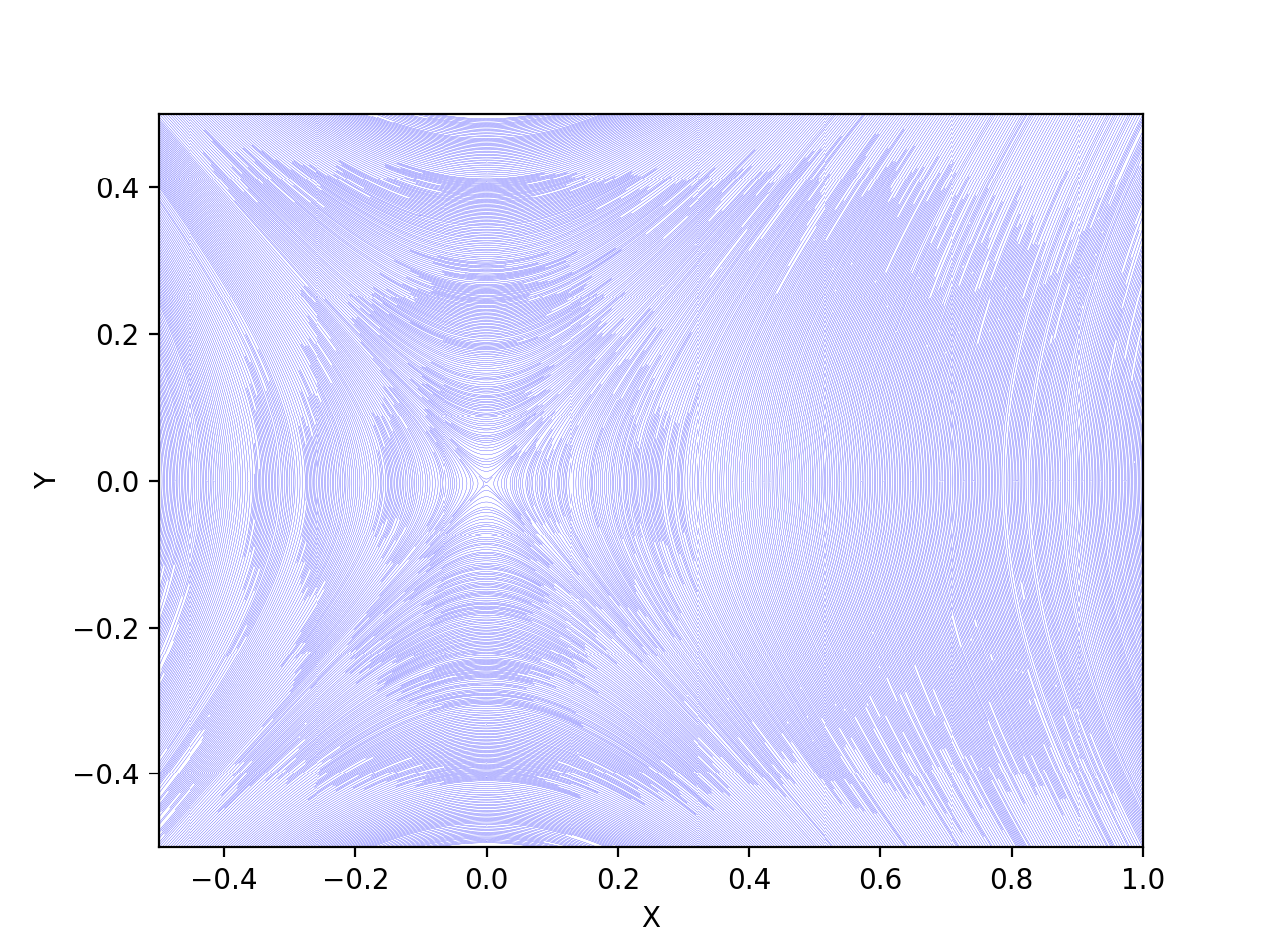}
            \hspace{0.2cm}
            \includegraphics[width=2 in]{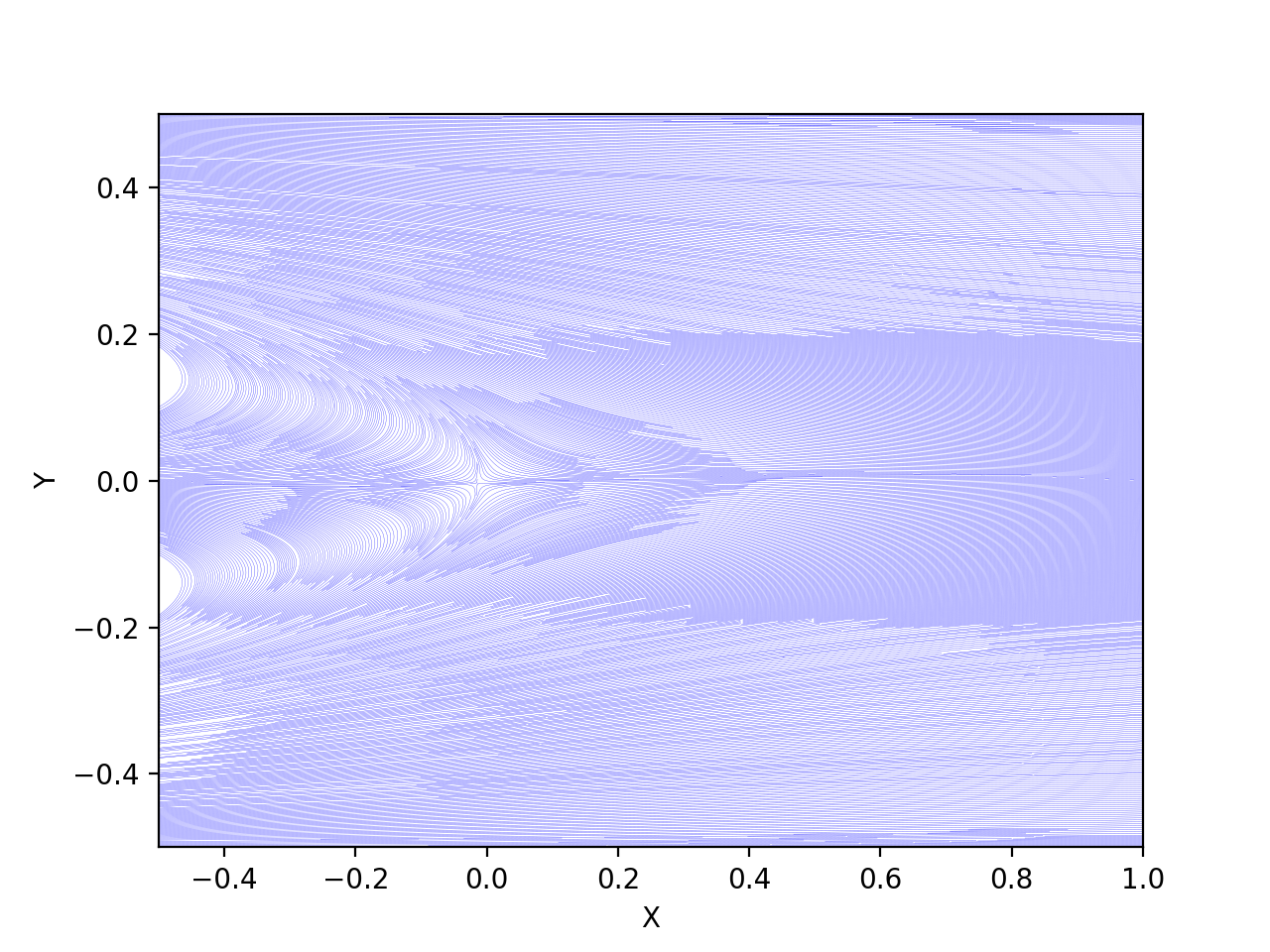}\\
            \hspace{0.2cm}
        \end{minipage}%
    }%
    \caption{2D steady case: the first row is the magnetic streamlines, and the second row is the velocity streamlines. (a) exact solutions, (b) $MHDnet-A_2$ solutions and (c) $PINN-A_2$ solutions.}
    \vspace{-0.2cm}
    \label{fig:2d-stream}
\end{figure}

In \cref{fig:2d-stream}, the first row depicts the magnetic streamlines, and the second row displays the velocity streamlines. The predicted solutions of $MHDnet-A_2$ and $PINN-A_2$ are compared with the exact solutions in $\Omega$. For the velocity streamlines, it is apparent that two vortices are concentrated near the y-axis on the left side of $\Omega$, and $MHDnet-A_2$ successfully captures more information related to velocity and magnetic field.

\begin{figure}[htb]
    \centering
    \subfigure[$MHDnet-A_2$]{\includegraphics[width=0.49\hsize]{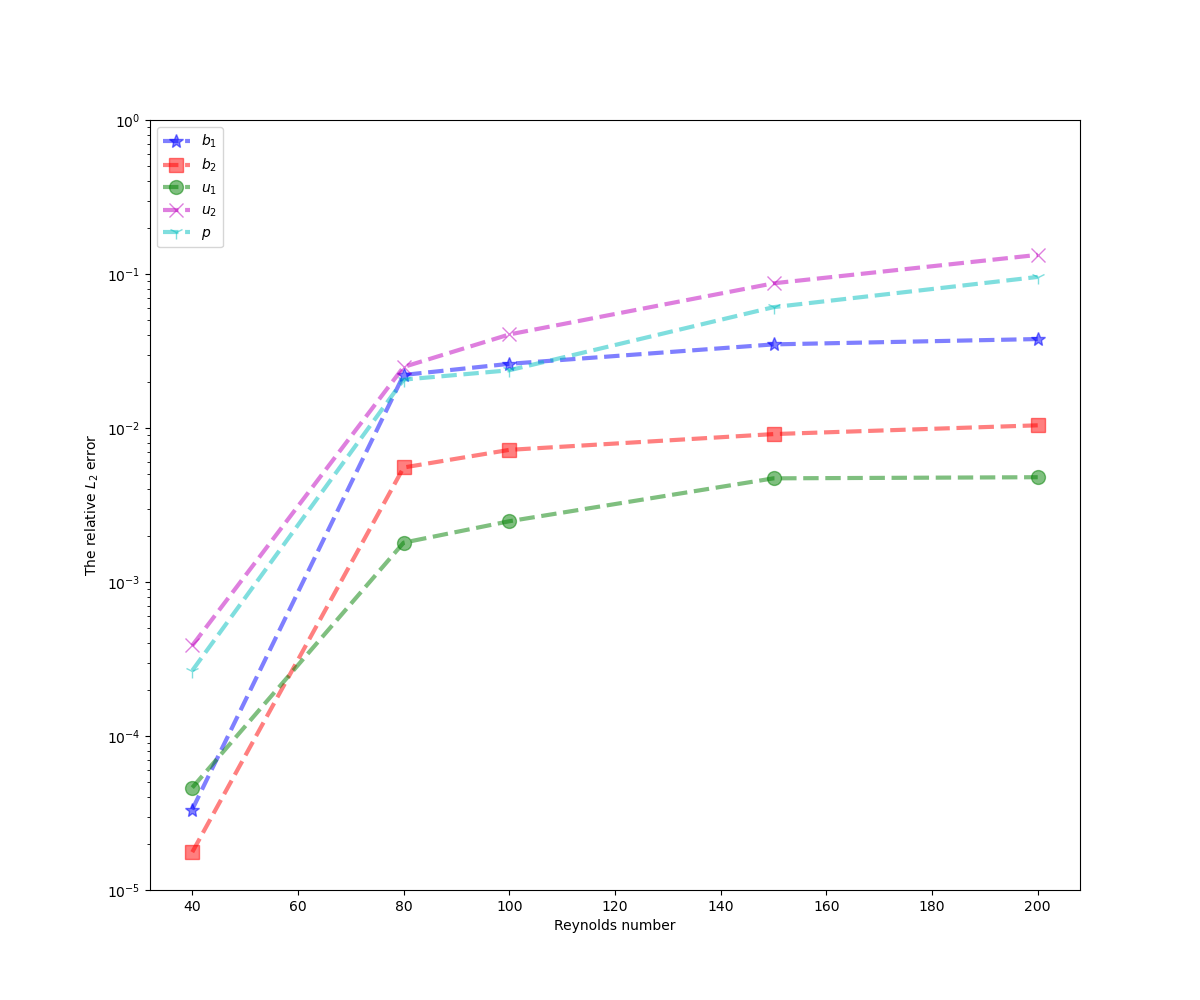}}
    \subfigure[$PINN-A_2$]{\includegraphics[width=0.49\hsize]{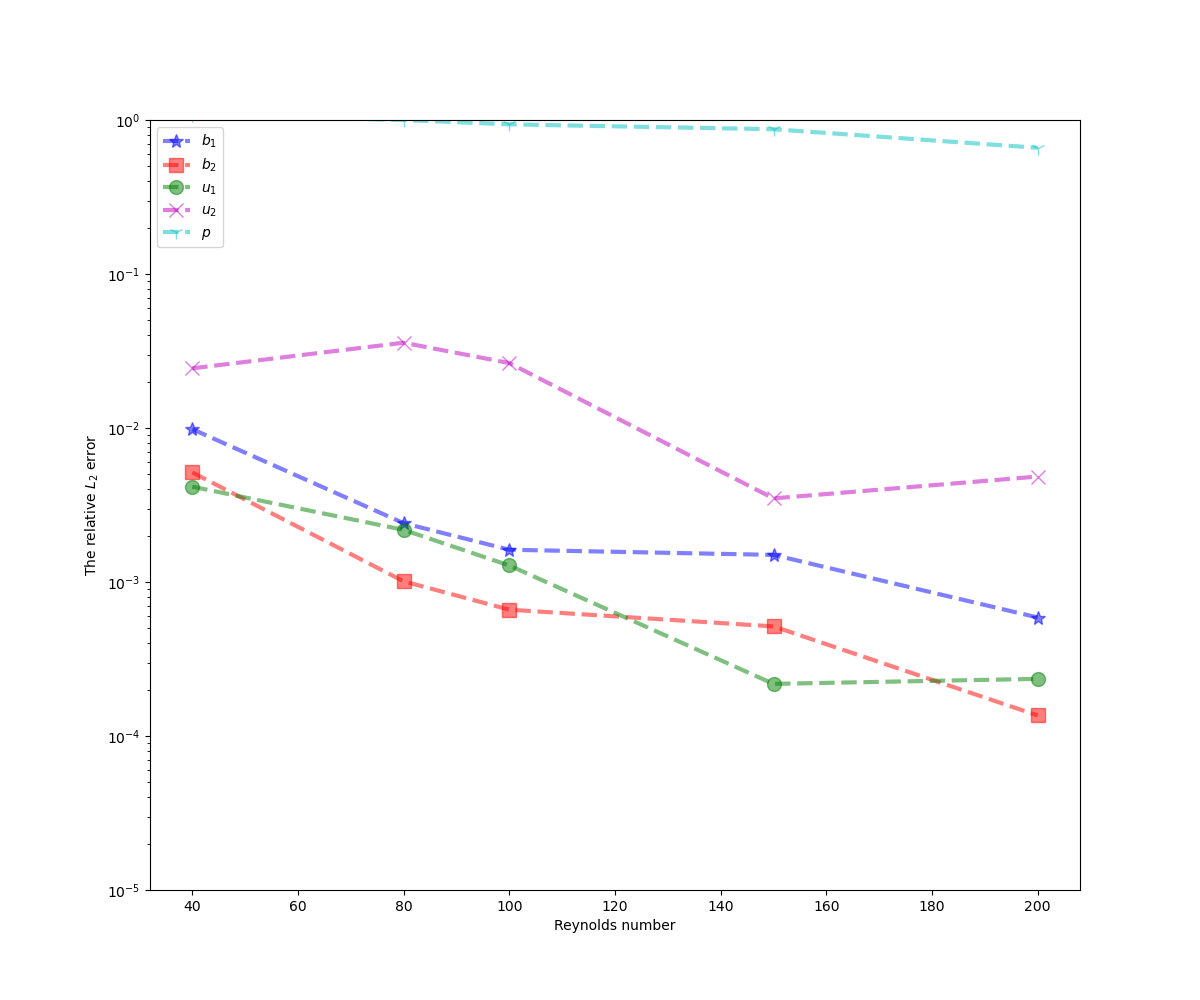}}
    \caption{2D steady case: Performance of $MHDnet-A_2$ and $PINN-A_2$ with different Reynolds Numbers.}
    \label{fig:2d-Re}
\end{figure}

\begin{table}[htb]
    \centering
    \begin{tabular}{llllllll}
        \hline I \& B points  & $\epsilon_{u_1}$ & $\epsilon_{u_2}$ & $\epsilon_{b_1}$ & $\epsilon_{b_2}$ & $\epsilon_{p}$ \\
        \hline {$2500/400$}  & $4.62e-5 $       & $3.87e-4$        & $3.33e-5$        & $1.76e-5 $       & $2.66e-4$      \\
        \hline {$2500/1000$} & $4.80e-5 $       & $4.34e-4$        & $3.04e-5$        & $1.93e-5 $       & $3.08e-4$      \\
        \hline{$5000/400$}   & $4.67e-5 $       & $4.14e-4$        & $4.05e-5$        & $2.22e-5 $       & $3.14e-4$      \\
        \hline{$5000/1000$}  & $2.52e-5 $       & $2.35e-4$        & $1.86e-5$        & $1.10e-5 $       & $1.77e-4$      \\
        \hline
    \end{tabular}
    \caption{2D steady case: relative $L_2$ errors of the velocity field, magnetic field, and pressure field for $MHDnet-A_2$. The subnetworks are trained with different numbers of interior and boundary points (I $\&$ B points).}
    \label{tab:2d-2-2}
\end{table}
 The effect of the Reynolds number and the number of residual points on the performance of MHDnet are also studied. The performance of $MHDnet-A_2$ is compared across different Reynolds numbers, as shown in \cref{fig:2d-Re}. The relative $L_2$ errors are all below $10\%$, and the relative $L_2$ error increases with the Reynolds number, which is consistent with previous research in this field. Moreover, \cref{tab:2d-2-2} exhibits the relative $L_2$ error of MHDnet for a specific number of sampling points, where the increasing sample points can not affect the accuracy of solutions.

\subsection{2D unsteady case}
\label{Sect:2Dunsteady}
In this example, we focus on using $A_2~formulation$ and $B~formulation$ for solving unsteady MHD equations. The domain is $\Omega = [0, 1] \times [0, 1]$. Assuming the following functions represent the exact solution: 

\begin{equation}
\begin{aligned}
& \boldsymbol{u}=(y^5+t^2,x^5+t^2), \\
&  p=10(2x-1)(2y-1) \\
&\boldsymbol{B}=(\sin(8*y)+t^2,\sin(8*x)+t^2)
\end{aligned}
\end{equation}
To satisfy the system requirements, we impose appropriate force fields and set the parameters $R_e=1$, $\kappa=1$, and $R_m=1$, and the time interval $[0,1]$. 100 initial points are added to train dataset, and the other settings are the same as 2D steady case in \cref{Sect:2Dsteady}.

\begin{table}[htb]
    \centering
    \begin{tabular}{lllllll}
        \hline                                    & $\epsilon$       & $t=0.25$     & $t=0.50$    & $t=0.75$    & $t=1.00$    \\
        \hline \multirow{3}{*}{$\mathrm{MHDnet-B}$} & $\epsilon_{u_1}$ & $1.02e-03$  & $6.03e-04$  & $6.25e-04$  & $9.51e-04$  \\
                                                 & $\epsilon_{u_2}$ & $9.36e-04  $ & $7.61e-04 $ & $5.18e-04$  & $1.02e-03$  \\
                                                  & $\epsilon_{b_1}$ & $2.10e-04 $  & $2.19e-04 $ & $2.26e-04 $ & $5.58e-04 $ \\
                                                  & $\epsilon_{b_2}$ & $2.46e-04$  & $2.15e-04 $ & $1.86e-04 $ & $5.11e-04$  \\
                                                  & $\epsilon_{p}$   & $6.73e-03$  & $7.70e-03$ & $9.99e-03 $ & $ 1.40e-02$ \\

        \hline \multirow{3}{*}{$\mathrm{PINN-B}$}   & $\epsilon_{u_1}$ & $3.46e-03 $  & $4.52e-03$  & $4.09e-03$  & $3.15e-03$  \\
                                                  & $\epsilon_{u_2}$ & $2.42e-03 $  & $4.33e-03 $ & $4.75e-03$  & $4.65e-03$  \\
                                                  & $\epsilon_{b_1}$ & $ 6.25e-04 $  & $5.39e-04$  & $5.05e-04 $ & $ 7.31e-04$  \\
                                                  & $\epsilon_{b_2}$ & $5.50e-04$  & $6.75e-04 $ & $6.66e-04$ & $ 9.09e-04$  \\
                                                  & $\epsilon_{p}$   & $8.24e-02$   & $9.67e-02$ & $1.15e-01$  & $1.53e-01$  \\

        \hline \multirow{3}{*}{$\mathrm{MHDnet-A_2}$} & $\epsilon_{u_1}$ & $6.81e-03 $   & $4.90e-03 $  & $4.30e-03 $ & $7.08e-03 $ \\
                                                    & $\epsilon_{u_2}$ & $6.84e-03   $ & $5.73e-03  $ & $3.59e-03 $ & $9.19e-03 $ \\
                                                    & $\epsilon_{b_1}$ & $5.23e-04 $   & $2.23e-04 $  & $3.29e-04 $ & $ 2.02e-04$ \\
                                                    & $\epsilon_{b_2}$ & $3.50e-04 $   & $4.03e-04 $  & $2.13e-04 $ & $5.52e-04$ \\
                                                    & $\epsilon_{p}$   & $1.32e-02 $   & $1.56e-02 $  & $2.59e-02 $ & $4.40e-02$ \\

        \hline \multirow{3}{*}{$\mathrm{PINN-A_2}$}   & $\epsilon_{u_1}$ & $5.01e-03 $   & $5.69e-03$   & $3.54e-03$  & $8.06e-03$  \\
                                                    & $\epsilon_{u_2}$ & $5.07e-03  $  & $6.55e-03 $  & $2.69e-03$  & $7.01e-03 $ \\
                                                    & $\epsilon_{b_1}$ & $9.88e-04 $   & $8.79e-04 $  & $ 7.72e-04$  & $9.23e-04$ \\
                                                    & $\epsilon_{b_2}$ & $8.62e-04 $   & $ 7.36e-04 $  & $7.74e-04 $ & $8.13e-04 $ \\
                                                    & $\epsilon_{p}$   & $ 1.77e-02$    & $2.37e-02$   & $ 2.13e-02$ & $3.54e-02
 $ \\
        \hline
    \end{tabular}
    \caption{The 2D unsteady case: the relative $L_2$ error of the $MHDnet-B$, $PINN-B$, $MHDnet-A_2$ and $PINN-A_2$ when t=0.25, 0.50, 0.75,1.00. }
    \label{tab:2d-1-2}
\end{table}
\cref{tab:2d-1-2} shows the relative $L_2$ error of the velocity field and magnetic field of MHDnet and PINN after training at different times. We observe the relative $L_2$ error remains small although they have a slight changes for the predicted values over time. On the other hand, $MHDnet-A_2$ shows the best performance when the condition of free divergence constraints is fully satisfied. It is worth noting that, for the hidden pressure field, the predictions of PINN significantly deviate from the accurate solutions, whereas MHDnet can directly yield and accurately predict the pressure field without extra data and computational cost. 

\begin{figure}[H]
    \centering
    \subfigure[Exact]{\includegraphics[width=0.33\hsize]{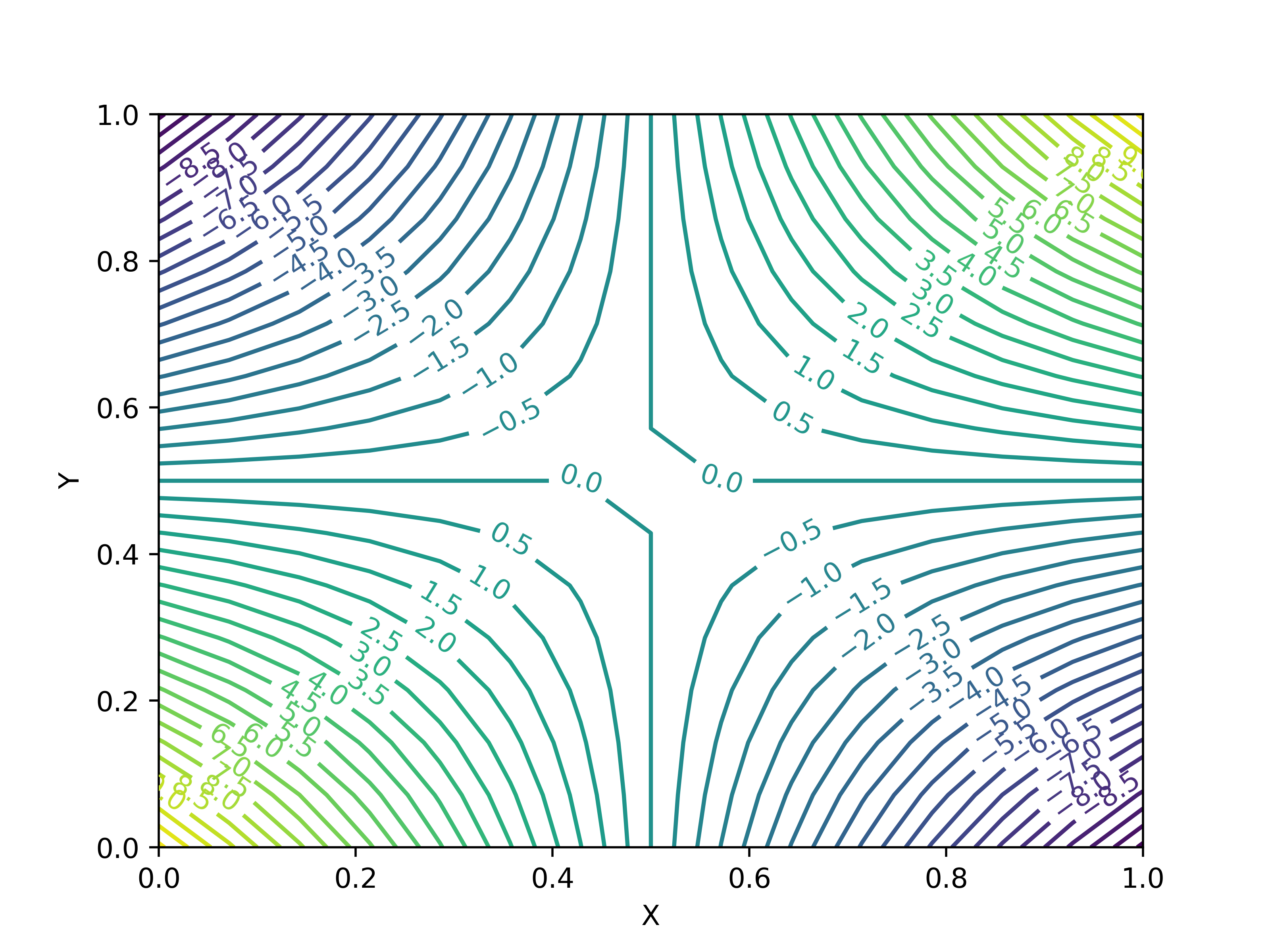}}
    \subfigure[$MHDnet-A_2$]{\includegraphics[width=0.33\hsize]{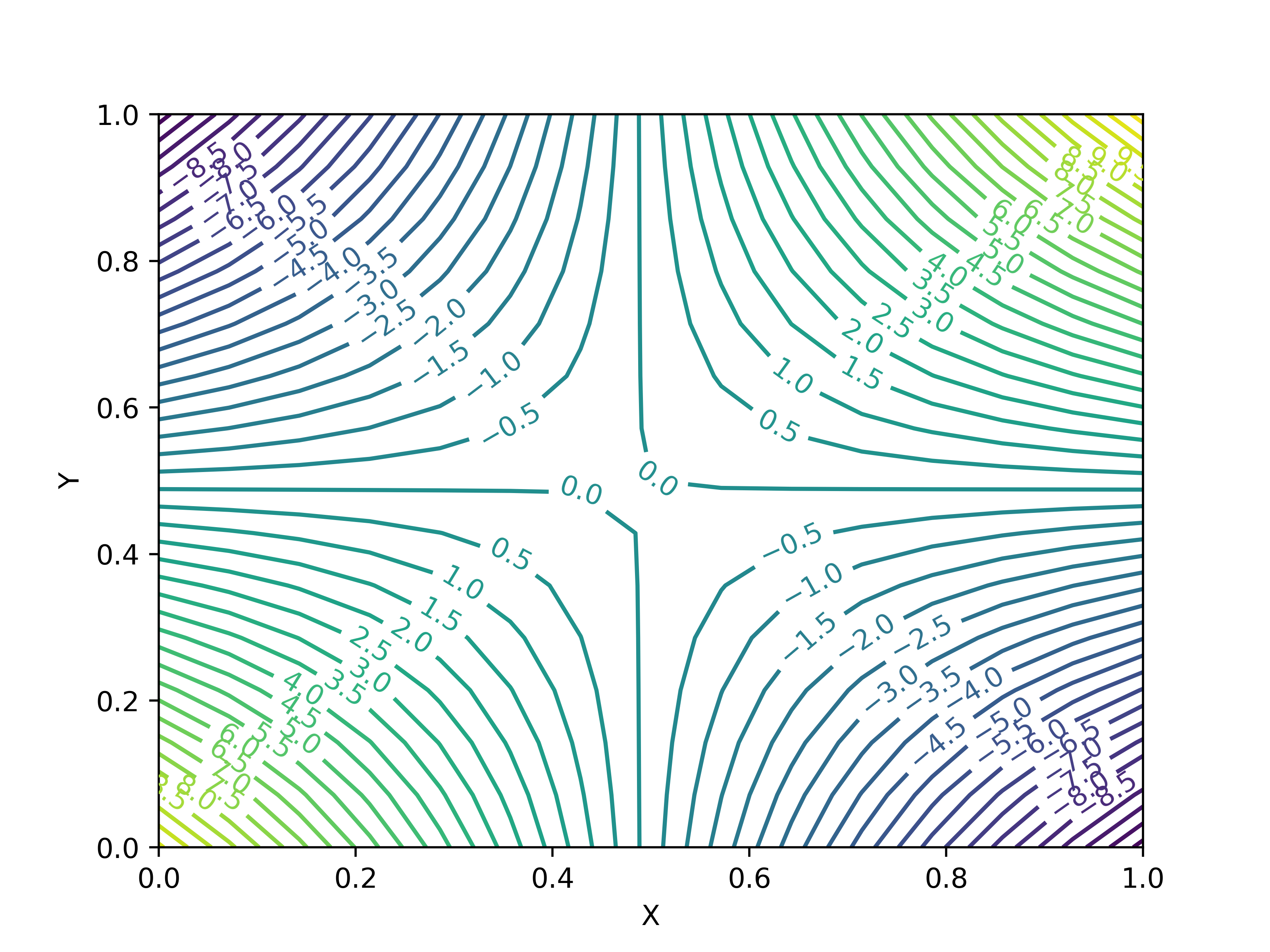}}
    \subfigure[$PINN-A_2$]{\includegraphics[width=0.33\hsize]{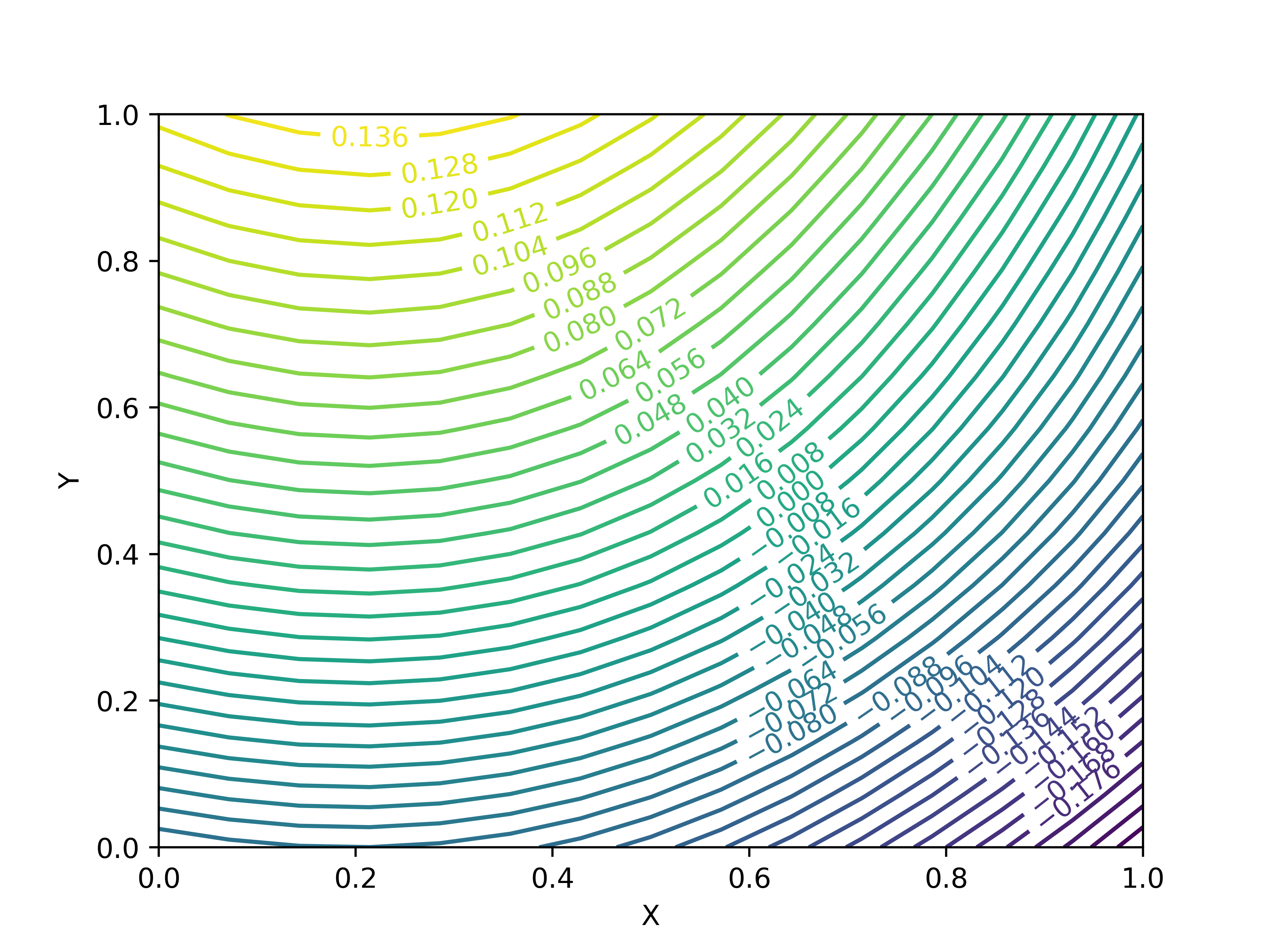}}\\
    \subfigure[Exact]{\includegraphics[width=0.33\hsize]{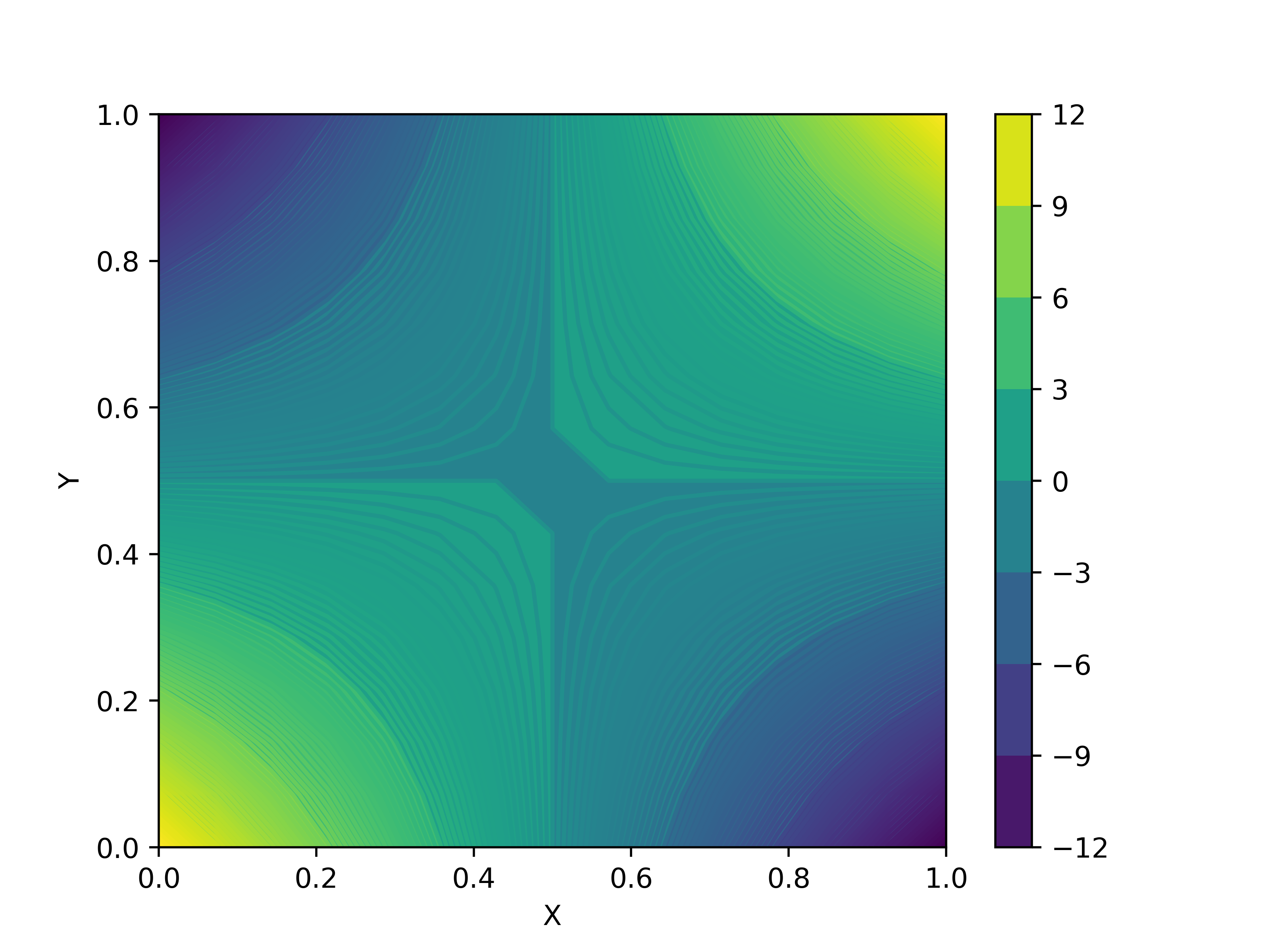}}
    \subfigure[$MHDnet-A_2$]{\includegraphics[width=0.33\hsize]{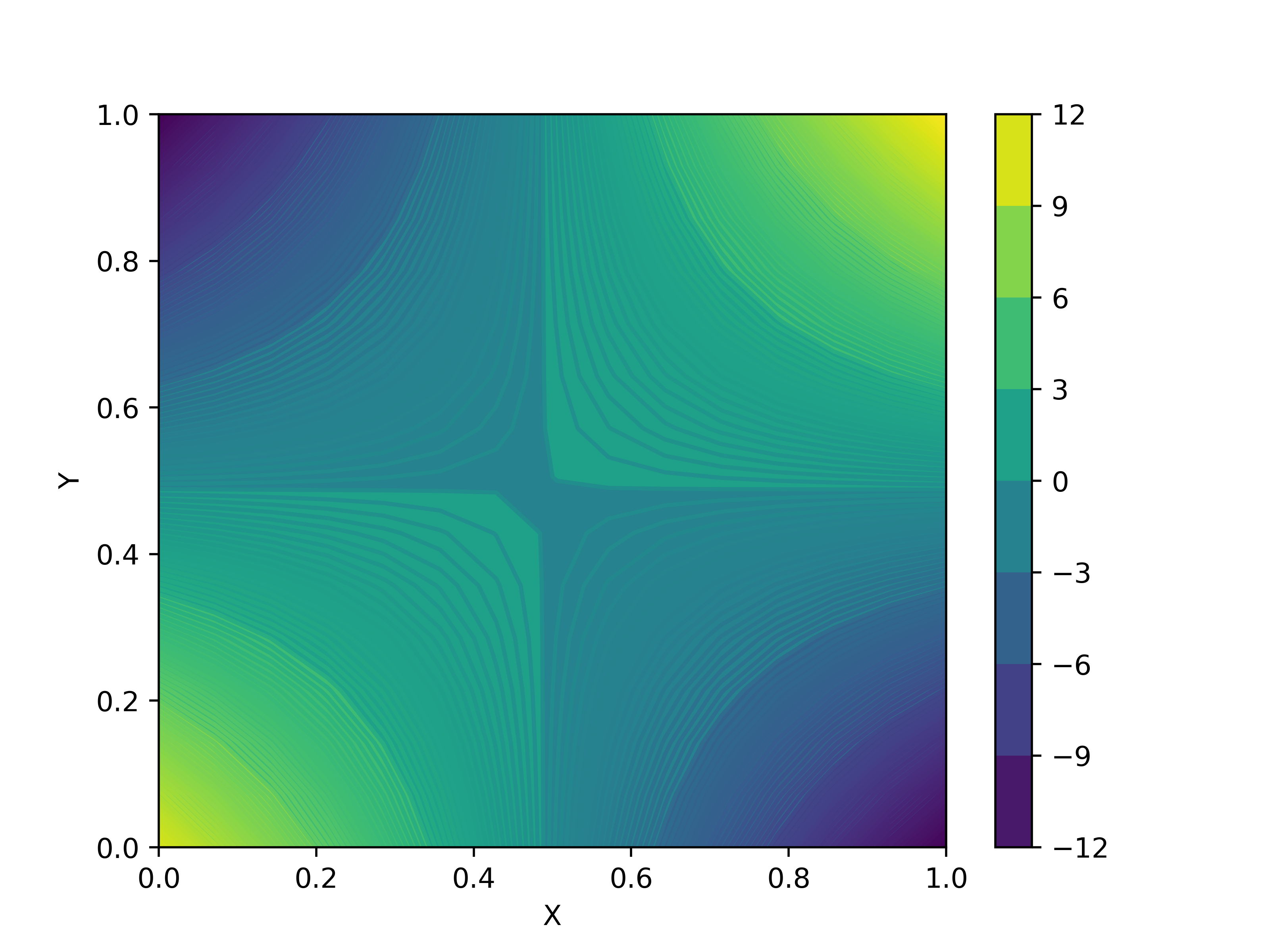}}
    \subfigure[$PINN-A_2$]{\includegraphics[width=0.33\hsize]{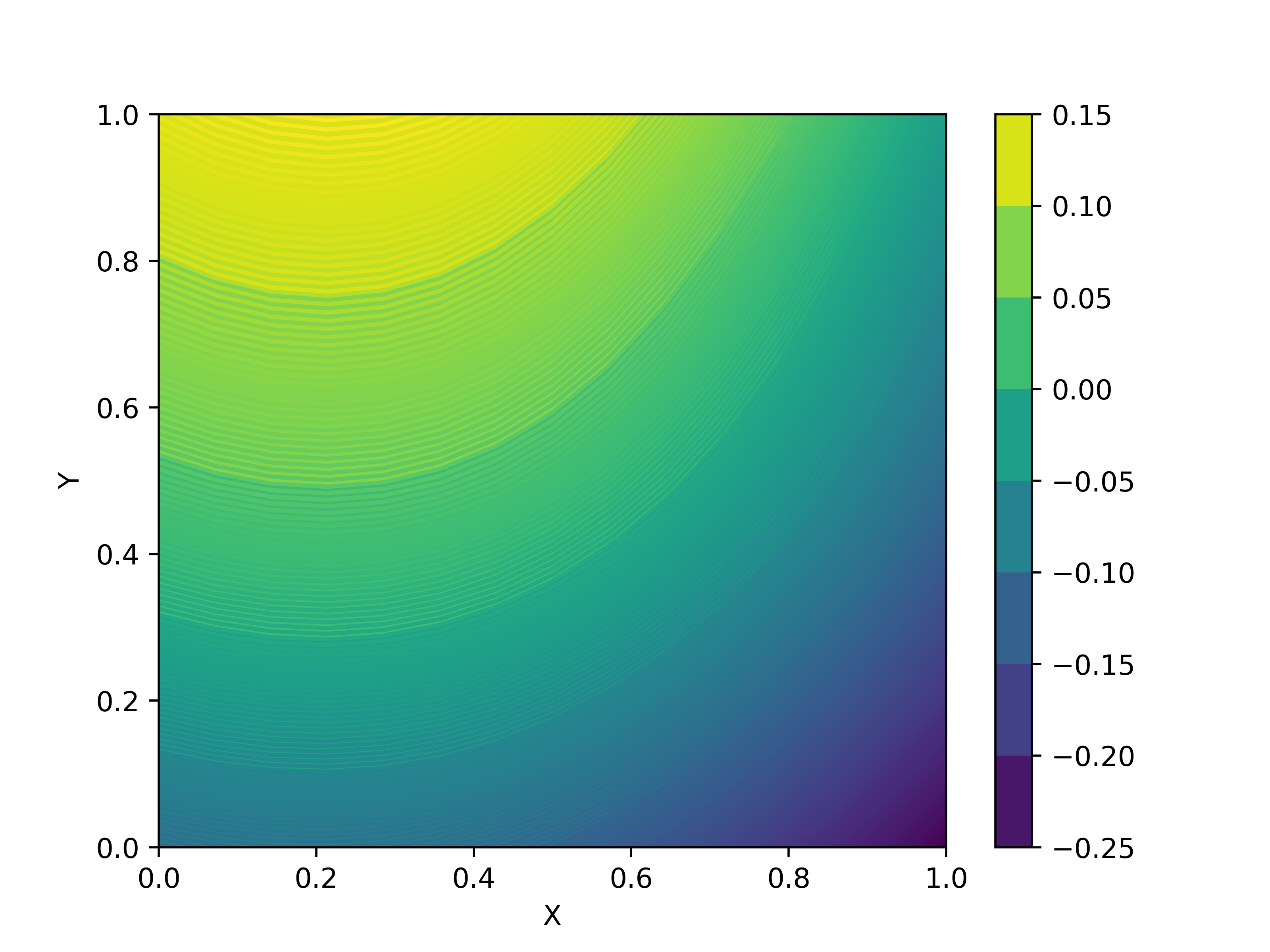}}
    \caption{The 2D unsteady case: pressure distribution obtained by $MHDnet-A_2$ and $PINN-A_2$.}
    \label{fig:2d-1-p}
\end{figure}


\begin{table}[htb]
    \centering
    \begin{tabular}{llllllll}
        \hline                       & $\epsilon_{u_1}$ & $\epsilon_{u_2}$ & $\epsilon_{b_1}$ & $\epsilon_{b_2}$ & $\epsilon_{p}$ \\
        \hline {$\mathrm{PINN-B}$}     & $3.76e-03$      & $4.35e-03$       & $ 5.74e-04$       & $ 6.94e-04$       & $1.22e-01  $   \\
        \hline {$\mathrm{MHDnet-B}$}   & $7.67e-04 $      & $7.69e-04$       & $3.12e-04$       & $2.93e-04$       & $9.91e-03  $   \\
                \hline {$\mathrm{PINN-A_2}$}   & $5.76e-03 $      & $6.91e-03$       & $7.19e-04$       & $9.91e-04$       & $8.60e-02 $    \\
        \hline {$\mathrm{MHDnet-A_2}$} & $5.83e-03 $      & $4.38e-03$       & $3.09e-04$       & $2.57e-04$       & $1.51e-02$ \\    
        \hline
    \end{tabular}
    \caption{2D unsteady case: time-averaged relative $L_2$ error of the $MHDnet-B$, $PINN-B$, $MHDnet-A_2$, $PINN-A_2$}
    \label{tab:2d-1-3}
\end{table}

\begin{figure}[H]
    \centering
    \subfigure[$PINN-B$]{\includegraphics[width=0.45\hsize]{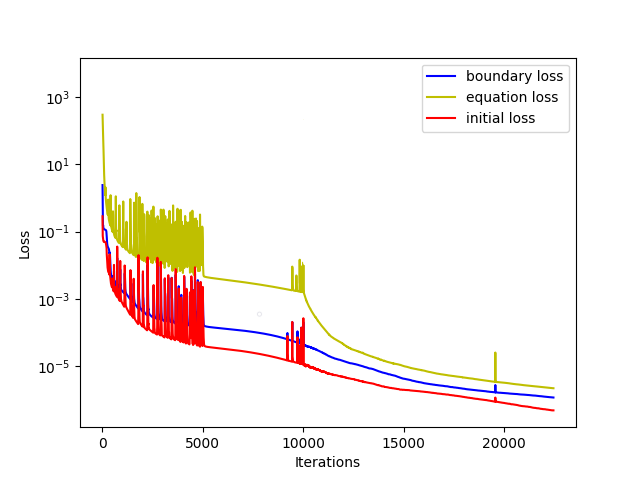}}
    \subfigure[$MHDnet-B$]{\includegraphics[width=0.45\hsize]{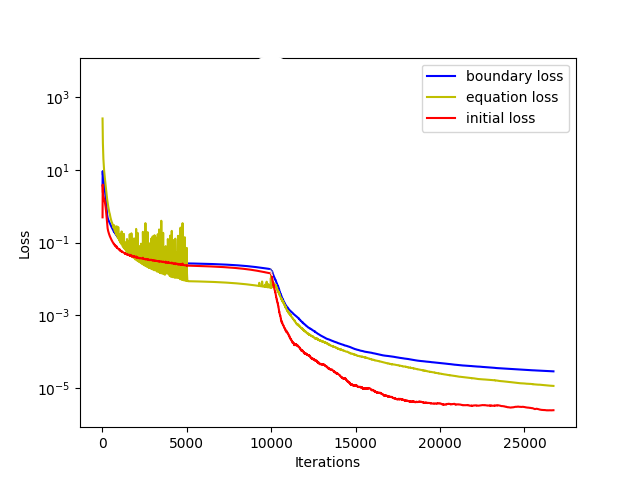}}
    \subfigure[$PINN-A_2$]{\includegraphics[width=0.45\hsize]{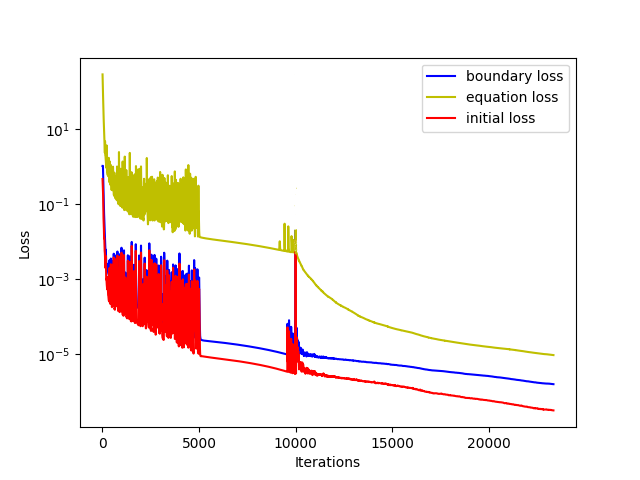}}
    \subfigure[$MHDnet-A_2$]{\includegraphics[width=0.45\hsize]{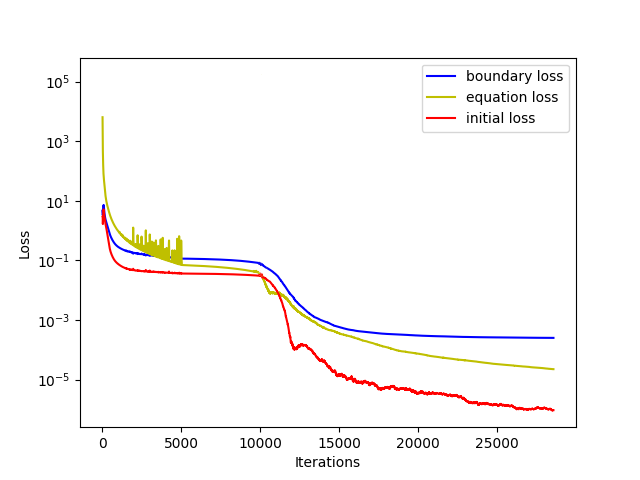}}
    \caption{2D unsteady case: convergence of loss function, equation,initial and boundary loss($L_{\boldsymbol{\hat{f}}},L_{\boldsymbol{\hat{g}}}, L_{\boldsymbol{\hat{h}}}$). }
    \label{fig:2d-loss}
\end{figure}

\cref{fig:2d-1-p} illustrates the pressure distribution obtained by $MHDnet-A_2$ and $PINN-A_2$. The time-averaged relative $L_2$ error of the magnetic field, velocity field, and pressure field can be found in \cref{tab:2d-1-3}. We observe that the relative $L_2$ error between  the velocity field  and the magnetic field  is different in order of magnitude, this is due to the automatic differentiation generated by different derivative orders.
\cref{fig:2d-loss} displays the convergence of loss functions for the four models, which can also be divided into two parts. The oscillating part is obtained by the Adam optimizer, and the steady decline part is obtained by the L-BFGS optimizer. Although all methods achieved a similar convergence at a level of $1e-5$, $MHDnet$ gives higher accuracy solutions, especially for the hidden pressure field.  

\subsection{Hartmann flow}

Hartmann flow is a typical benchmark problem in MHD which describes the internal flow of an incompressible conducting fluid under an external transverse magnetic field $\boldsymbol{B}^d
$. In this example, we focus on using  $A_2~formulation$ . We consider 2D domain $\Omega = [0, L_x]\times[−1, 1], \boldsymbol{B}^d = (0, 1)$, and impose the following boundary conditions:
$$
\begin{array}{ll}
\boldsymbol{u}=\boldsymbol{0} & \text { on } y= \pm 1, \\
\left(p \mathcal{I}-\frac{1}{Re} \nabla \boldsymbol{u}\right) \boldsymbol{n}=p_d \boldsymbol{n} & \text { on } x=0, L_x, \\
\boldsymbol{B} \times \boldsymbol{n}=\boldsymbol{B}^d \times \boldsymbol{n} & \text { on } \partial \Omega .
\end{array}
$$
Then, there is an explicit solution to the MHD equations(see \cite{hartmann}):
$$
\boldsymbol{u}(x, y)=\left(u_1(y), 0\right), \quad p(x, y)=-G x-\frac{s B_1^2(y)}{2}+p_0, \quad \boldsymbol{B}(x, y)=\left(b_1(y), 1\right)
$$
where
$$
\begin{aligned}
& u_1(y)=\frac{G R_e}{H_a \tanh \left(H_a\right)}\left(1-\frac{\cosh \left(y H_a\right)}{\cosh \left(H_a\right)}\right), \quad \\ & b_1(y)=\frac{G}{s}\left(\frac{\sinh \left(y H_a\right)}{\sinh \left(H_a\right)}-y\right), \\
& H_a=\sqrt{s R_e R_m},
\end{aligned}
$$
$H_a$ is Hartmann number, $s =\frac{\boldsymbol{B}^2}{ \mu \rho \boldsymbol{U}^2 }$ denote the coupling number. In this case, we choose $L_x = 4, G = 0.1, p_0 = 0$. and set four sets of parameters (Re = Rm = 1, s = 1; Re = Rm = 20, s = 4; Re = Rm = 50, s = 2; Re = Rm = 50, s = 2).

\begin{figure}[H]
    \centering
    \subfigure[$b_1(y), Re=Rm=s=1$]{\includegraphics[width=0.23\hsize]{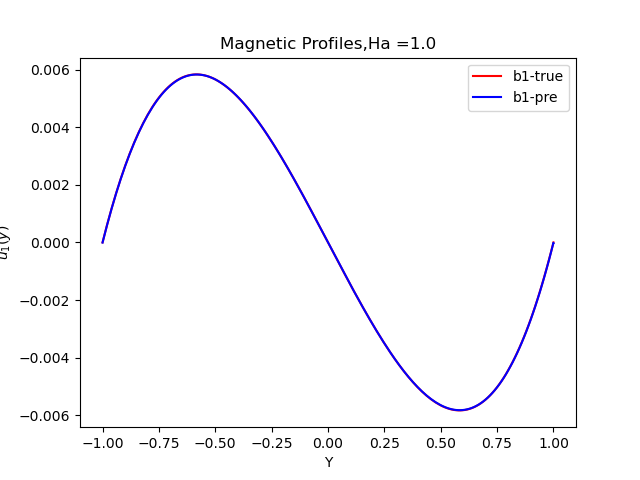}}
    \subfigure[$u_1(y), Re=Rm=s=1$]{\includegraphics[width=0.23\hsize]{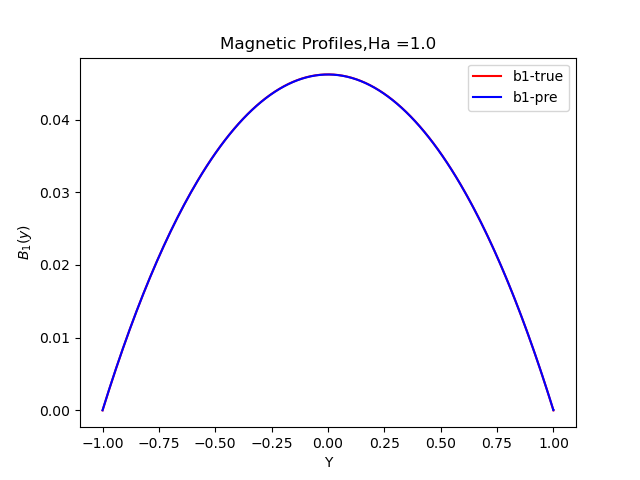}}
    \subfigure[$ b_1(y), Re=Rm=20,s=4$]{\includegraphics[width=0.23\hsize]{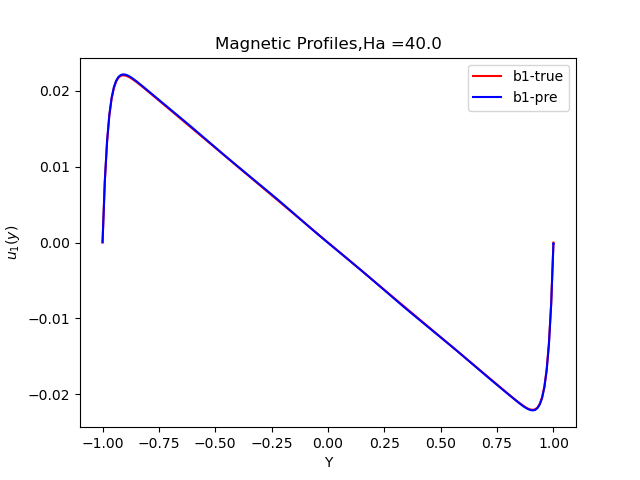}}
    \subfigure[$u_1(y), Re=Rm=20,s=4,$]{\includegraphics[width=0.23\hsize]{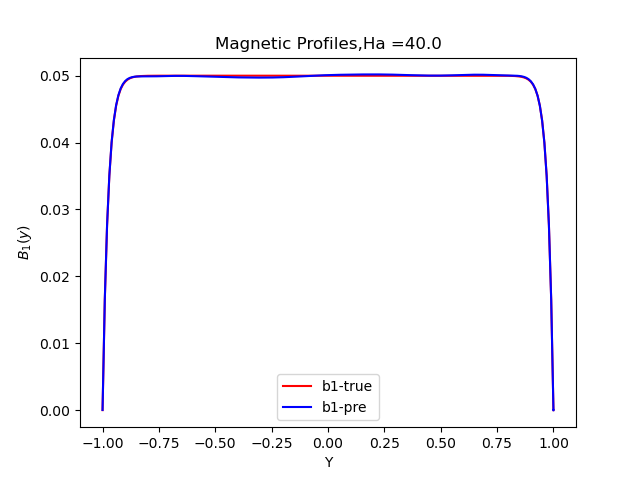}}    \subfigure[$b_1(y), Re=Rm=40,s=2$]{\includegraphics[width=0.23\hsize]{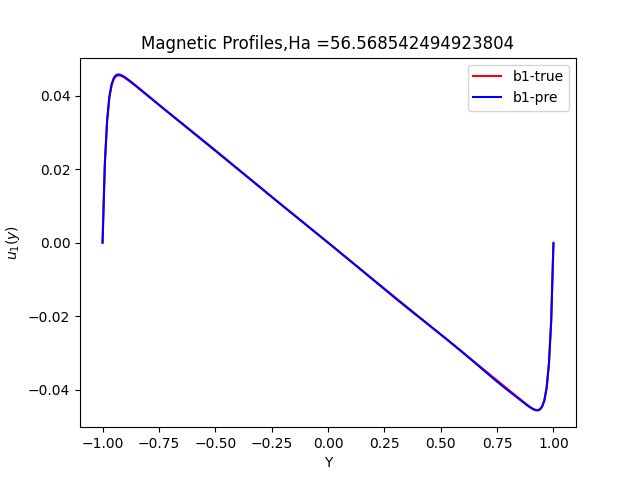}}
    \subfigure[$u_1(y), Re=Rm=40,s=2$]{\includegraphics[width=0.23\hsize]{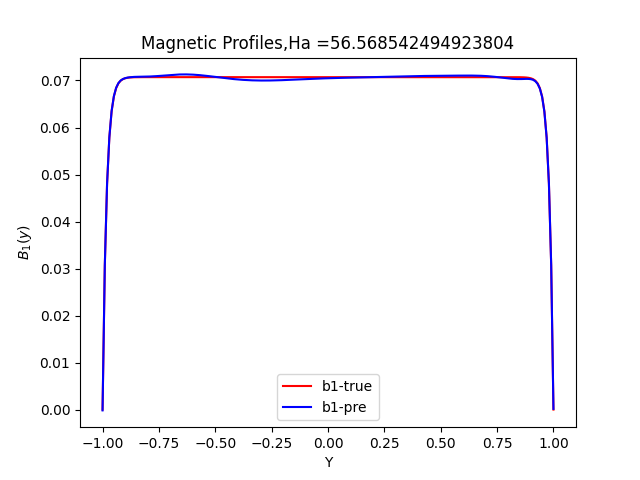}}
    \subfigure[$b_1(y), Re=Rm=50,s=2$]{\includegraphics[width=0.23\hsize]{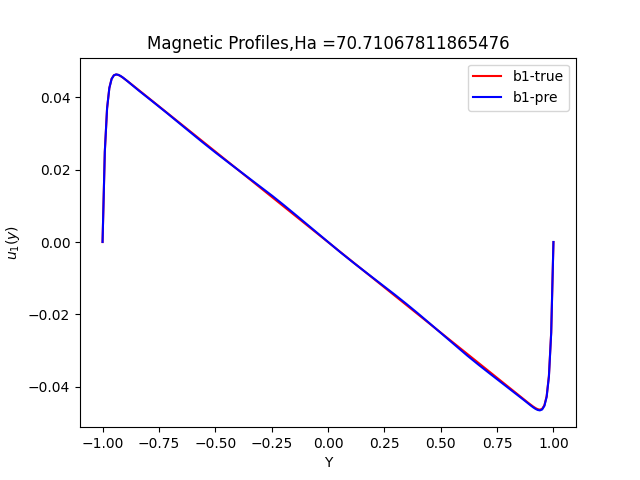}}
    \subfigure[$u_1(y), Re=Rm=50,s=2$]{\includegraphics[width=0.23\hsize]{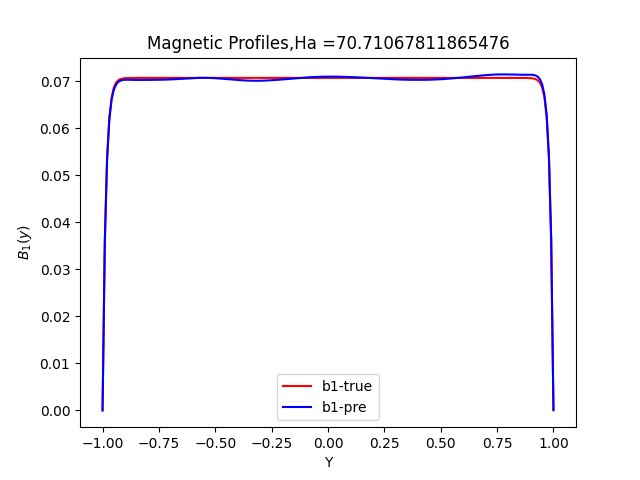}}
    \caption{Hartmann flow: the exact profile(red line) and predicted profile(blue line) of the solution when x=2 under various parameters.}
    \label{fig:2d-3-line}
\end{figure}

The magnetic force near the wall, generated by the electric field, plays a significant role in accelerating the fluid flow. Moreover, a stronger force leads to a reduction in the thickness of the boundary layer.
Fig\ref{fig:2d-3-Re=1}-Fig\ref{fig:2d-3-Re=50} represent distribution of  the true solution, predicted solution, and absolute errors for the Hartmann flow under various parameters. The first row corresponds to the true solution, the second row shows the predicted solution, and the third row displays the absolute error, the details are listed in \ref{appendix A}. \cref{tab:2d-3-1} presents the relative $L_2$ error between the true and predicted solutions obtained using MHDnet in this particular case. The subnetwork has 100 neurons per layer, while other settings remain unchanged. \cref{fig:2d-3-line} illustrates the predicted results of MHDnet at the cross-section $x=2$. The red line represents the actual velocity, while the blue line represents the predicted velocity. The predicted solution still aligns with the exact solution with an increase in the Hartmann number.

\begin{table}[htb]
    \centering
    \begin{tabular}{lllllllll}
        \hline                          & $\epsilon_{u_1}$ &  $\epsilon_{b_1}$ & $\epsilon_{b_2}$ & $\epsilon_{p}$ \\
        \hline {$Re=Rm=1,s=1$}                   & $4.138051e-04$       & $1.303448e-03$        & $6.706587e-06$        & $6.833642e-04 $ \\
        \hline {$Re=Rm=20,s=4$}                   & $3.428867e-03$       & $5.512170e-03$        & $1.885394e-04$        & $4.626782e-03 $ \\
        \hline {$Re=Rm=40,s=2$}                   & $3.972328e-03$       & $6.224121e-03$        & $2.421003e-04$        & $ 3.871101e-03 $ \\
        \hline {$Re=Rm=50,s=2$}                   & $5.405405e-03$       & $8.725427e-03$        & $2.620381e-04$        & $3.851215e-03$ \\
        \hline
    \end{tabular}
    \caption{Hartmann flow: relative $l_2$ error of MHDnet under different parameters}
    \label{tab:2d-3-1}
\end{table}


\subsection{3D unsteady case}
\label{Sect:3Dsteady}
This example is applied to verify the approximation accuracy of the MHDnet on the three-dimensional MHD equations, and the exact solution can be chosen as
\begin{equation}
      \left\{ 
    \begin{array}{llllllll}\vspace{1mm}
        u_1(x, y, z, t)=        -a\left[e^{a x} \sin (a y+d z)+e^{a z} \cos (a x+d y)\right] e^{-d^2 t},                              \\  \vspace{1mm}
        u_2(x, y, z, t)=        -a\left[e^{a y} \sin (a z+d x)+e^{a x} \cos (a y+d z)\right] e^{-d^2 t},                              \\ \vspace{1mm}
        u_3(x, y, z, t)=        -a\left[e^{a z} \sin (a x+d y)+e^{a y} \cos (a z+d x)\right] e^{-d^2 t},                              \\ \vspace{1mm}
        p(x, y, z, t)=        -\frac{1}{2} a^2\left[e^{2 a x}+e^{2 a y}+e^{2 a z}+2 \sin (a x+d y) \cos (a z+d x) e^{a(y+z)}\right. \\ \vspace{1mm}
                             \left. +2 \sin (a y+d z) \cos (a x+d y) e^{a(z+x)}          +2 \sin (a z+d x) \cos (a y+d z) e^{a(x+y)}\right] e^{-2 d^2 t},                                \\ \vspace{1mm}
        b_1(x,y,z,t)       =  sin(z),                                                                                               \\ \vspace{1mm}
        b_2(x,y,z,t)       =  sin(x),                                                                                               \\ \vspace{1mm}
        b_3(x,y,z,t)       =  sin(t+y),                                                                                             \\
    \end{array}
    \right.
\end{equation}
where the domain is defined by $[-1,1]^3$, the time interval is $[0,1]$, $a=d=1,~Re = 40,~\kappa = 1$ and $Rm = 1$. Here 500 interior points, 100 initial points and 400 boundary points are sampled for training MHDnet. 


\begin{table}[htb]
    \centering
    \begin{tabular}{c|cccccc}
        \hline
        z and t  & \text { $\epsilon_{b_1}$ } & \text { $\epsilon_{b_2}$ } & \text { $\epsilon_{b_3}$ } & \text { $\epsilon_{u_1}$ } & \text { $\epsilon_{u_2}$ } & \text { $\epsilon_{u_3}$ } \\
        \hline
        0.25 & $9.96e-03$                 & $3.79e-03 $                & $4.60e-03 $                & $2.52e-02$                 & $1.51e-02$                 & $2.95e-02$                 \\
        0.50 & $4.97e-03$                 & $4.71e-03 $                & $3.98e-03 $                & $2.07e-02$                 & $4.22e-02$                 & $4.03e-02$                 \\
        0.75 & $4.22e-03$                 & $4.52e-03 $                & $4.94e-03 $                & $3.58e-02$                 & $6.95e-02$                 & $2.95e-02$                 \\
        1.0  & $1.00e-02$                 & $1.48e-02 $                & $1.41e-02 $                & $2.55e-02$                 & $6.19e-02$                 & $3.98e-02$                 \\
        temporal mean & $8.77e-03$                 & $6.81e-03 $                & $7.41e-03 $                & $1.90e-02$                 & $1.43e-02$                 & $1.28e-02$                 \\
        \hline
    \end{tabular}
    \caption{3D unsteady case: relative $L_2$ error of $MHDnet-A_2$, when z and t equal 0.25, 0.50, 0.75, 1.00. }
    \label{tab:3d1}
\end{table}

\begin{figure}[H]
    \centering
    \subfigure[$MHDnet-A_2$]{\includegraphics[width=0.38\hsize]{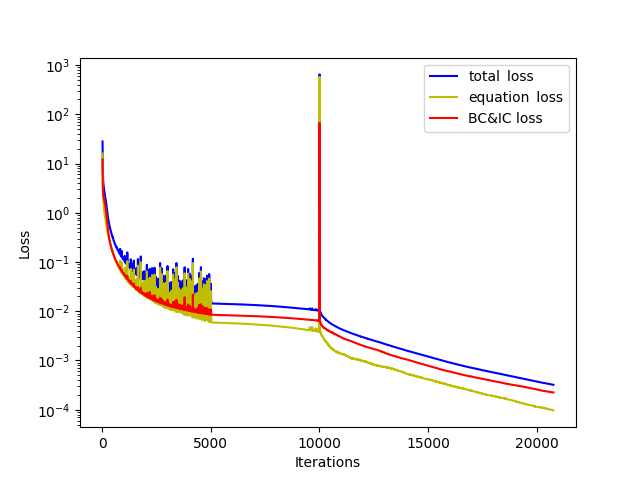}}
    \subfigure[ Vector diagram of magnetic field]{\includegraphics[width=0.3\hsize]{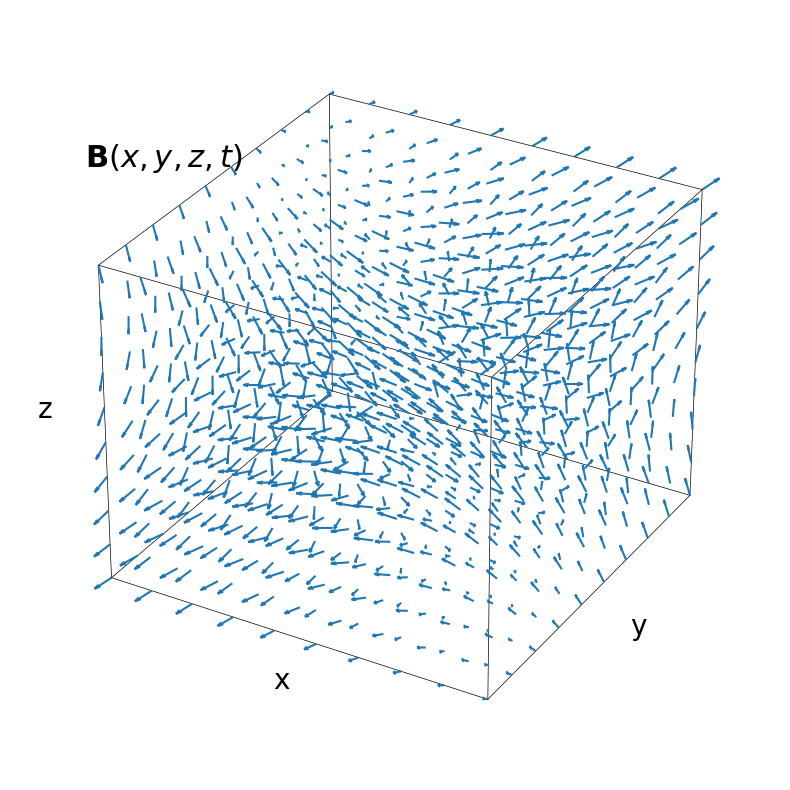}}
    \subfigure[ Vector diagram of velocity field]{\includegraphics[width=0.3\hsize]{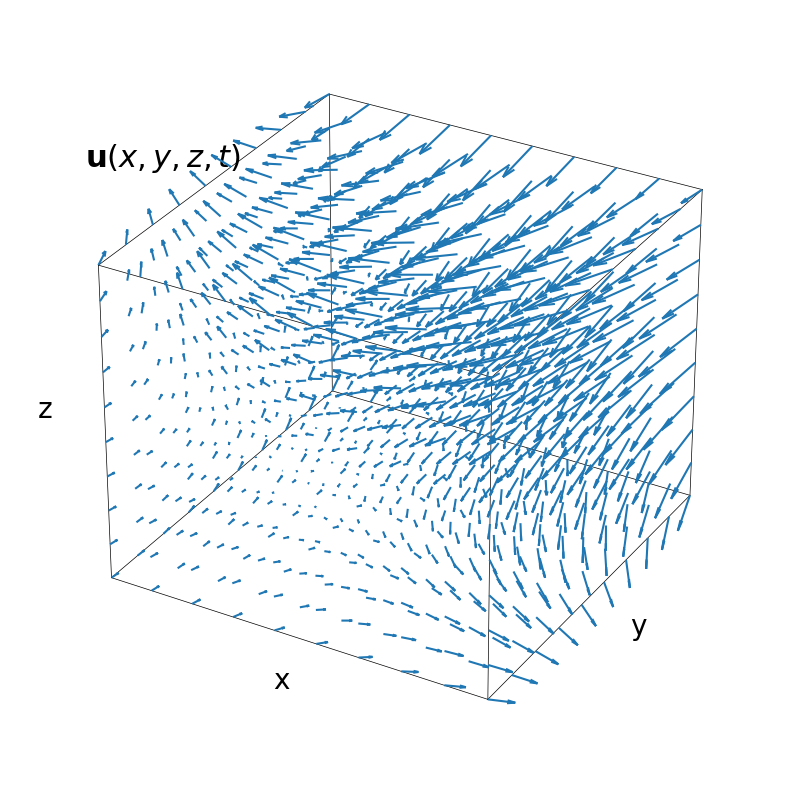}}
    \caption{3D unsteady case: (a) convergence of loss function, (b) magnetic field at t=1, (c) velocity field at t=1}
    \label{fig:3d-loss-vector}
\end{figure}

\cref{tab:3d1} shows the relative $L_2$ error of $MHDnet-A_2$ on different x-y cross-sections. We can observe that $MHDnet-A_2$ demonstrates remarkable performance in prediction of magnetic field, velocity field and pressure field. The loss curve is presented in \cref{fig:3d-loss-vector} (a) with the final value of $3.20e-4$. \cref{fig:3d-loss-vector}(b)$\sim$ (c) displays the vector diagram of the magnetic field and velocity field by $MHDnet-A_2$ at t=1, respectively.

\begin{figure}[H]
    \centering
    \subfigure[Exact ]{
        \begin{minipage}[t]{0.33\linewidth}
            \centering
            \includegraphics[width=1.7in]{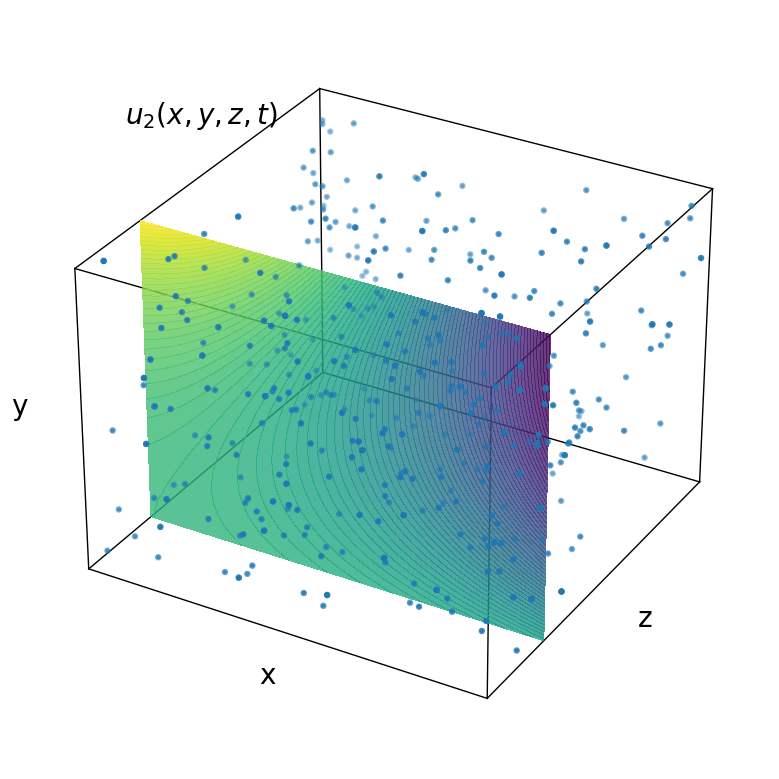}
            \hspace{0.1cm}
            \includegraphics[width=1.7in]{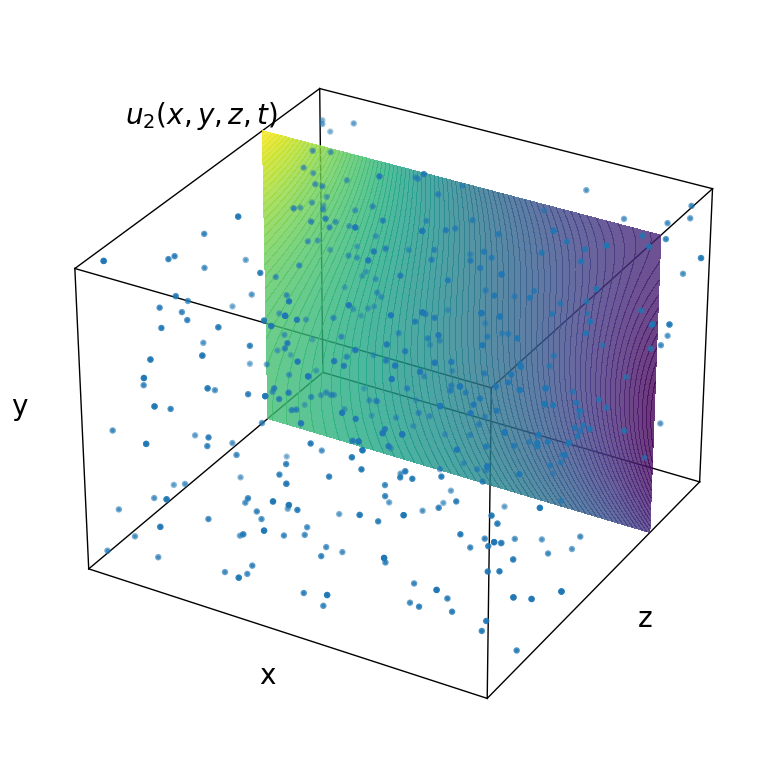}
            \hspace{0.1cm}
            \includegraphics[width=1.7in]{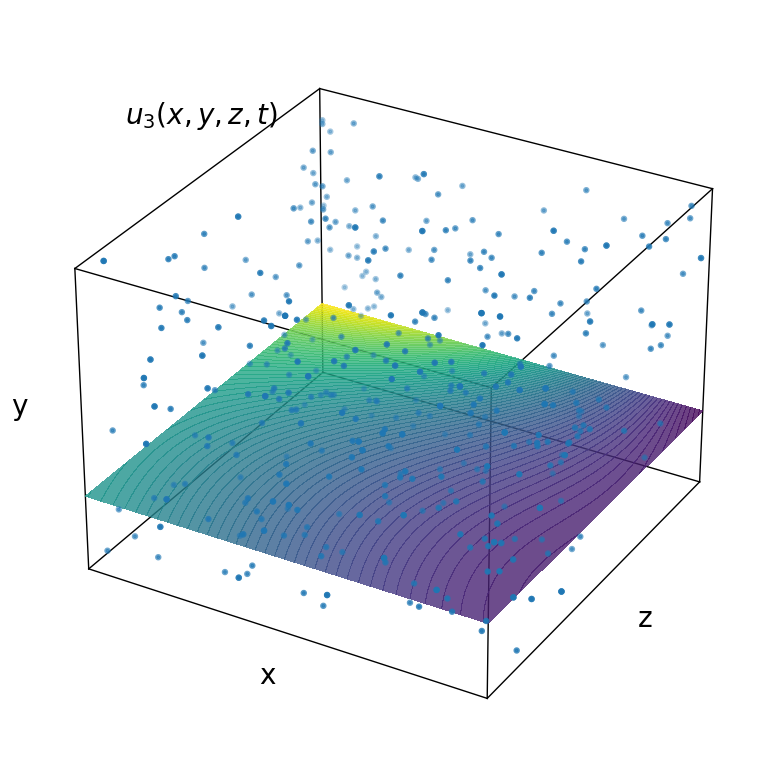}
            \hspace{0.1cm}
            \includegraphics[width=1.7in]{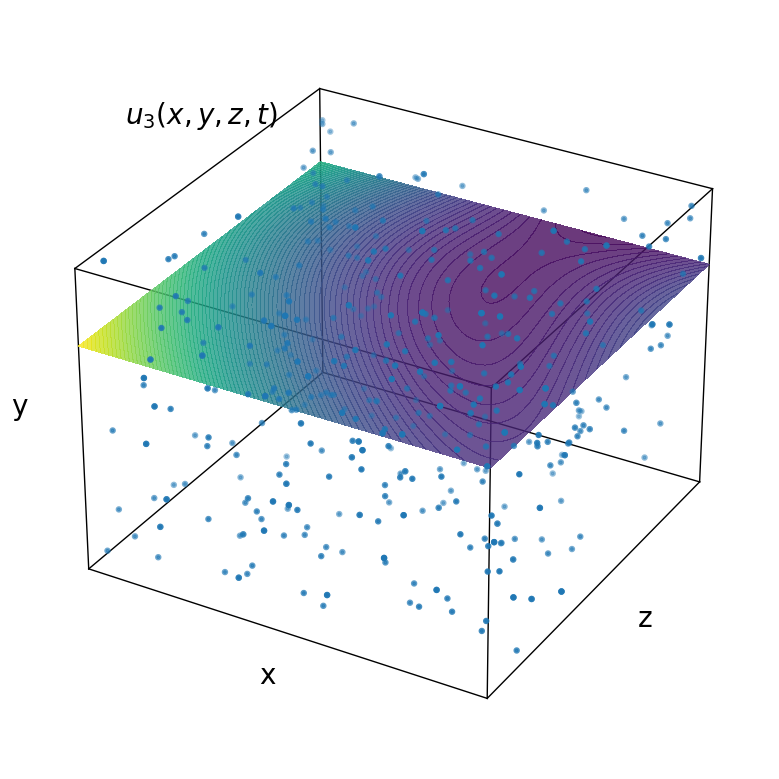}
            \hspace{0.1cm}
        \end{minipage}%
    }%
    \subfigure[MHDnet ]{
        \begin{minipage}[t]{0.33\linewidth}
            \centering
            \includegraphics[width=1.7in]{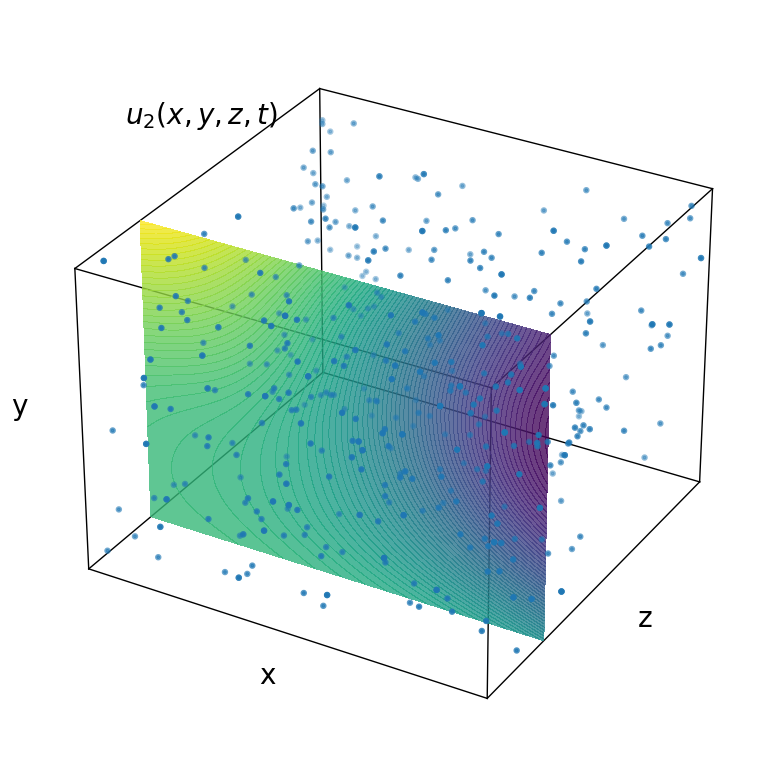}
            \hspace{0.1cm}
            \includegraphics[width=1.7in]{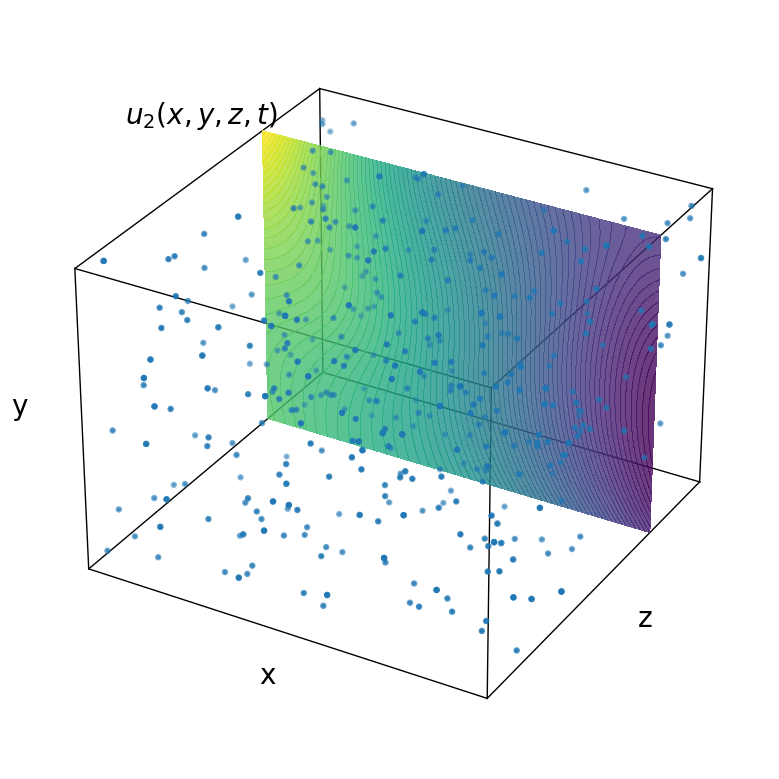}
            \hspace{0.1cm}
            \includegraphics[width=1.7in]{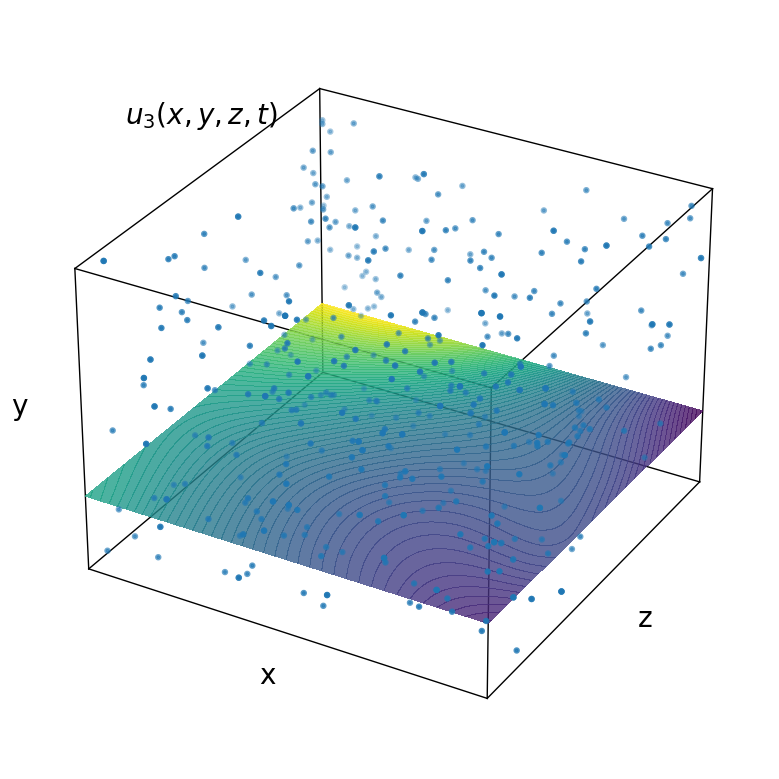}
            \hspace{0.1cm}
            \includegraphics[width=1.7in]{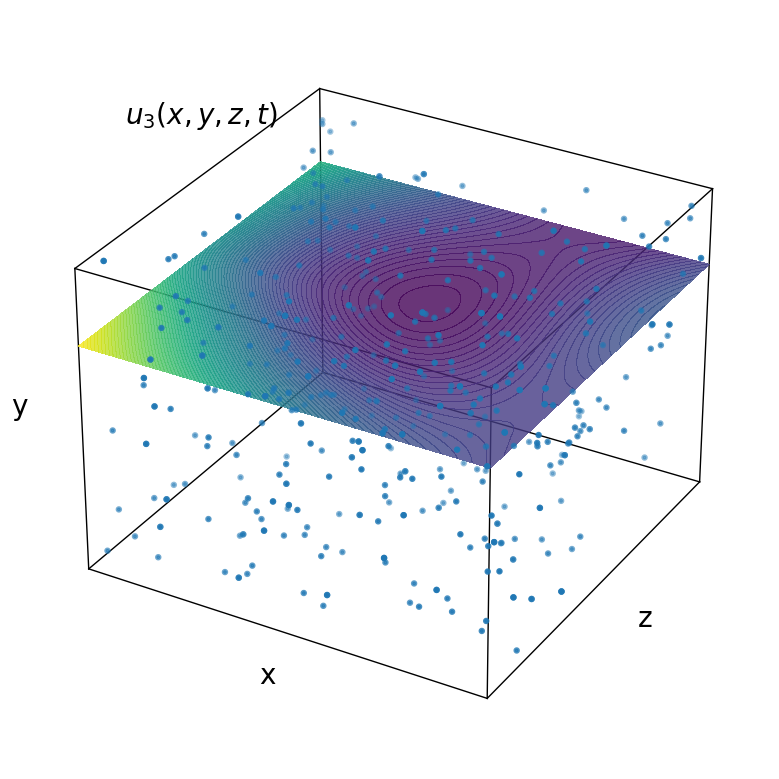}
            \hspace{0.1cm}
        \end{minipage}%
    }%
    \subfigure[Absolute error]{
        \begin{minipage}[t]{0.33\linewidth}
            \centering
            \includegraphics[width=1.7in]{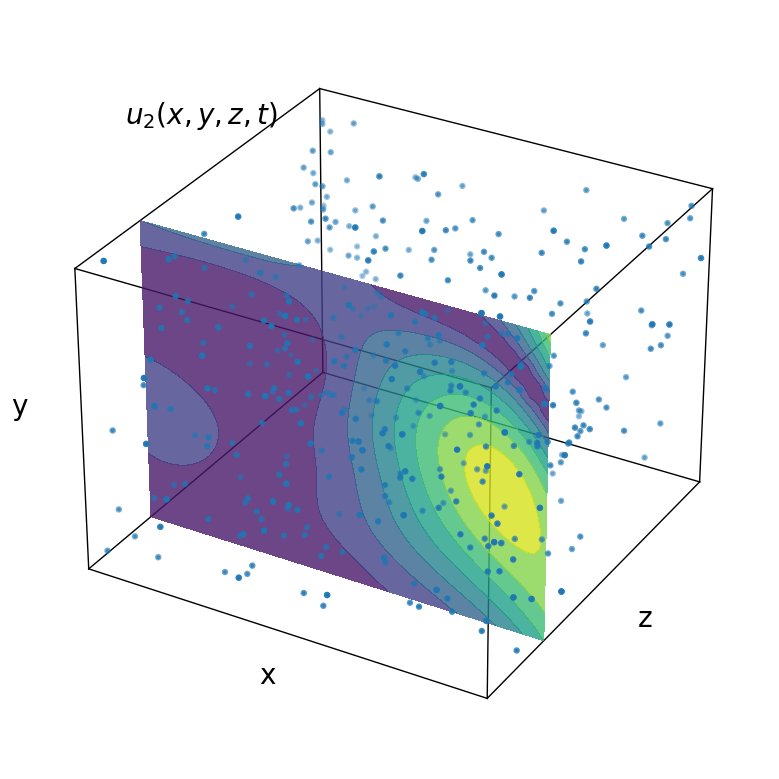}
            \hspace{0.1cm}
            \includegraphics[width=1.7in]{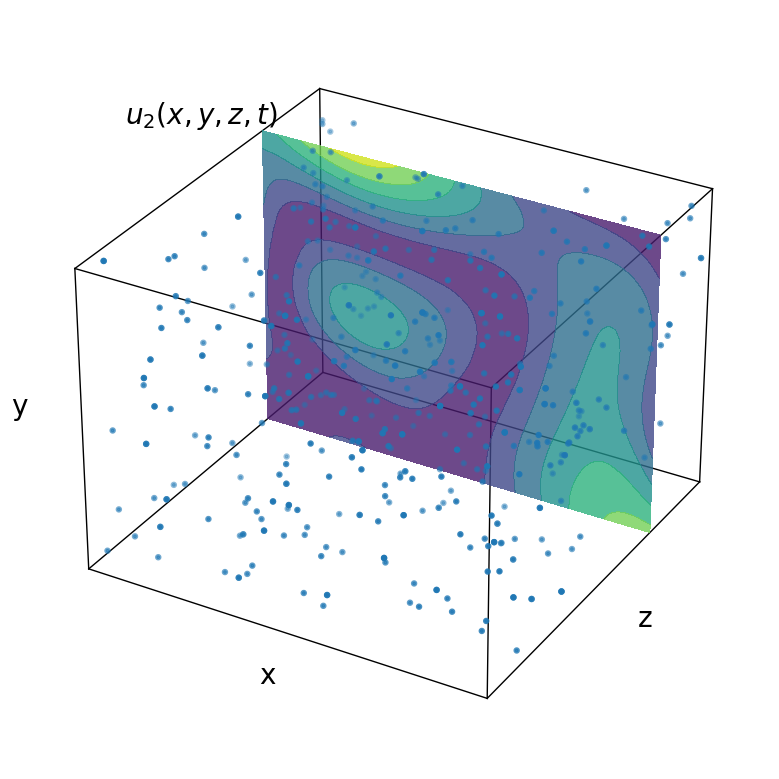}
            \hspace{0.1cm}
            \includegraphics[width=1.7in]{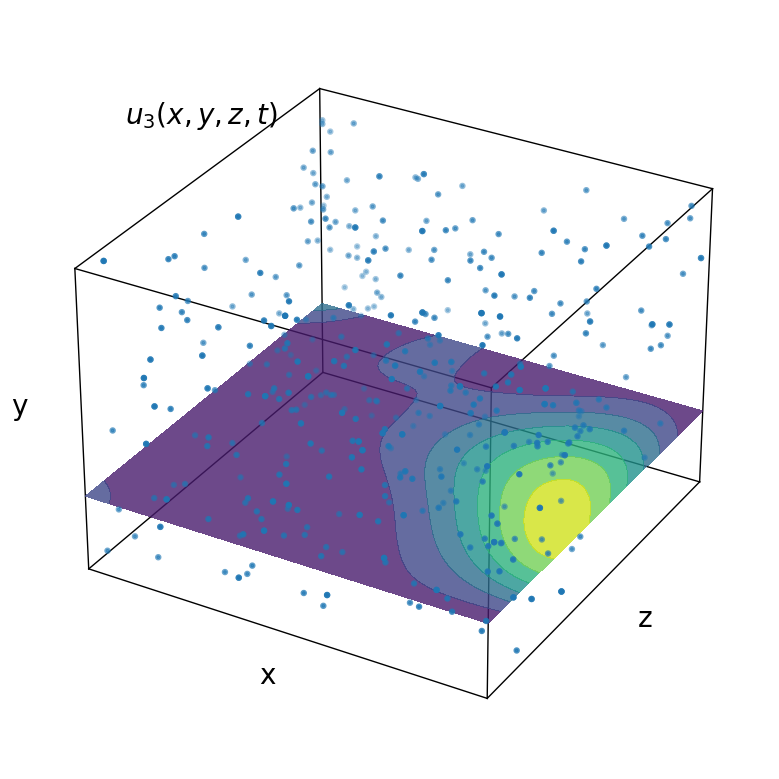}
            \hspace{0.1cm}
            \includegraphics[width=1.7in]{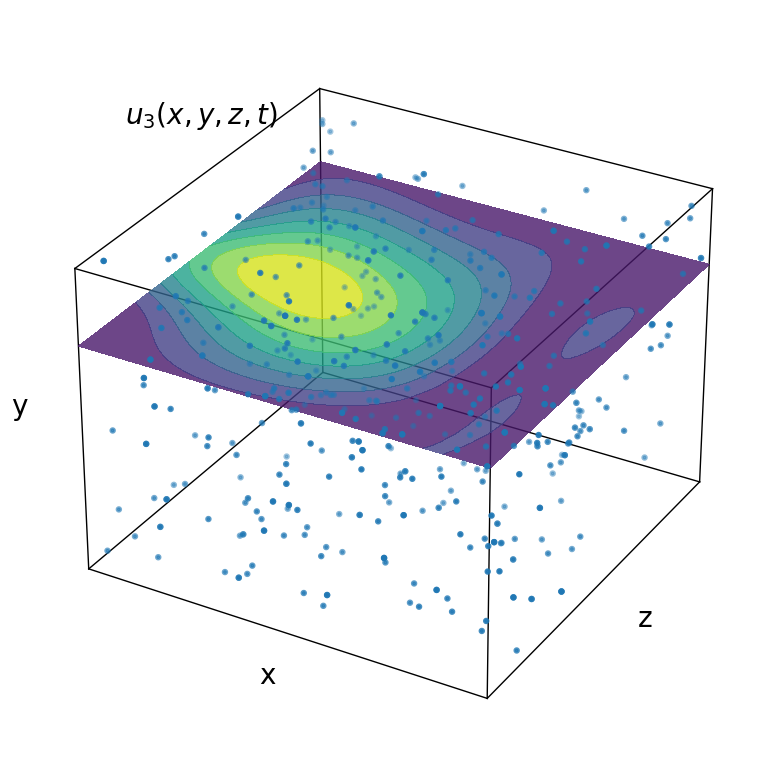}
            \hspace{0.1cm}
        \end{minipage}%
    }%
    \centering
    \caption{3D unsteady case: exact solution, MHDnet solution, and absolute error in different cross-sections for the velocity field.}
    \vspace{-0.2cm}
    \label{fig:3d-1-2}
\end{figure}

\cref{fig:3d-1-2} shows the exact solution, MHDnet solution, and absolute errors of the magnetic field and velocity field. The first two rows display the x-y plane of $z=\pm0.5$ and $t=1.00 $, the last two rows are the x-z plane of $y=\pm0.5$ and $z=1$. Moreover, all of the relative $L_2$ errors are less than $10\%$, which indicates MHDnet can successfully capture the relatively complex distribution of the velocity field. It is worth stating that adaptive sampling methods could be used to achieve even better results. However, due to the focus and scope of this paper, we did not delve into the details of it. Nonetheless, the results have shown that the MHDnet can effectively predict the physical fields on different cross-sections, and a promising approach is provided in studying and analyzing MHD systems.

\subsection{2D steady case with missing or noisy boundary conditions}
\label{Sect:2Dsteady1}
Because of its powerful fitting ability and convenient implementation, MHDnet can easily overcome the solving difficulties of traditional numerical methods with missing or noisy boundary conditions, which are obviously ill-conditioned. 2D steady MHD equations are chosen as model experiments, and all hyperparameter settings, training optimization strategies, and data sets are also chosen as same as those of the 2D steady case in \cref{Sect:2Dsteady}.


\begin{figure}[H]
    \centering
    \subfigure[standard]{\includegraphics[width=0.23\hsize]{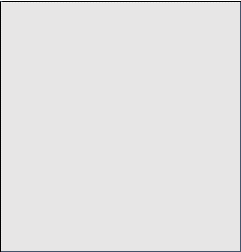}}
    \subfigure[right]{\includegraphics[width=0.23\hsize]{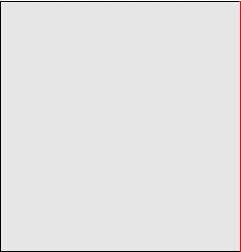}}
    \subfigure[upper right]{\includegraphics[width=0.23\hsize]{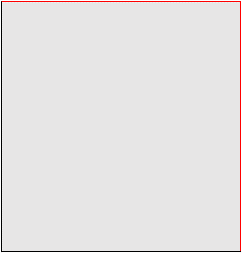}}
    \subfigure[upper down]{\includegraphics[width=0.23\hsize]{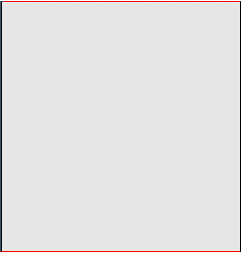}}
    \subfigure[middle]{\includegraphics[width=0.23\hsize]{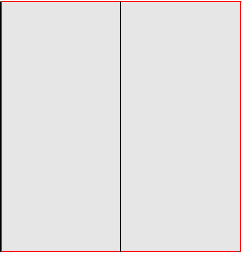}}
    \subfigure[noisy]{\includegraphics[width=0.23\hsize]{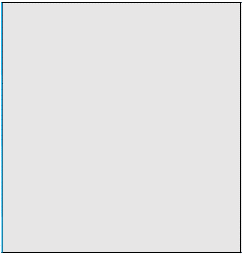}}
    \subfigure[middle noisy]{\includegraphics[width=0.23\hsize]{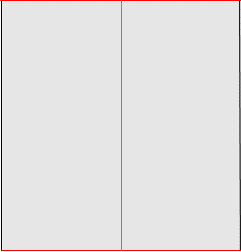}}
    \subfigure[stagger]{\includegraphics[width=0.23\hsize]{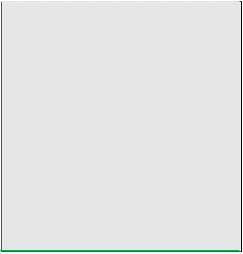}}
    \caption{Schematic diagram of missing and noisy boundary}
    \label{3d-plane}
\end{figure}

\begin{figure}[H]
    \centering
    \subfigure[standard]{\includegraphics[width=0.23\hsize]{new/2d-2-MHD-nodiv-loss.png}}
    \subfigure[right]{\includegraphics[width=0.23\hsize]{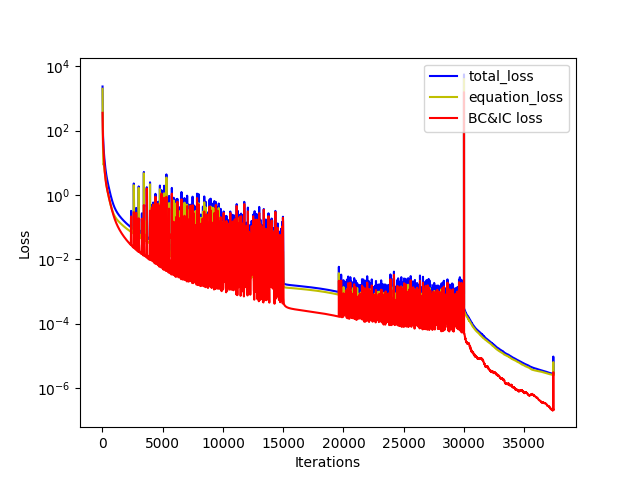}}
    \subfigure[upper right]{\includegraphics[width=0.23\hsize]{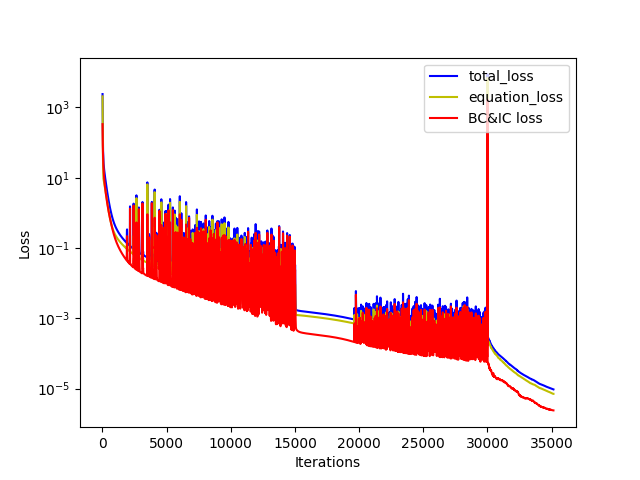}}
    \subfigure[upper down]{\includegraphics[width=0.23\hsize]{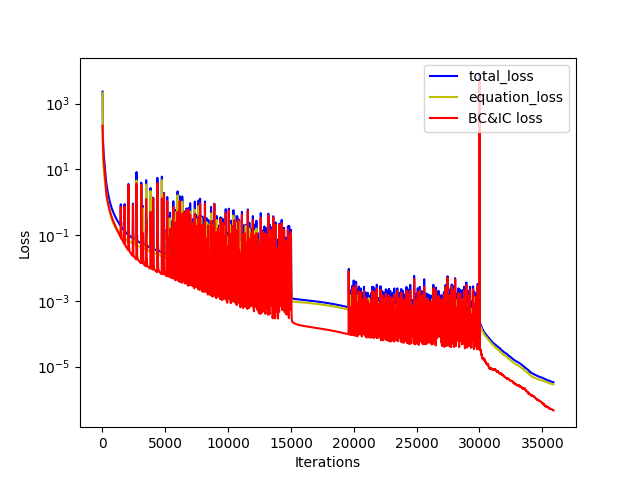}}
    \subfigure[noisy]{\includegraphics[width=0.23\hsize]{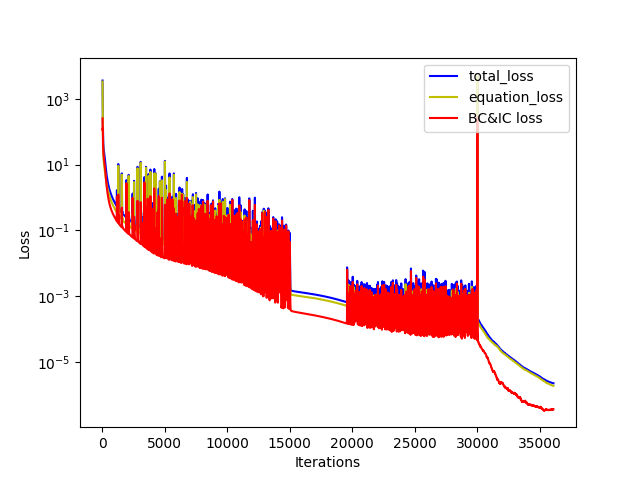}}
    \subfigure[middle noisy]{\includegraphics[width=0.23\hsize]{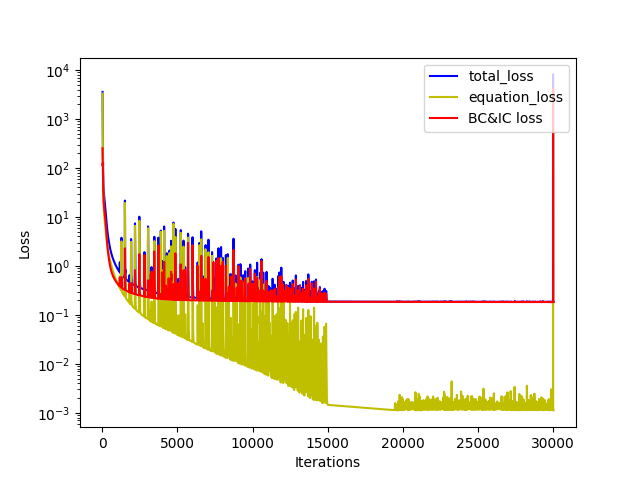}}
    \subfigure[noisy]{\includegraphics[width=0.23\hsize]{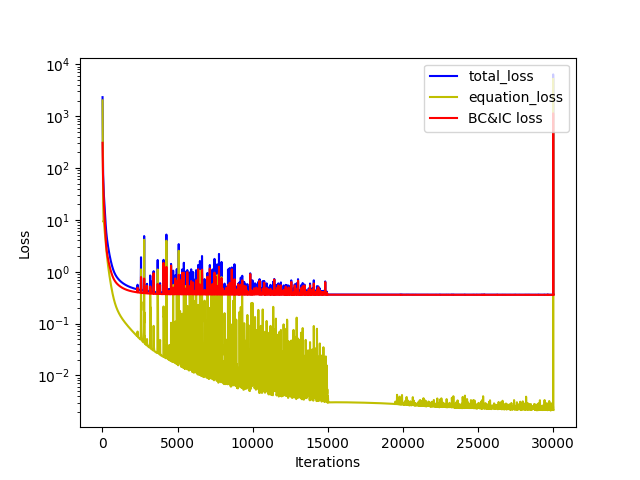}}
    \subfigure[stagger]{\includegraphics[width=0.23\hsize]{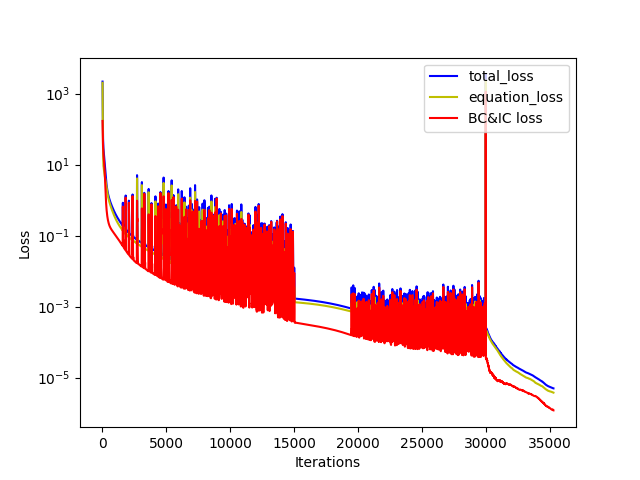}}
    \caption{Convergence of loss function under different boundary conditions}
    \label{fig:2d-missloss}
\end{figure}

\cref{3d-plane} provides a schematic diagram of settings for missing or noisy boundary conditions in this example. The black lines denote the known boundary condition, the red lines are the unknown, and the blue lines represent the noisy boundary condition, where the noise amplitude is $10\%$ of the standard data deviation. The 'stagger' in the bottom right corner refers to the missing staggered boundary, that is, the upper boundary data of the magnetic field is missing, while the lower boundary data of the velocity field is missing. 

\cref{fig:2d-missloss} shows the convergence of the loss functions under different boundary conditions. In the absence of the boundary condition, the convergence result is not as good as before. The data containing noise also has a significant effect on the training process. In \cref{fig:2d-missloss}(f),(g), it can be observed that the boundary loss is the main component of the total loss when stopping the training.
\begin{figure}[H]
    \centering
    \subfigure[standard ]{\includegraphics[width=0.23\hsize]{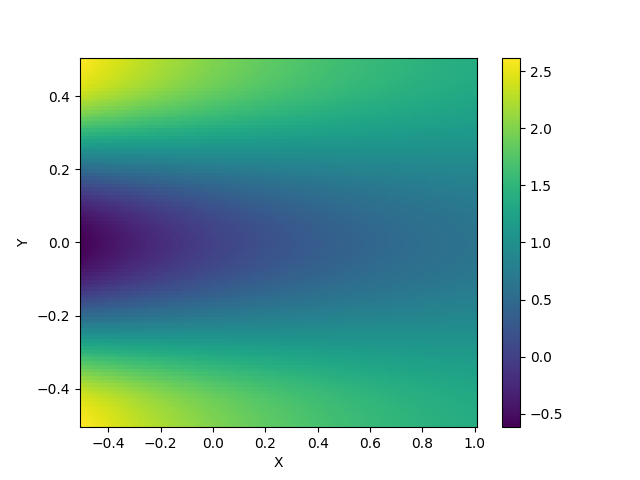}}
    \subfigure[right ]{\includegraphics[width=0.23\hsize]{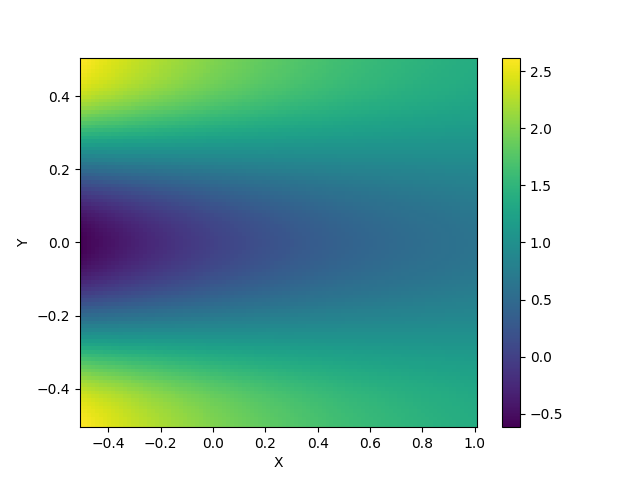}}
    \subfigure[upperright ]{\includegraphics[width=0.23\hsize]{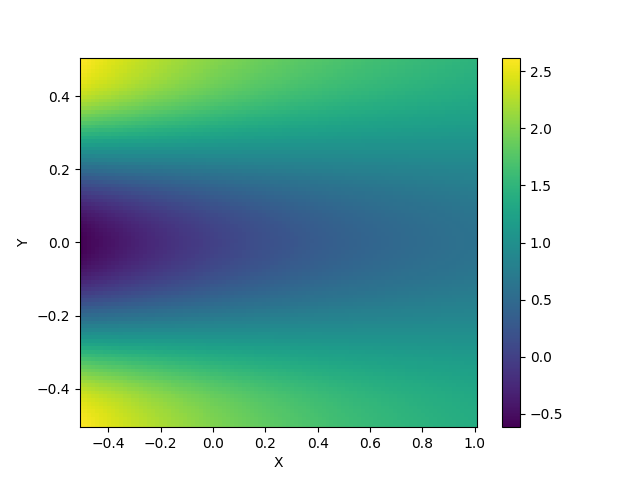}}
    \subfigure[upperdown ]{\includegraphics[width=0.23\hsize]{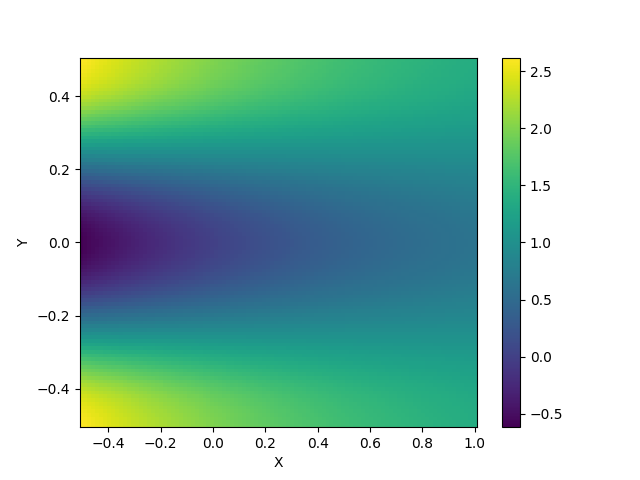}}
    \subfigure[noisy ]{\includegraphics[width=0.23\hsize]{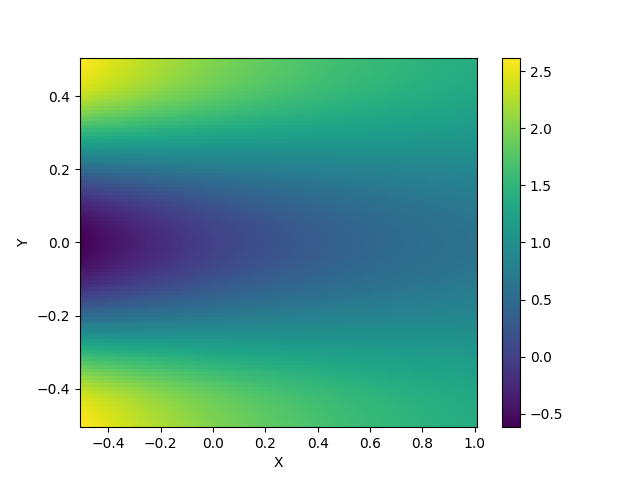}}
    \subfigure[middle noisy ]{\includegraphics[width=0.23\hsize]{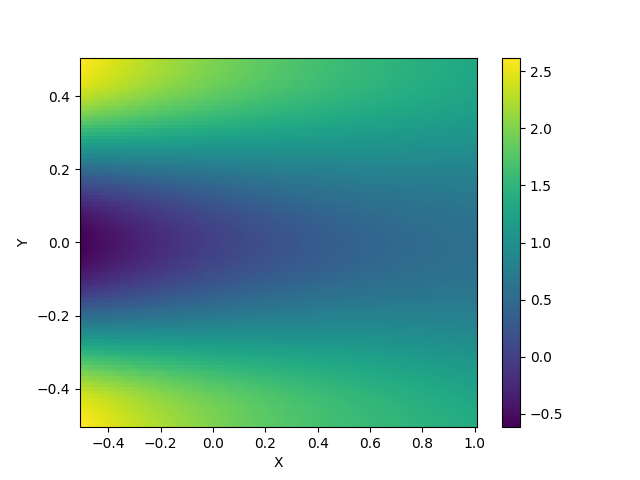}}
    \subfigure[noisy ]{\includegraphics[width=0.23\hsize]{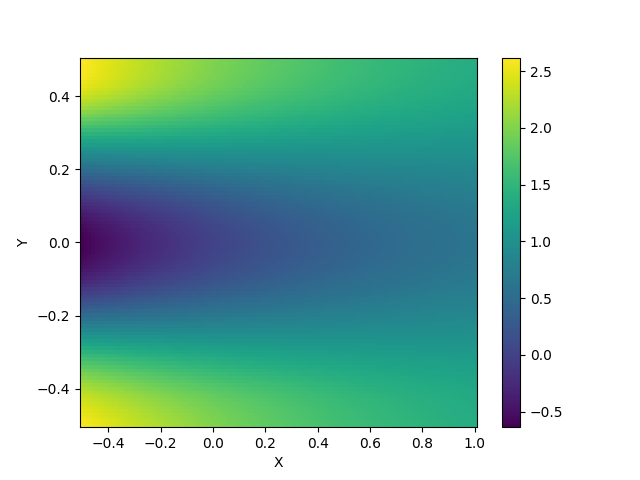}}
    \subfigure[stagger ]{\includegraphics[width=0.23\hsize]{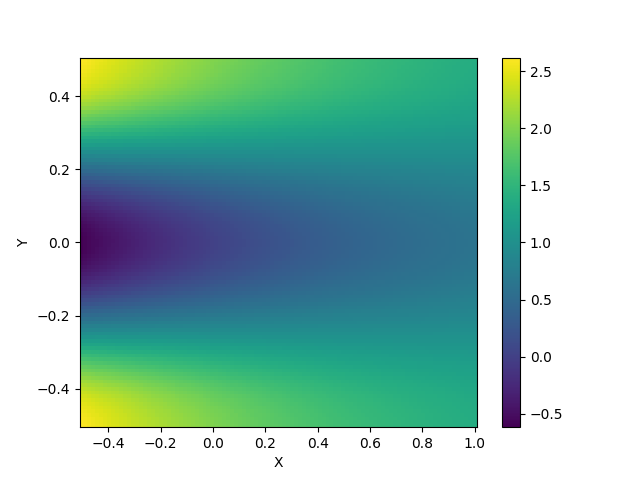}}
    \caption{The predicted solutions of horizontal velocity components $u_1$ under different boundary conditions}
    \label{fig:2d-missu1}
\end{figure}

\begin{figure}[H]
    \includegraphics[width=0.99\textwidth]{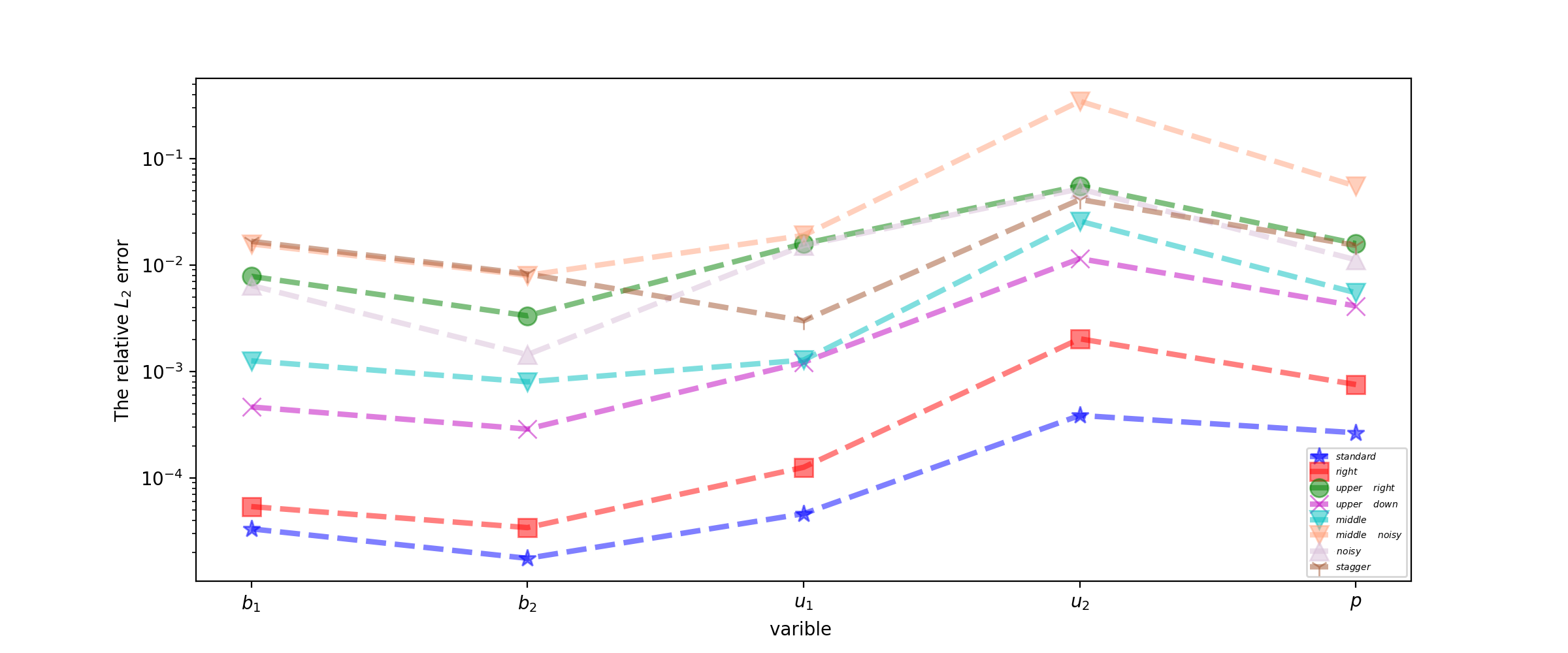}
    \caption{Relative $L_2$ error of $MHDnet-A_2$ with missing or noisy boundaries.    }
    \label{fig:2d-2bd}
\end{figure}


\cref{fig:2d-missu1} shows the horizontal velocity distribution predicted by MHDnet, and \cref{fig:2d-2bd} shows the relative $L_2$ errors of MHDnet under different boundary conditions. Standard (see \cref{3d-plane} (a)) boundary conditions achieve the highest accuracy, but MHDnet also obtains good results under most cases of missing boundary conditions. As for the noisy boundary conditions, MHDnet will report a slightly higher error in the vertical direction of velocity when the upper and lower boundary data are missing at the same time (see \cref{3d-plane} (f)). when random noise is added to the inlet boundary (see \cref{3d-plane} (g)), the relative $L_2$ error of MHDnet is around $1e-3$. This shows that MHDnet is highly reliable and robust in handling various boundary conditions. Moreover, traditional numerical methods may fail or struggle in such scenarios, but the MHDnet can effectively handle these issues. Its ability to work with incomplete or faulty data adds to its utility and makes it an attractive option for scientific research.

\section{Conclusion}\label{sec:con}
Magnetohydrodynamics system has been extensively applied to practical engineering and everyday life. The real-time prediction of the magnetic, velocity, and pressure field of the MHD system for the principled manipulation of MHD is not trivial because of their complex nonlinear, multi-modes feature. So far, researches in this area have been mainly based on experimental tests and traditional numerical simulation. Therefore, the development of reliable and efficient deep learning methods is of great importance for engineering practice. 

In this paper, a novel MHDnet is developed for efficiently simulating complex coupling multiphysics behaviors of MHD system with the multi-modes feature embedding into multiscale neural network architecture. The contributions of this paper are threefold. First, MHDnet can efficiently balance the coupling loss terms of MHD equations with the high-accuracy approximation. Second, the MHDnet naturally preserves the important physical properties  $\operatorname{div}\boldsymbol{B}=0$ and $\operatorname{div}\boldsymbol{u}=0 $, which can accelerate the convergence of the NN by alleviating the interaction of magnetic fluid coupling across different frequency modes. Third, this method can be robust for MHD problems with missing or noisy data, and the pressure fields as the hidden state can be obtained without extra data and computational cost with high accuracy. Numerical experiments show that the proposed MHDnet is stable and effective, and can accurately capture the nonlinear and multiscale information and provide enough numerical accuracy. Despite some simplifying hypotheses adopted in this study, the framework of MHDnet presented here is quite general and can be easily extended to other complex coupled physical problems. 

\section*{Acknowledgement}
This work was supported by the National Natural Science Foundation of China (Nos. 12271409, 12271514, and 12161141017), the Natural Science Foundation of Shanghai (No. 21ZR1465800), the Science Challenge Project (No. TZ2018001), the Interdisciplinary Project in Ocean Research of Tongji University and the Fundamental Research Funds for the Central Universities.
\section*{CRediT authorship contribution statement}
\textbf{Xiaofei Guan:} Conceptualization, Formal analysis, Investigation, Methodology, Validation, Writing - Original Draft, Writing - Review $\&$ Editing.
\textbf{Boya Hu:} Investigation, Methodology, Software, Validation, Visualization, Writing - Original Draft.
\textbf{Shipeng Mao:} Conceptualization, Formal analysis, Investigation, Methodology, Validation, Writing - Review $\&$ Editing.
\textbf{Xintong Wang:} Software, Visualization.
\textbf{Zihao Yang:} Methodology, Software, Writing - Review $\&$ Editing.

\appendix

\subsection{The Details of Hartmann Flow}
\label{appendix A}

The following figure shows the distribution of the true solution, predicted solution, and absolute errors for the Hartmann flow
under various parameters. The first row corresponds to the true solution, the second row shows the predicted
solution, and the third row displays the absolute error. 

\begin{figure}[H]
    \centering
    \subfigure[$B_1(y)$]{\includegraphics[width=0.25\hsize]{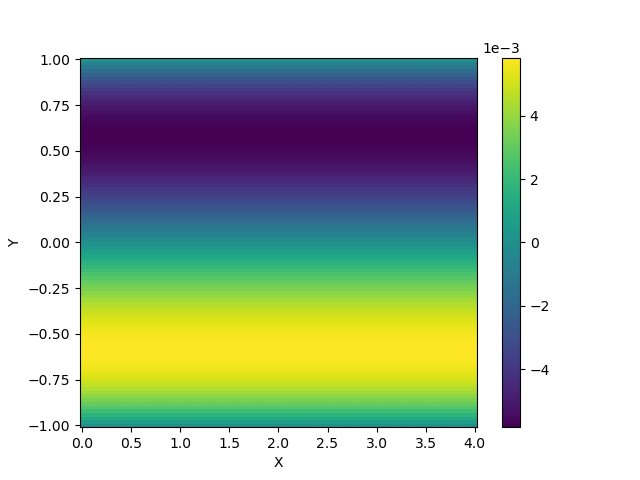}}
    \subfigure[$u_1(y)$]{\includegraphics[width=0.25\hsize]{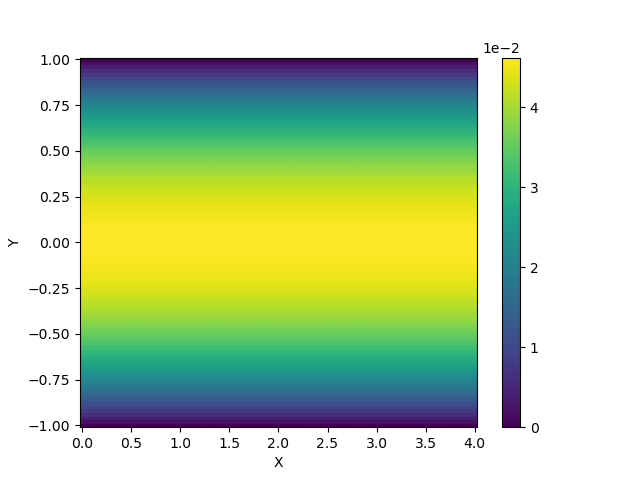}}
    \subfigure[$p$]{\includegraphics[width=0.25\hsize]{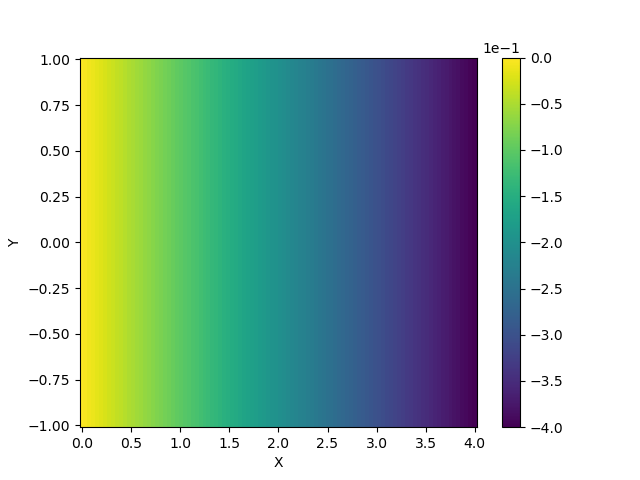}}
    \subfigure[$\hat{B}_1(y)$]{\includegraphics[width=0.25\hsize]{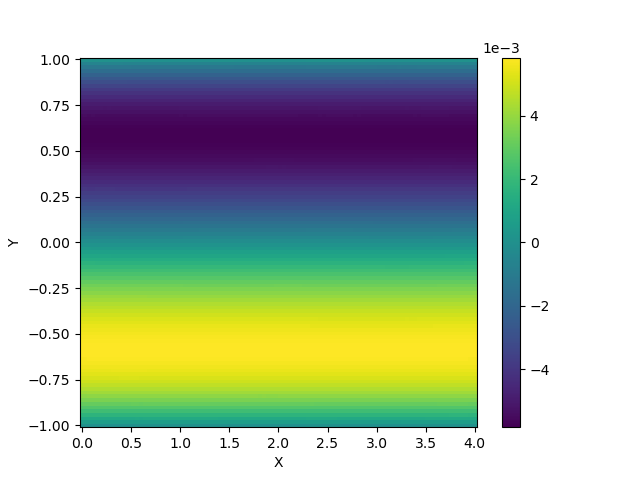}}
    \subfigure[$\hat{u}_1(y)$]{\includegraphics[width=0.25\hsize]{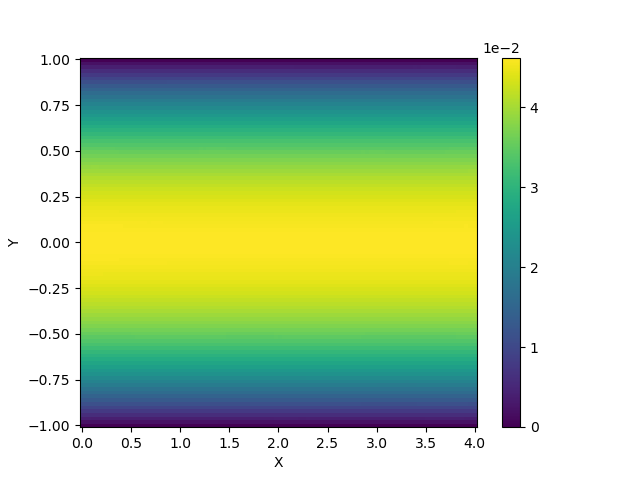}}
    \subfigure[$\hat{p}$]{\includegraphics[width=0.25\hsize]{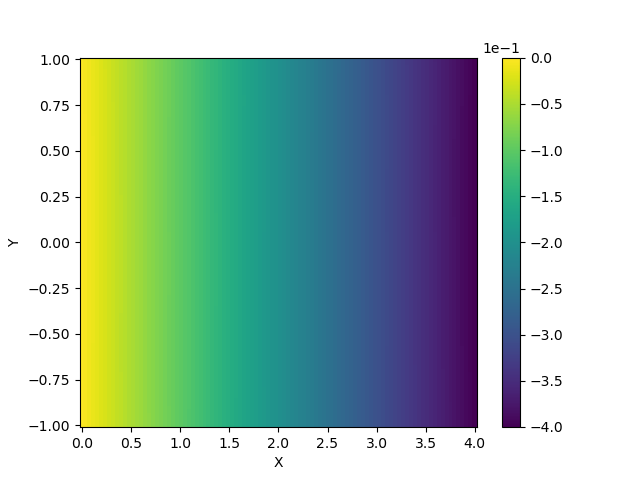}}
    \subfigure[$|\hat{B}_1(y)-B_1(y)|$]{\includegraphics[width=0.25\hsize]{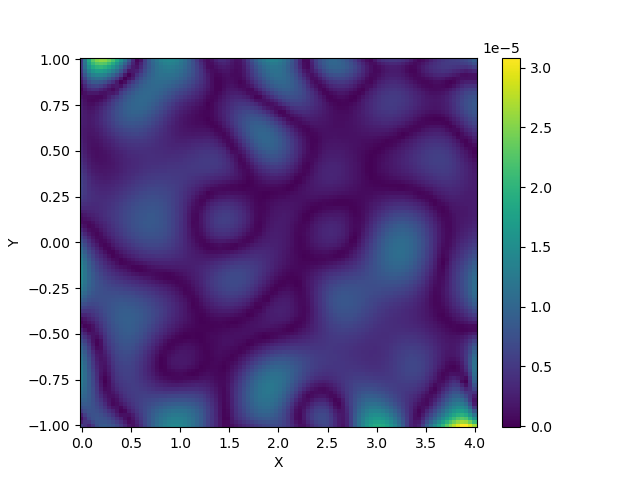}}
    \subfigure[$|\hat{u}_1(y)-u_1(y)|$]{\includegraphics[width=0.25\hsize]{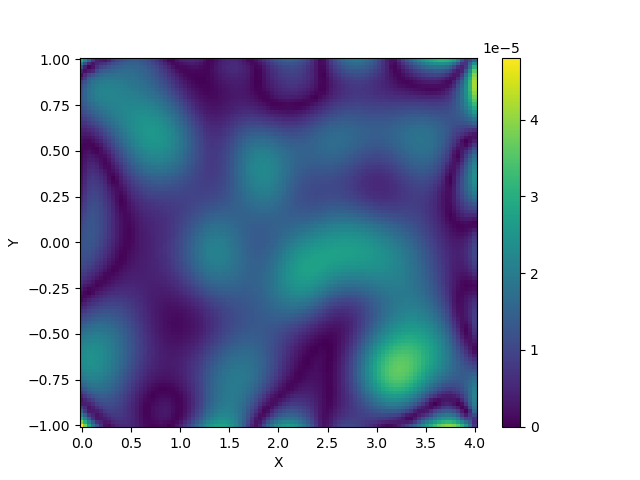}}
    \subfigure[$|\hat{p}-p|$]{\includegraphics[width=0.25\hsize]{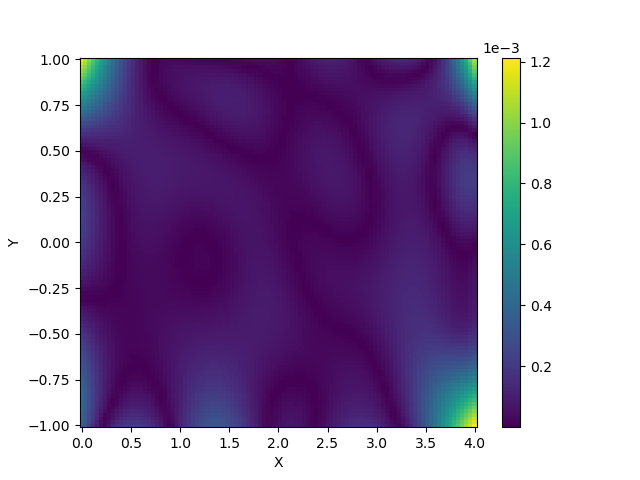}}
    \caption{Hartmann flow:Re=Rm=s=1}
    \label{fig:2d-3-Re=1}
\end{figure}

\begin{figure}[H]
    \centering
    \subfigure[$B_1(y)$]{\includegraphics[width=0.25\hsize]{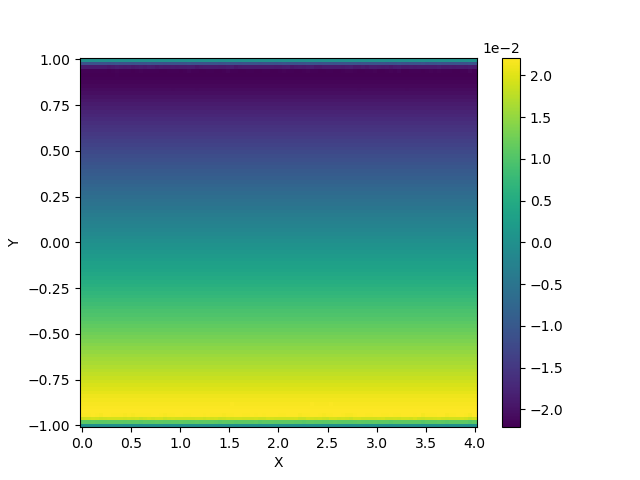}}
    \subfigure[$u_1(y)$]{\includegraphics[width=0.25\hsize]{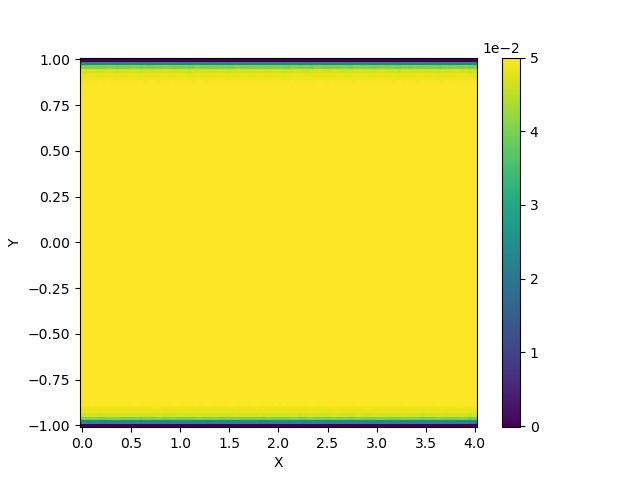}}
    \subfigure[$p$]{\includegraphics[width=0.25\hsize]{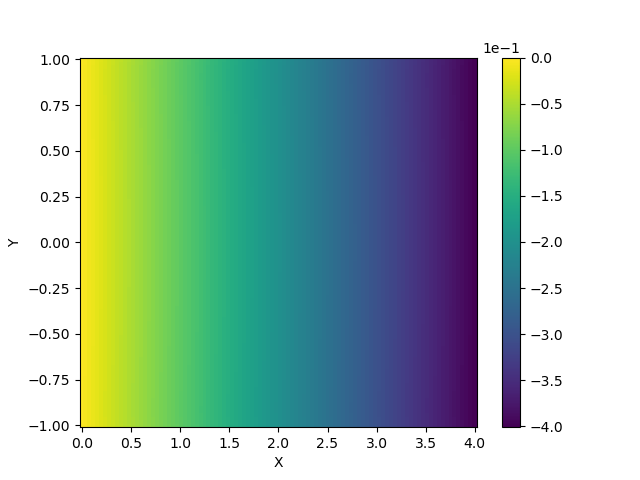}}
    \subfigure[$\hat{B}_1(y)$]{\includegraphics[width=0.25\hsize]{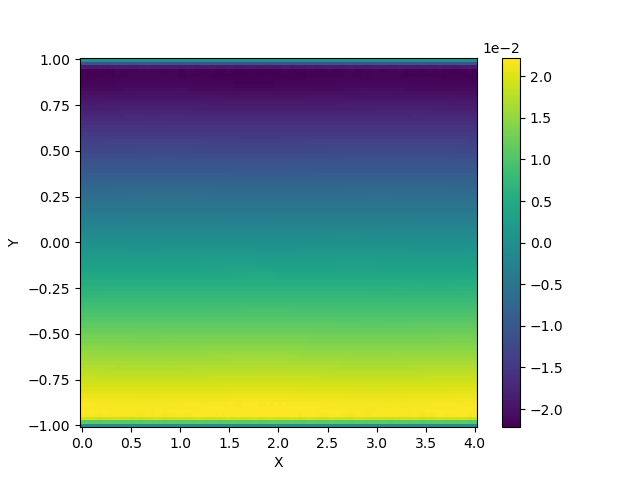}}
    \subfigure[$\hat{u}_1(y)$]{\includegraphics[width=0.25\hsize]{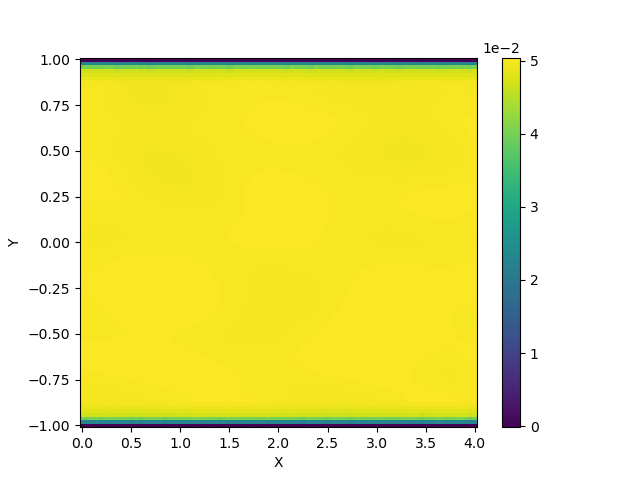}}
    \subfigure[$\hat{p}$]{\includegraphics[width=0.25\hsize]{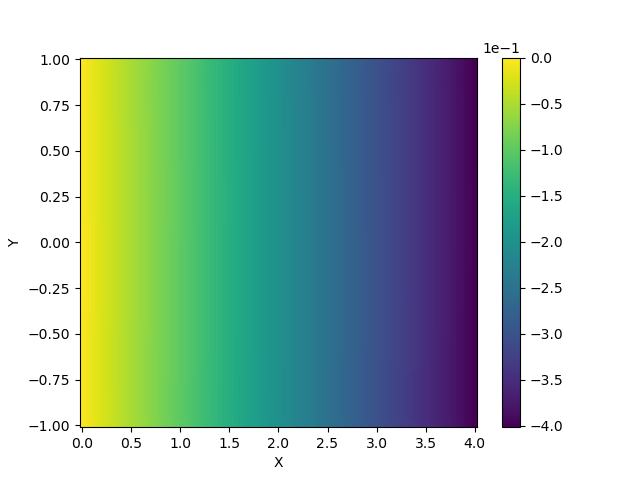}}
    \subfigure[$|\hat{B}_1(y)-B_1(y)|$]{\includegraphics[width=0.25\hsize]{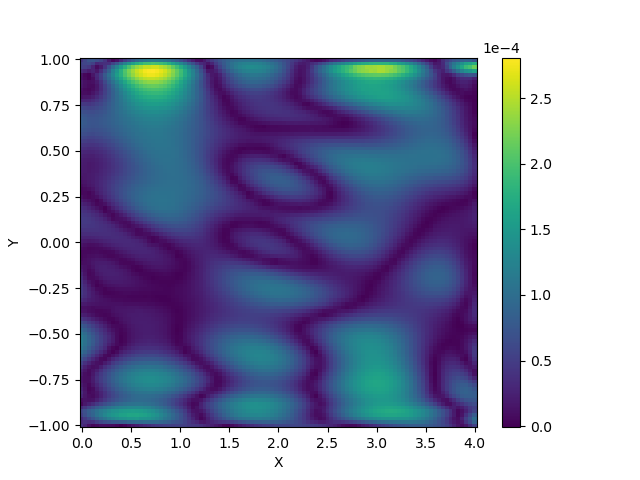}}
    \subfigure[$|\hat{u}_1(y)-u_1(y)|$]{\includegraphics[width=0.25\hsize]{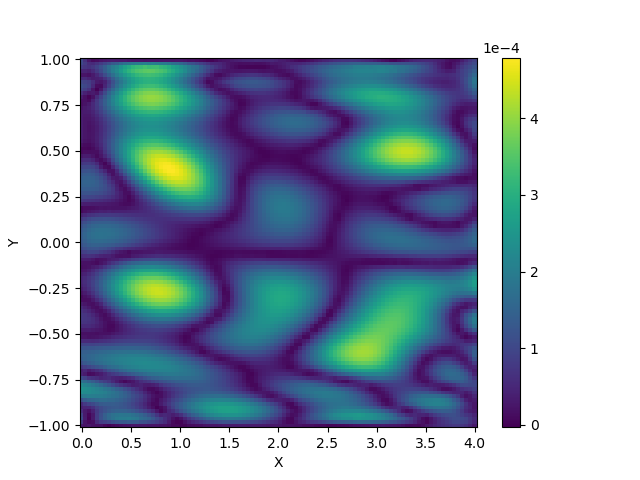}}
    \subfigure[$|\hat{p}-p|$]{\includegraphics[width=0.25\hsize]{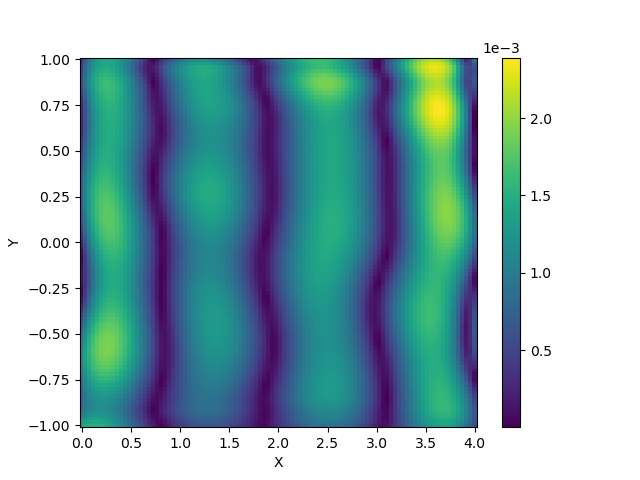}}
    \caption{Hartmann flow:Re=Rm=20,,s=4}
    \label{fig:2d-3-Re=20}
\end{figure}

\begin{figure}[H]
    \centering
    \subfigure[$B_1(y)$]{\includegraphics[width=0.25\hsize]{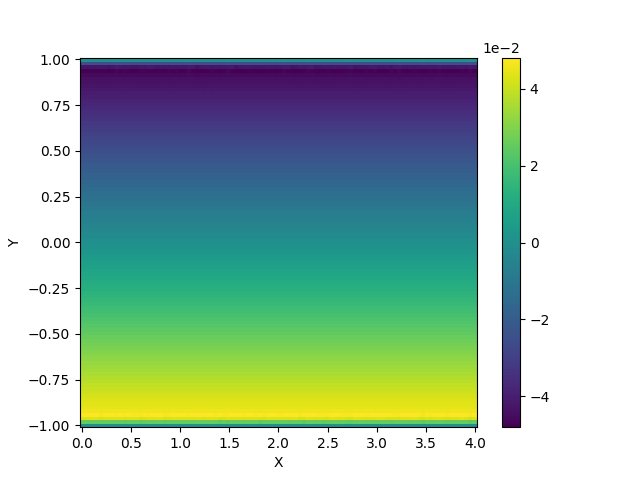}}
    \subfigure[$u_1(y)$]{\includegraphics[width=0.25\hsize]{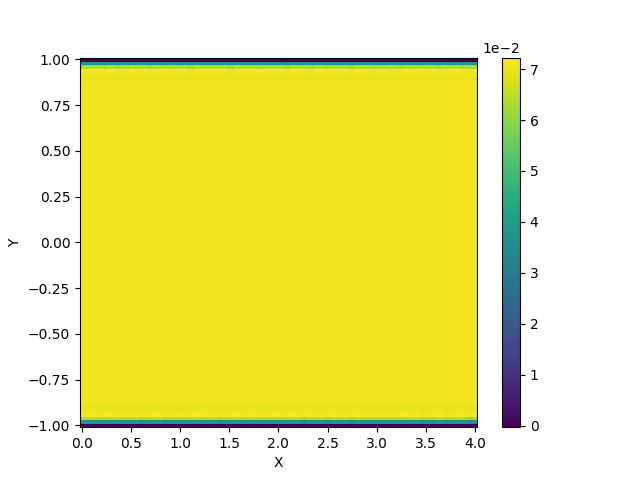}}
    \subfigure[$p$]{\includegraphics[width=0.25\hsize]{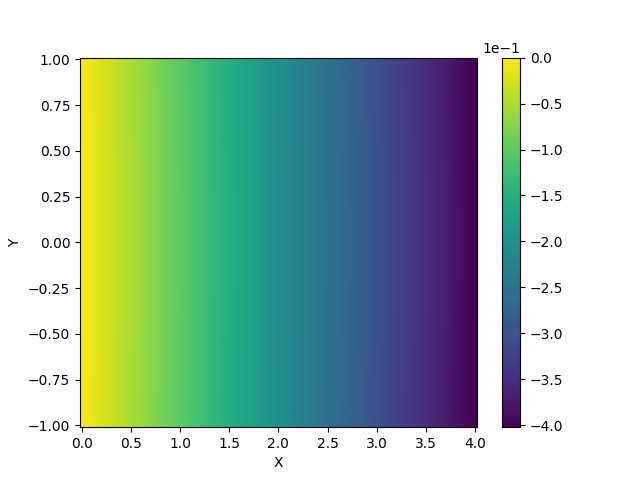}}
    \subfigure[$\hat{B}_1(y)$]{\includegraphics[width=0.25\hsize]{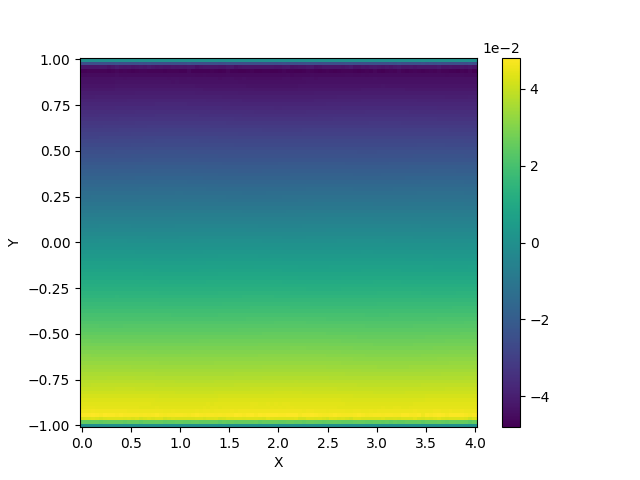}}
    \subfigure[$\hat{u}_1(y)$]{\includegraphics[width=0.25\hsize]{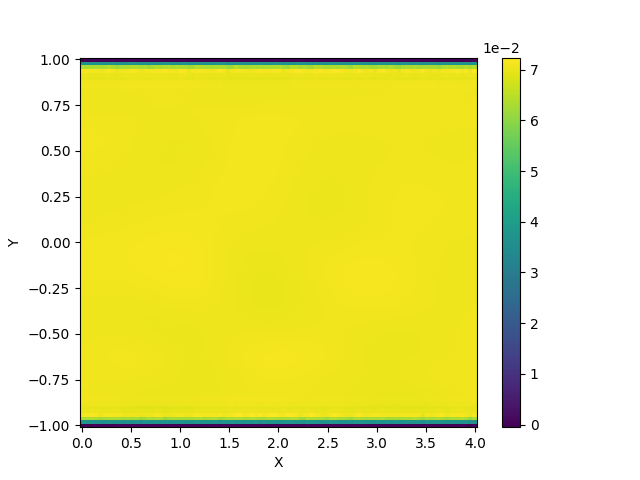}}
    \subfigure[$\hat{p}$]{\includegraphics[width=0.25\hsize]{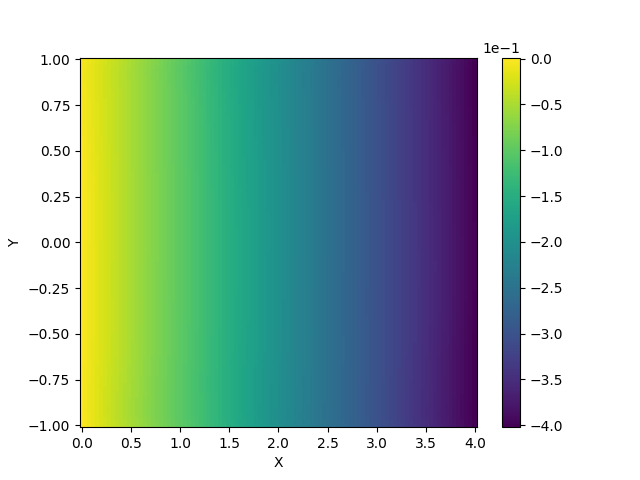}}
    \subfigure[$|\hat{B}_1(y)-B_1(y)|$]{\includegraphics[width=0.25\hsize]{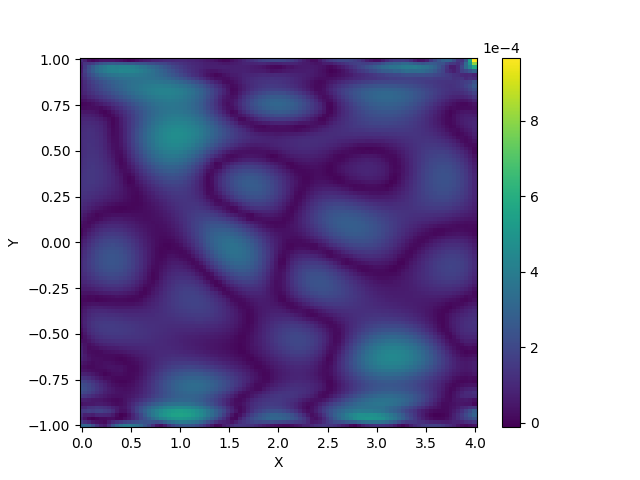}}
    \subfigure[$|\hat{u}_1(y)-u_1(y)|$]{\includegraphics[width=0.25\hsize]{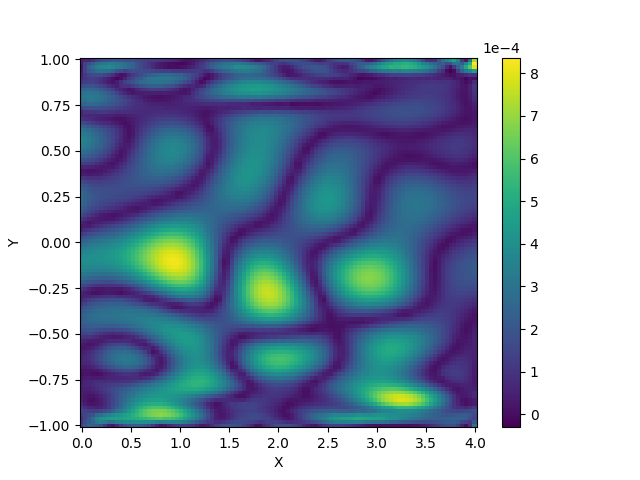}}
    \subfigure[$|\hat{p}-p|$]{\includegraphics[width=0.25\hsize]{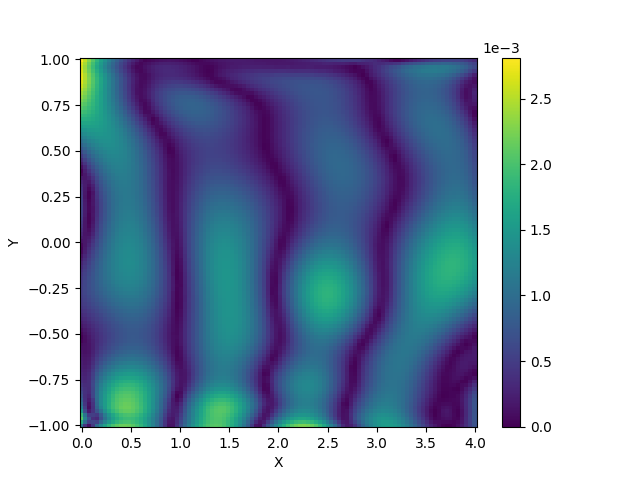}}
    \caption{Hartmann flow:Re=Rm=40,s=2, }
    \label{fig:2d-3-Re=40}
\end{figure}

\begin{figure}[H]
    \centering
    \subfigure[$B_1(y)$]{\includegraphics[width=0.25\hsize]{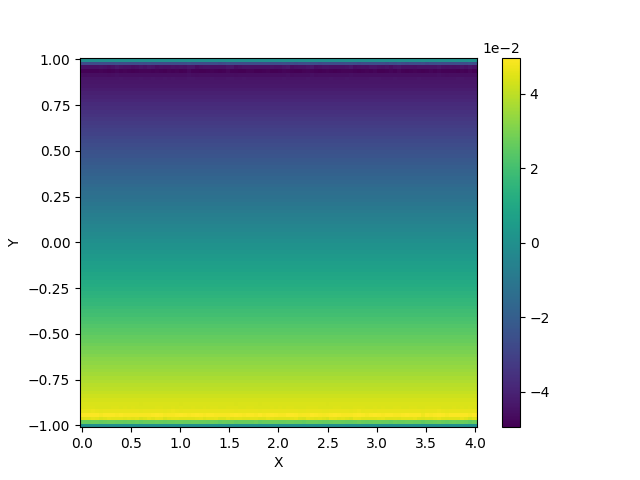}}
    \subfigure[$u_1(y)$]{\includegraphics[width=0.25\hsize]{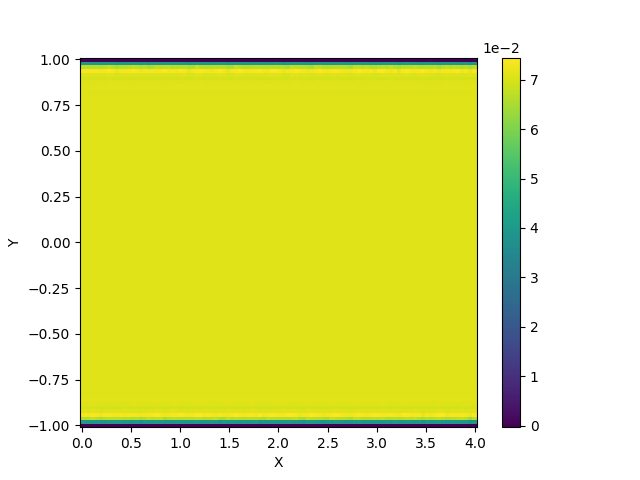}}
    \subfigure[$p$]{\includegraphics[width=0.25\hsize]{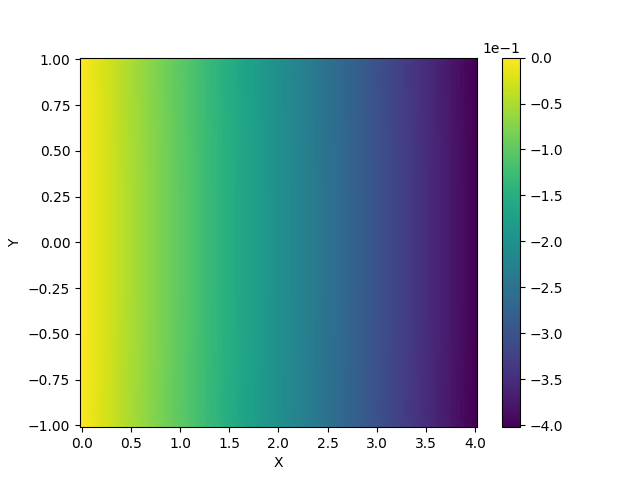}}
    \subfigure[$\hat{B}_1(y)$]{\includegraphics[width=0.25\hsize]{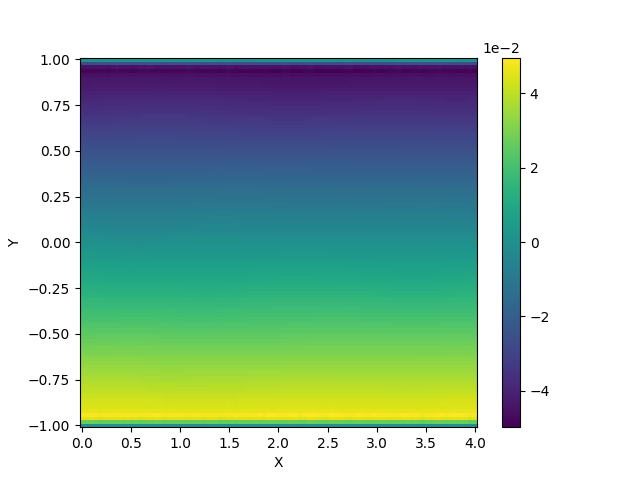}}
    \subfigure[$\hat{u}_1(y)$]{\includegraphics[width=0.25\hsize]{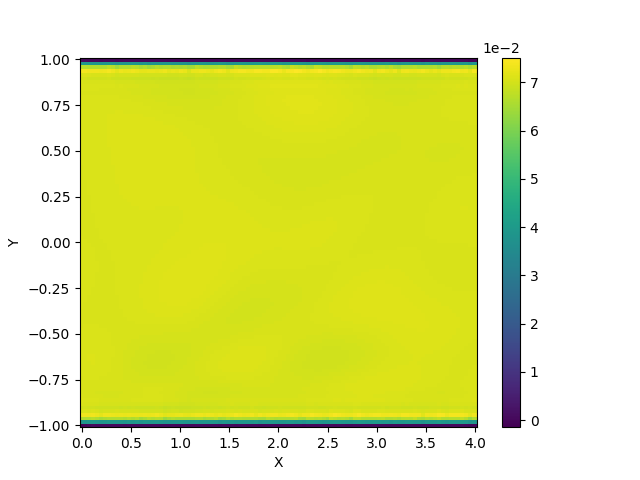}}
    \subfigure[$\hat{p}$]{\includegraphics[width=0.25\hsize]{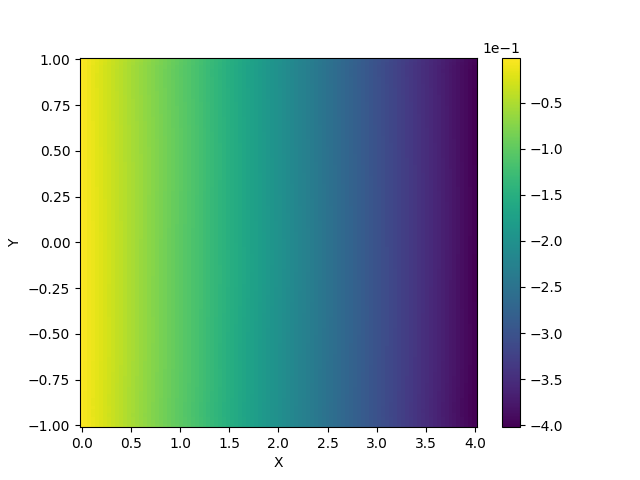}}
    \subfigure[$|\hat{B}_1(y)-B_1(y)|$]{\includegraphics[width=0.25\hsize]{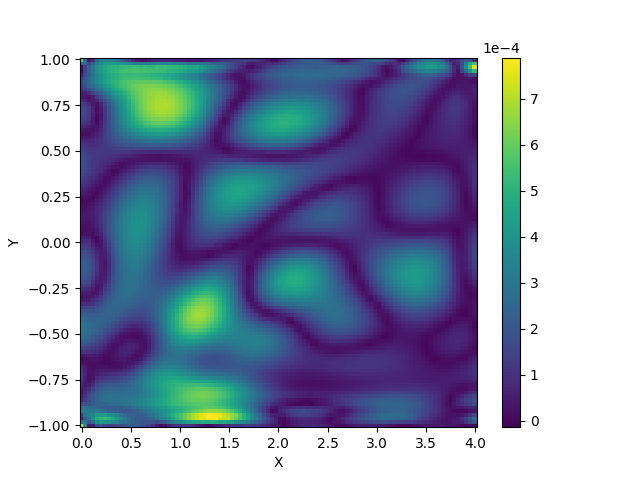}}
    \subfigure[$|\hat{u}_1(y)-u_1(y)|$]{\includegraphics[width=0.25\hsize]{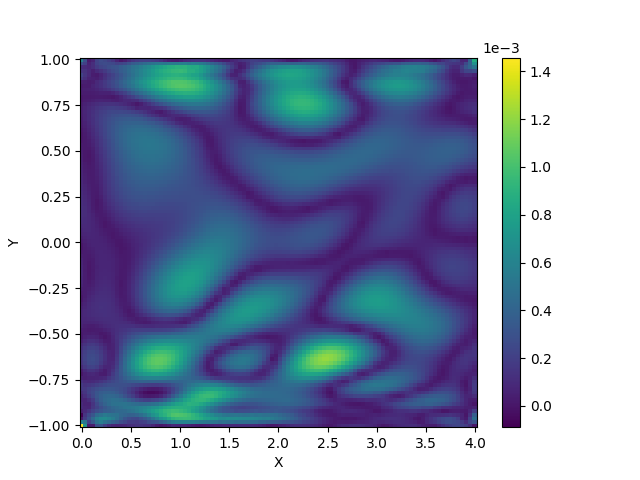}}
    \subfigure[$|\hat{p}-p|$]{\includegraphics[width=0.25\hsize]{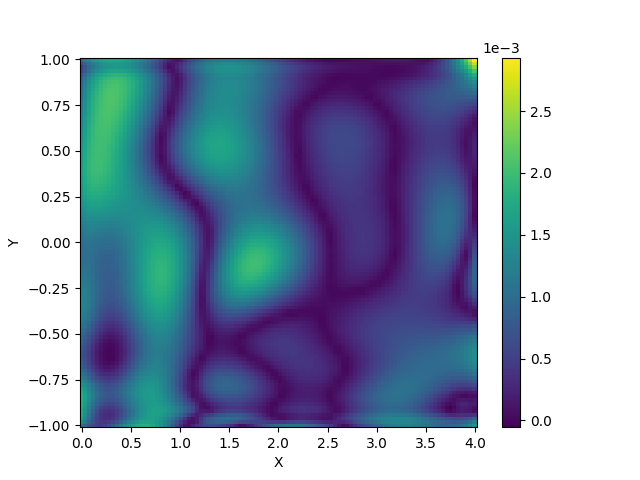}}
    \caption{Hartmann flow:Re=Rm=50,s=2}
    \label{fig:2d-3-Re=50}
\end{figure}
\bibliography{mhdpinn}

\end{CJK}
\end{document}